\documentclass[12pt]{article}
\usepackage{graphicx}
\usepackage{natbib} %comment out if you do not have the package
\usepackage{url} % not crucial - just used below for the URL 
\usepackage{enumitem}
\usepackage[utf8]{inputenc}
\usepackage{amssymb, psfig, epsfig, amsmath, multirow, caption}
\usepackage[colorlinks = true, citecolor = blue, urlcolor = blue]{hyperref}
\usepackage{graphicx}
\usepackage{subfiles}
\usepackage{epstopdf}
\usepackage[ruled,vlined]{algorithm2e}
\usepackage{float}
\usepackage{booktabs,xltabular,float,subfig}
\usepackage[english]{babel}
\usepackage{amsthm}
\newtheorem{example}{Example}%
\newtheorem{remark}{Remark}%
\newtheorem{assumption}{Assumption}
\newtheorem{lemma}{Lemma}
\newtheorem{theorem}{Theorem}

\newtheorem{proposition}[theorem]{Proposition}% 

\DeclareGraphicsExtensions{.jpg,.fig,.eps}

\newcommand{\EE}{\mathbb{E}}
\newcommand{\PP}{\mathbb{P}}

\newcommand{\ind}[1]{{\bf 1}\left(#1\right)}

\newcommand{\bom}{\boldsymbol{\omega}}
\newcommand{\KK}{\mathcal{K}}
\newtheorem{thm}{Theorem}
\newcommand{\bx}{\boldsymbol{x}}
\newcommand{\bH}{\boldsymbol{H}}
\newcommand{\bZ}{\boldsymbol{Z}}
\newcommand{\blind}{0}

\begin{document}

\def\spacingset#1{\renewcommand{\baselinestretch}%
{#1}\small\normalsize} \spacingset{1}

\date{}
%%%%%%%%%%%%%%%%%%%%%%%%%%%%%%%%%%%%%%%%%%%%%%%%%%%%%%%%%%%%%%%%%%%%%%%%%%%%%%

\if0\blind
{
  \title{\bf Online Learning and Optimization for Queues with Unknown Demand Curve and Service Distribution}
  \author{Xinyun Chen\\
  	School of Data Science, School of Management and Economics,\\ Chinese University of Hong Kong (Shenzhen),\\
  		Guiyu Hong\\
  	College of Business, Shanghai University of Finance and Economics,\\
    Yunan Liu\\
    Dept. of Industrial and Systems Engineering, North Carolina State University\\
}
  \maketitle
} \fi

\if1\blind
{
  \bigskip
  \bigskip
  \bigskip
  \begin{center}
    {\LARGE\bf Title}
\end{center}
  \medskip
} \fi

\bigskip
\begin{abstract}
	We investigate an optimization problem in a queueing system where the service provider selects the optimal service fee $p$ and service capacity $\mu$ to maximize the cumulative expected profit (the service revenue minus the capacity cost and delay penalty). 
%Due to the complex nature of the queueing dynamics, such a problem has no analytic solution so that previous research often resorts to heavy-traffic analysis where both the arrival rate and service rate are sent to infinity. In this work we propose an online learning framework designed for solving this problem which does not require the system’s scale to increase. 
%The main challenge is that neither the arrival rate or the service-time distribution is available. 
The conventional \textit{predict-then-optimize} (PTO) approach takes two steps: first, it estimates the model parameters (e.g., arrival rate and service-time distribution) from data; second, it optimizes a model taking these parameters as input. A major drawback of PTO is that its solution accuracy can  often be highly sensitive to the parameter estimation errors because PTO is unable to effectively account for how these errors (step 1) will impact  the solution quality of the downstream optimization (step 2). 
To remedy this issue, we develop an online learning framework that automatically incorporates the
aforementioned parameter estimation errors in the  optimization process; it is an end-to-end
approach that can learn the optimal solution without needing to set up the parameter estimation
as a separate step as in PTO.
%The ingenuity of this approach lies in its online nature, which allows the service provider do better by interacting with the environment. 
Effectiveness of our online learning approach is substantiated by (i) theoretical results including the algorithm convergence and analysis of the \textit{regret} (``cost" to pay over time for the algorithm to learn the optimal policy), and (ii) engineering confirmation via simulation experiments of a variety of representative examples. We also provide careful comparisons between PTO and our online learning method.
\end{abstract}

\noindent%
\vfill

\newpage
\spacingset{1.25} % DON'T change the spacing!
%%%%%%%%%%%%%%%%%%%%%%%%%%%%%%%%%%%%%%%%%%%%%%%%%%%%%%%%%%%%%%%%%%%%%%%%%%%%%%%%%%%%%%%%%%%%%%%%%%%%%%%%%\part{title}

\section{Introduction}\label{sec:intro}
The conventional performance analysis and optimization in queueing systems require the precise knowledge of certain distributional information of the arrival process and service times. For example, consider the $M/GI/1$ queue having Poisson arrivals and  general service times, the expected steady-state workload $W(\lambda,\mu, c_{s}^2)$ is a function of the arrival rate $\lambda$, service rate $\mu$ and second moment or \textit{squared coefficient of variation} (SCV) $c_{s}^2\equiv \text{Var}(S)/\mathbb{E}[S]^2$ of the service time $S$. In particular, according to the famous Pollaczek–Khinchine (PK) formula \citep{pollaczek1930aufgabe}, we have
\begin{align}\label{PK}
	\mathbb{E}[W(\lambda, \mu, c_{s}^2)] = \frac{\rho}{1-\rho}\frac{1+c_s^2}{2}, \qquad \text{with}\quad \rho \equiv \frac{\lambda}{\mu}.
\end{align}
%\red{This seems to be PK for workload? PK for queue shall be $\rho+\frac{\rho^2(1+c_s^2)}{2(1-\rho)}$. For $M/M/1$ they are equal, but not for $M/G/1$. }
%\red{[C: should be $\EE[W_\infty(\lambda, \mu)]$ ?]\\}
One can never overstate the power of the PK formula because it has such a nice structure that insightfully ties the system performance to all model primitives $\lambda$, $\mu$ and $c_s^2$. Indeed, the PK formula has been predominantly used in practice and largely extended to several more general settings such as the $GI/GI/1$ queue with non-Poisson arrivals \citep{abate1993calculation} and $M/GI/n$ queue with multiple servers \citep{cosmetatos1976some}. 

To optimize desired queueing performance, it is natural to follow the \textit{predict-then-optimize}
(PTO) approach, where ``predict'' means the estimation of required model parameters (e.g., $\lambda$, $\mu$ and $c_s^2$) from data (e.g., arrival times and service times) and ``optimize'' means the optimization of certain queueing decisions using formulas such as \eqref{PK} with the predicted parameters treated as the true parameters. See panel (a) in Figure 1 for a flow chart of PTO. 
A potential issue of PTO is that the required queueing formulas can be highly sensitive to the estimation errors of the input parameters (e.g., $\lambda$ and $\mu$), especially when the system's congestion is critical. 
%the premise is that the model inputs $\lambda$, $\mu$ and $c_s^2$ ought to be available and precise. In practice, one always needs to statistically calibrate all input parameters using existing data  of the arrival and service times. Unfortunately, 
%clear-cut structure, despite of its elegance, may become a disadvantage rather than an advantage. 
%   making the its fidelity questionable. 
%Such an %backfiring 
%effect aggravates when the queueing system is in heavy traffic. 
For example, when $c_s=\mu = 1$ and $\lambda = 0.99$, the PK formula \eqref{PK} yields that $\mathbb{E}[W(\lambda, \mu, c_{s}^2)] = 99$. But a $0.5\%$ increase of the demand rate $\lambda$  will yield $\mathbb{E}[W(\lambda, \mu, c_{s}^2)] = 197$, resulting in a $99\%$ relative error in the predicted workload. %Hence, %We call this the \textit{curse of inaccuracy in heavy traffic}.
%In such a situation, this ``backfires'' clear-cut structure, despite of its elegance, may become a disadvantage rather than an advantage. 
%   making the its fidelity questionable. 
%Such an %backfiring 
Consequently, the practical effectiveness of PTO heavily relies on the accuracy of the prediction step to provide near-perfect estimates of the input parameters. Without such precision, solution methods based on these convenient formulas may prove counterproductive or even fail to deliver the desired outcomes.%the convenient structures of queueing formulas (such as \eqref{PK}), despite of their elegance, 
% (becoming a disadvantage rather than an advantage). %which would necessarily require a large number of data samples. 

\begin{figure}[h]
	\centering
	\includegraphics[width=0.95\linewidth]{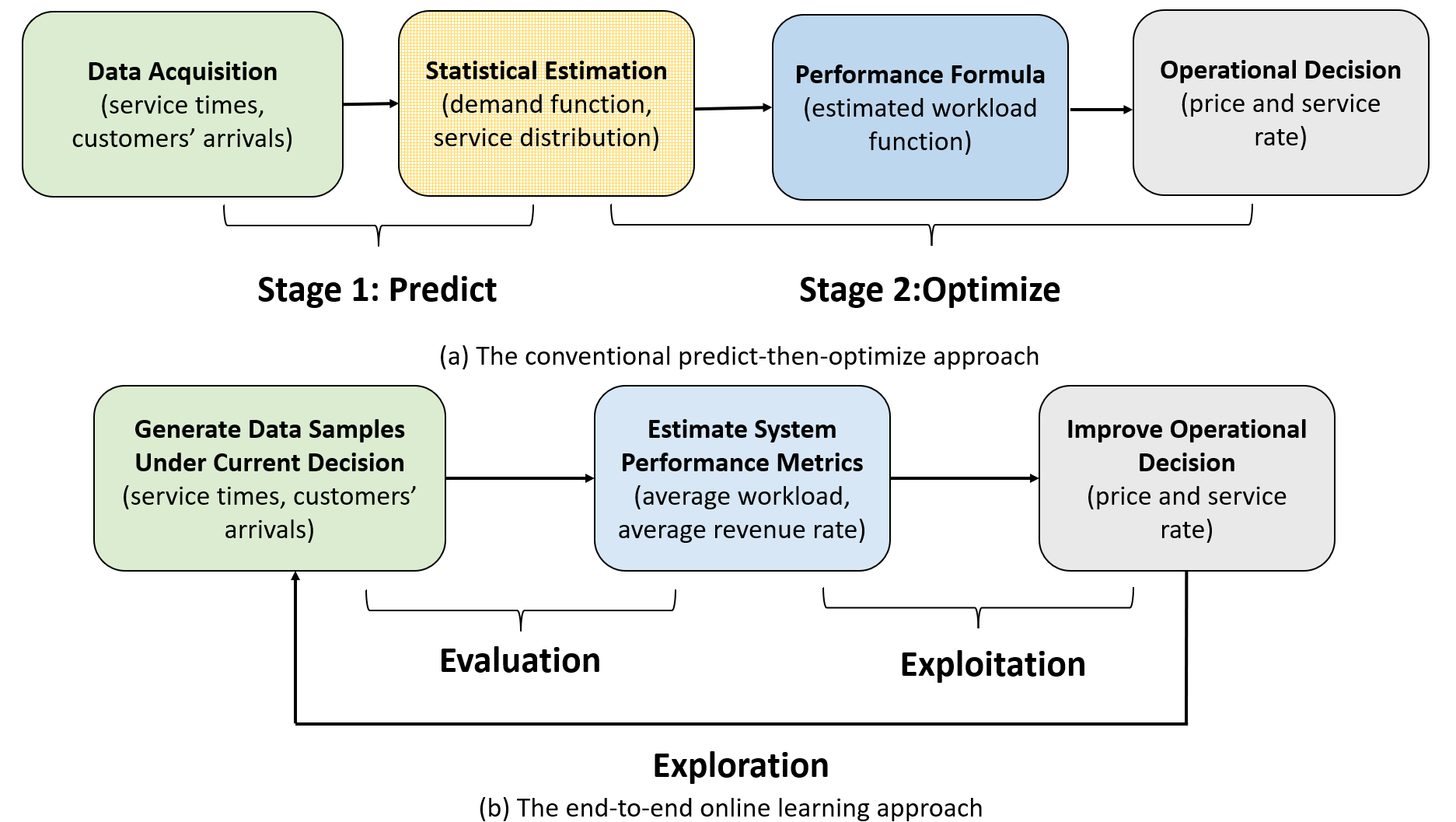}
	\caption{Schematic presentations for (a) the two-step conventional  \textit{predict-then-optimize} scheme and (b) the end-to-end \textit{online learning} scheme. }\label{fig:flow}
\end{figure}

The performance shortcomings of PTO, particularly in heavy-traffic conditions, stem from its inability to adequately account for parameter estimation errors and the substantial impact these errors have on the quality of the resulting ``optimized" decision variables.
To help remedy this issue, %one may of course aim to keep the errors of the input parameters under control. But this may put too much pressure on the statistical fitting 
%that can take into account the impact system parameter estimation errors to decision making. 
we propose \textit{an online learning  framework that automatically incorporates the aforementioned parameter estimation errors in the solution prescription process; it is an end-to-end approach that can learn the optimal solution more directly from data,  so that we no longer need to set up the parameter estimation as a separate stage as in PTO.}
In this paper, we solve a pricing and capacity sizing problem in an $M/GI/1$ queue, where the service provider seeks the optimal service fee $p$ and service rate $\mu$ so as to maximize the long-term profit, which is the revenue minus the staffing cost and the queueing penalty, namely,
\begin{align}\label{objective}
	\max_{\mu,p} \ \mathcal{P}(\mu,p) \equiv \lambda(p)p - h_0\mathbb{E}[W] - c(\mu),
\end{align}
where $W$ is the system's steady-state workload, $c(\mu)$ is the cost for providing service capacity $\mu$ and $h_0$ is a holding cost per job per unit of time. Problems in this framework have a long history, see for example \cite{Kumar2010}, \cite{LeeWard2014}, \cite{LeeWard2019}, \cite{Maglaras2003}, \cite{Nair2016}, \cite{Kim2018}, \cite{ChenLiuHong1}  and the references therein. %Following \cite{LeeWard2014}. %, we consider a $GI/GI/1$ queue
The major distinction is that in the present paper, we assume that neither the arrival rate $\lambda(p)$ (as a function of $p$) or the service-time distribution is explicitly available to the service provider. %, so the PK formula is not directly applicable. 
(As showed In Section \ref{sec:PTOSensitivity}, we will see that the PTO approach for solving Problem \eqref{objective} indeed suffer from unaccountable estimation errors in the model parameters.) % due to the estimation errors of the model parameters.

Our online learning approach operates in successive cycles in which  the service provider's operational decisions are being continuously evolved using newly generated data. Data here include { customers' arrival and service times} under the policy presently in use. See panel (b) in Figure \ref{fig:flow} for an illustration of the online learning approach. 
{In each iteration $k$, the service provider evaluates the current decision $(\mu_k,p_k)$ based on the newly generated data. Then, the decision is updated to $(\mu_{k+1},p_{k+1})$ according to the evaluation result (the exploitation step). In the next iteration, the service provider continues to operate the system under $(\mu_{k+1},p_{k+1})$ to generate more data (the exploration step).} %At the beginning of each cycle $k$, the service provider’s decisions $(p_k,\mu_k)$ will be computed and enforced throughout the cycle. 
%At the heart of the above procedure for updating decision $(p_k,\mu_k)$ is to obtain an effective estimator $H_{k}$ for the gradient of the profit function. 
%In particular, our online algorithm will update $(p_k,\mu_k)$ according to 
%$$(\mu_k,p_k)\leftarrow (\mu_{k-1},p_{k-1})-\eta_{k-1}H_{k-1},$$ 
%where $\eta_k$ is the updating step size. 
We call this algorithm \textit{Learning in Queue with Unknown Arrival Rate} (LiQUAR).

\subsection{{Advantages and challenges.}}
First, the conventional queueing control problem builds heavily on formulas such as \eqref{PK} and requires the precise knowledge of certain distributional information which may not always be readily available. For example, the acquisition of an accurate estimate of the function $\lambda(p)$ across the entire spectrum of the price $p$ is not straightforward and can be both time consuming and costly. In contrast, the online learning approach eliminates the need for such prior information, excelling at ``learning from scratch''. Second, unlike the two-step PTO procedure, the online learning approach is an integrated method that inherently accounts for estimation errors in observed data during the decision-making process. This allows it to utilize data more effectively, leading to improved decisions that are more robust and effective. In contrast to PTO's ``static'' learning, where prediction and optimization are distinctly separate steps, LiQUAR employs a reactive learning approach, characterized by its continuous and dynamic interaction with data.
%more robust with respect to the sensitivity of system performance on the unknown system parameters.   uccessfully ``skips" the statistical estimation step (i.e., first step in the conventional solution procedure). In other words, because this method directly translates data into desired queueing functions which in turn prescribes improved decisions, it should not suffer from the curse of inaccuracy in heavy traffic.

On the other hand, the development of online learning methodologies in queue systems is far from a straightforward extension of their use in other fields, as it must address the unique characteristics of queueing dynamics. 
First, when the control policy is updated at the beginning of a cycle, the previously established near steady-state dynamics are disrupted, and the system enters a transient phase. The dynamics during this period are endogenously influenced by the updated control policy, giving rise to the so-called \textit{regret of nonstationarity}.
%First, as the control policy is updated at the beginning of a cycle, the previously established near steady-state dynamics breaks down and the system undergoes a new transient period  of which the dynamics endogenously depends on the control policy. Such an effect gives rise to the so-called \textit{regret of nonstationarity}. 
%Next, convergence of the decision iterations heavily relies on statistical efficiency in the evaluation step and properties of queueing data. Towards this, the unique features of queueing dynamics brings new challenges. %On the one hand, the recurring transient behavior due to policy updates is the primary source of increased bias.  On
%the other hand, 
%Unlike the standard online leanring settings (e.g., stochastic bandit), queueing data
%such as waiting time and queue length are biased, unbounded and temporally correlated. 
Second, the convergence of decision iterations depends heavily on the statistical efficiency of the evaluation step and the specific properties of the queueing data. This introduces new challenges due to the distinctive nature of queueing dynamics. Unlike standard online learning settings (e.g. stochastic bandits), queueing data such as waiting times and queue lengths are often \textit{biased, unbounded, and temporally correlated.}
These unique features of queueing models present significant obstacles to the design and analysis of online learning methodologies, necessitating novel approaches that account for these complexities.
Finally, our algorithm operates without requiring knowledge of the arrival rate function or the service-time distribution. This makes our research problem more challenging because we cannot take advantage of the detailed structure of the underlying model. Therefore, we are motivated to develop a conceptually simple, model-free learning framework in order to address the above-mentioned challenges. % based on FD gradient estimator. Such model-free design requires more careful control on the bias and variance of gradient estimator. Then, we discuss the feature of queueing data. }

\subsection{Contributions and organization}
Our paper makes the following contributions.
\begin{itemize}
\item We are the first to develop an %model-free 
online learning scheme for the $M/GI/1$ queue with unknown demand function and service-time distribution. The effectiveness of our algorithm stems from its well-integrated  queueing features, encompassing both the overall algorithm design and the optimization of hyperparameters.
For our online learning algorithm, we establish a regret bound of $O(\sqrt{T}\log(T))$. In comparison with the standard $O(\sqrt{T})$ regret for model-free \textit{stochastic gradient descent} (SGD) methods assuming unbiased and independent reward samples, 
%which will yield regret of $O(\sqrt{T})$ if applied online \citep{BoradieZeevi2011}. The 
our regret analysis exhibits an extra $\log(T)$ term which rises from the nonstationary queueing dynamics due to the policy updates. For the  $M/M/1$ model, we derived a more detailed regret bound expressed explicitly as a function of the traffic intensity.

\item At the heart of our regret analysis is to properly link the estimation errors from queueing data to the algorithm's hyperparameters and regret bound. For this purpose, we develop new results that establish useful statistical properties of data samples generated by a $M/G/1$ queue. 
%such as the bound of auto-covariance of queueing data, 
%based on ergodicity analysis.
Besides serving as building blocks for our regret analysis in the present paper, these results are of independent research interest %for statistical inference of queueing data. 
and may be used %We develop a new framework 
to analyze the estimation errors of data in sequential decision making in ergodic  queues. Hence, the theoretic analysis and construction of the gradient estimator may be extended to other queueing models which share similar ergodicity properties. %\red{Methodology contribution:  This framework can be applied to general queueing models with similar ergodicity results. }

%	\item At the heart of our algorithm design is the development and analysis of an efficient gradient estimator which is easily computed from our queueing data. For this purpose, we develop new results that establish useful statistical properties of data samples generated by a $M/G/1$ queue 
%such as the bound of auto-covariance of queueing data, 
%based on ergodicity analysis.
%Besides serving as cornerstones for our regret analysis, these results are of independent research interest for statistical inference of queueing data. 
%We develop a new framework to analyze the estimation errors of queueing data in sequential decision making utilizing ergodicity analysis of queues. Hence, the theoretic analysis and construction of the gradient estimator may be extended to other queueing models which share similar ergodicity properties. %\red{Methodology contribution:  This framework can be applied to general queueing models with similar ergodicity results. }

\item  Supplementing the theoretical results, we evaluate the practical
effectiveness of our method by conducting comprehensive numerical experiments. In particular, our numerical results confirm that the online learning algorithm is efficient and robust to several model and algorithm parameters such as service distributions and updating step sizes; we also generalize our algorithm to the $GI/GI/1$ model. 
Next, we conduct a systematic analysis and experiments to compare LiQUAR to (i) PTO and (ii) gradient-based reinforcement learning methods.                     
\end{itemize}

\paragraph{\textbf{Organization of the paper.}} 
In Section \ref{sec: liter}, we review the related literature. 
In Section \ref{sec: ModelandAssumptions}, we introduce the model and its assumptions. 
In Section \ref{sec: algorithm}, we present LiQUAR and describe how the queueing data is processed in our algorithm. 
{In Section \ref{sec: regret upper bound}, we conduct the convergence and regret analysis for LiQUAR. The key steps of our analysis form a quantitative explanation of how estimation errors in queueing data propagate through our algorithm flow and how they influence the quality of the LiQUAR solutions. We analyze the total regret by separately treating \textit{regret of non-stationarity} - the part of regret stemming from transient system dynamics, \textit{regret of suboptimality} - the part aroused by the errors due to suboptimal decsions, and \textit{regret of {finite difference}} - the part originating from the need of estimation of gradient.}
In Section \ref{sec: heavy traffic}, we report a regret bound explicitly expressed as a function of the traffic intensity for $M/M/1$.
%\ref{sec: gradient estimator analysis}, we provide the detailed analysis for the gradient estimator and confirm its statistical effectiveness by establishing bounds for its (conditional) bias and variance, which are of independent interest for statistical inference of queueing data. 
In Section \ref{sec: Numerical}, we conduct numerical experiments to confirm the effectiveness and robustness of LiQUAR. 
In Sections \ref{sec:PTO_compare} and \ref{sec: RL}, We compare LiQUAR to PTO and gradient-based reinforcement learning methods.
%In Section \ref{sec: proofs} we give the proofs of our main results.
We provide concluding remarks in Section \ref{sec: Conclusions}.
Technical proofs and supplementary results are given in the e-Companion.

\section{Related Literature}\label{sec: liter}
The present paper is related to the following four streams of literature.
\paragraph{Pricing and capacity sizing in queues.} There is rich literature on pricing and capacity sizing for service systems under various settings. \cite{Maglaras2003} studies pricing and capacity sizing problem in a processor sharing queue motivated by internet applications; \cite{Kumar2010} considers a single-server system with nonlinear delay cost; \cite{Nair2016} studies $M/M/1$ and $M/M/k$ systems with network effect among customers; \cite{Kim2018} considers a dynamic pricing problem in a single-server system. The specific problem that we consider here is related to \cite{LeeWard2014}, which considers joint pricing and capacity sizing for $GI/GI/1$ queues with known demand. Later, they further extend their results to the $GI/GI/1+G$ model with customer abandonment in \cite{LeeWard2019}. 
Although the present work is motivated by the pricing and capacity sizing problem for service systems, unlike the above-cited works, we assume no knowledge of the demand rate and service distribution.
%Even though this paper is related to pricing and capacity sizing, the existing results usually consider known demand or conventional predict-then-optimize scheme. Our methodology is very different from these papers.

\paragraph{Demand Learning.} \cite{Broder2012Dynamic} considers a dynamic pricing problem for a single product with an unknown parametric demand curve and establishes an optimal minimax regret in the order of $O(\sqrt{T})$. \cite{keskin2014dynamic} investigates a pricing problem for a set of products with an unknown parameter of the underlying demand curve. \cite{BesbesZeevi2015} studies demand learning using a linear curve as a local approximation of the demand curve and establishes a minimax regret in the order of $O(\sqrt{T})$. Later, \cite{cheung2017dynamic} solves a dynamic pricing and demand learning problem with limited price experiments. We draw distinctions from these papers by studying a pricing and capacity sizing problem with demand learning in a queueing setting where our algorithm design and analysis need to take into account unique features of the queueing systems.

\paragraph{Machine learning in queueing systems} Our paper is related to a small but booming literature on machine earning in queueing systems. \cite{Dai2021} studies an actor-critic algorithm for queueing networks. \cite{Xie2019} and \cite{Xie2020} develop reinforcement learning techniques to treat the unboundedness of the state space of queueing systems.   \cite{Krishnasamy2021learning} develops bandit methods for scheduling problems in a multi-server queue with unknown service rates. \cite{zhong2022learning} proposes an online learning method to study a scheduling problem for a multiclass $M_t/M/N+M$ system with unknown service rates and abandonment rates. \cite{ChenLiuHong1} studies the joint pricing and capacity sizing problem for $GI/GI/1$ with known demand. 
%\cite{Shi2020}, \cite{jia2022online} \red{cite} study a price-based revenue management problem in an $M/M/c$ queue with discrete/continuous price space using the bandit framework.
See \cite{walton2021learning} for a review of the role of information and learning in queueing systems. Recent research has explored the application of deep learning methods to predict queueing performance: \cite{BKSS23} proposes a deep-learning-based steady-state predictor for the $GI/GI/1$ queue; \cite{GLBC24a,GLBC24b} develop recurrent neural network models to predict transient performance in nonstationary queues.
Our paper is most closely related to \cite{Shi2020} which studies a price-based revenue management problem in an $M/M/c$ queue with unknown demand and discrete price space, under a multi-armed bandit framework.
Later, \cite{jia2022online} extends the results in \cite{Shi2020} to the problem setting with a continuous price space and {considers linear demand functions}. %mappings between the rates and the price with} linear demand.  %Based on this linear structure, they adapt the algorithms in \cite{Shi2020} into their linear bandit version and show the efficiency of the algorithms by proving the regret bounds.
Similar to \cite{Shi2020,jia2022online}, we also study a queueing control problem with unknown demand {and continuous decision variables}. The major distinction is that in addition to maximizing the service profit as by \cite{Shi2020}, the present paper also includes a \textit{queueing penalty} in our optimization problem as a measurement of the quality  of service
%. On the one hand, this is a practical consideration because a holding cost of the queue length or waiting is widely used in the queueing theory literature, e.g., 
\citep{Kumar2010,LeeWard2014,LeeWard2019}. %, to measure the quality of service. On the other hand,
However, this introduces new technical challenges in algorithm design and regret analysis, such as addressing the bias and autocorrelation inherent in queueing data. Besides, the present paper considers more general service distributions and demand functions. %, which also brings new challenges to the analysis. To deal with the challenges aroused by these difference, 
{ %\cite{huh2009nonparametric} investigates a stochastic inventory planning problem with zero lead time and unknown demand distribution.
Online learning with continuous decision-making has also been explored in inventory systems. For instance, \cite{huh2009adaptive} proposed an SGD-based algorithm to optimize base-stock policies for inventory systems with positive lead times. Later, \cite{zhang2020closing} developed a simulation-based algorithm for the same problem, achieving an optimal regret bound of $O(T^{1/2})$. More recently, \cite{Cong2021Marrying} integrated stochastic gradient descent with bandit algorithms to address convexity challenges and optimize $(s,S)$ policies.
While our approach also employs gradient-based methods, the fundamental differences between queueing system dynamics and inventory models lead to distinct algorithm designs, particularly in the construction of gradient estimators, as well as in the theoretical analysis.
}

\section{Model and Assumptions}\label{sec: ModelandAssumptions}
We study an $M/GI/1$ queueing system having customer arrivals according to a Poisson process (the $M$), \textit{independent and identically distributed} (I.I.D.) service times following a general distribution (the $GI$), and a single server that provides service following  the \textit{first-in-first-out} (FIFO) discipline.
Each customer upon joining the queue is charged by the service provider a fee $p>0$. The demand arrival rate (per time unit) depends on the service fee $p$ and is denoted as $\lambda(p)$. To maintain a service rate $\mu$, the service provider continuously incurs a staffing cost at a rate $c(\mu)$ per time unit.

For $\mu\in[\underline{\mu},\bar{\mu}]$ and $p\in[\underline{p}, \bar{p}]$, {we have $\lambda(p)\in[\underline{\lambda},\bar{\lambda}]\equiv[\lambda(\bar{p}),\lambda(\underline{p})]$}, and the service provider's goal is to determine the optimal service fee $p^*$ and service capacity $\mu^*$ with the objective of  maximizing the steady-state expected profit,
or equivalently minimizing the objective function $f(\mu,p)$ as follows
\begin{equation}\label{eqObjmin}
\min_{(\mu,p)\in\mathcal{B}} f(\mu, p) \equiv h_0 \EE[W_{\infty}(\mu,p)] + c(\mu)- p\lambda(p),\qquad \mathcal{B}\equiv [\underline{\mu},\bar{\mu}]\times[\underline{p}, \bar{p}].
\end{equation}
Here $W_\infty(\mu,p)$ is the stationary workload process observed in continuous time under control parameter $(\mu,p)$. In detail,  under control parameter $(\mu,p)$, customers arrive according to a Poisson process with rate $\lambda(p)$. Let $V_n$ be an I.I.D. sequence corresponding to customers' \textit{workloads} {under unit service rate} ({under service rate $\mu$,  customer $n$ has service time $V_n/\mu$}). We have $\EE[V_n]=1$ so that the mean service time is $1/\mu$ under service rate $\mu$. Denote by $N(t)$ the number of arrivals by time $t$. The total amount of workload brought by customers at time $t$ is denoted by $J(t)=\sum_{k=1}^{N(t)}V_{k}$.
Then the  workload process $W(t)$ follows the \textit{stochastic differential equation} (SDE)
$$dW(t)=dJ(t)-\mu\ind{W(t)>0}dt.$$ 
In particular, given the initial value of $W(0)$, we have
$$W(t) = R(t) - 0\wedge \min_{0\leq s\leq t}R(s), \quad R(t)\equiv W(0) + J(t)-\mu t.$$ The difference $(W(t) - R(t))/\mu$
is the total idle time of the server by time $t$. It is known in the literature \cite[Corollary 3.3, Chapter X]{AsmussenBook} that under the stability condition $
\lambda(p)<\mu$, the workload process $W(t)$ has a unique stationary distribution and we denote by $W_\infty(\mu,p)$ the stationary workload under parameter $(\mu,p)$.

We impose the following assumptions on the $M/GI/1$ system throughout the paper.
\begin{assumption}\label{assmpt: uniform}\textbf{$($Demand rate, staffing cost, and uniform stability$)$}
\begin{enumerate}
\item[$(a)$] The arrival rate $\lambda(p)$ is continuously differentiable in the third order and non-increasing in $p$. Besides,
\begin{align*}%\label{eq: lambda prime}
	C_1 < \lambda'(p) < C_2,
\end{align*}
where
{\footnotesize
	\begin{align*}
		C_1 &\equiv 2\max\left(g(\bar{\mu})\frac{\lambda''(p)}{\lambda'(p)},g(\underline{\mu})\frac{\lambda''(p)}{\lambda'(p)}\right)\lambda(p)-\frac{4\lambda(p)(\underline{\mu}-\lambda(p))}{h_0C},\quad
		C_2 	\equiv - \max\left(\sqrt{\frac{0\vee \left(-\lambda''(p)(\bar{\mu}-\lambda(p)\right))}{2}}~,~ \frac{p\lambda''(p)}{2}\right),
\end{align*}}
$g(\mu)=\frac{{\mu}}{{\mu}-\lambda(p)}-\frac{p({\mu}-\lambda(p))}{h_0C}$ and $C=(1+c_s^2)/2$. 
\item[$(b)$] The staffing cost $c(\mu)$ is continuously differentiable in the third order,  non-decreasing and convex in $\mu$.
\item[$(c)$] The lower bounds $\underline{p}$ and $\underline{\mu}$ satisfy that  $\lambda(\underline{p})<\underline{\mu}$ so that the system is uniformly stable for all feasible choices of  $(\mu,p)$.
\end{enumerate}
\end{assumption}
Although Condition (a) looks complicated, it essentially requires that the derivative of $\lambda(p)$ be not too large or too small.	Condition (a) will be used to ensure that the objective function $f(\mu,p)$ is convex  in the convergence analysis of our gradient-based online learning algorithm in Section \ref{subsec: convergence to optimality}. The two inequalities hold for a variety of commonly used demand functions, including both convex functions and concave functions. Examples include (1) linear demand $\lambda(p)=a-bp$ with $0<b<4\underline{\lambda}(\underline{\mu}-\bar{\lambda})/h_0C$; (2) quadratic demand $\lambda(p)=c-ap^2$ with $a,c>0$, and $\frac{\bar{\mu}-c}{3\underline{p}^2}<a<\left(\frac{3(\underline{\mu}-\bar{\lambda})\underline{p}}{h_0C}-\frac{\underline{\mu}}{\underline{\mu}-\bar{\lambda}}\right)\frac{\underline{\lambda}}{\bar{p}^2}$; (3) exponential demand $\lambda(p)= \exp(a-bp)$ with $0<b<2/\bar{p}$; (4) logit demand $\lambda(p)=M_0\exp(a-bp)/(1+\exp(a-bp))$ with $a-b\bar{p}<\log(1/2)$ and $0<b<2/\bar{p}$. See Section \ref{appx: demand function} for detailed discussions.   

Condition $(c)$ of Assumption \ref{assmpt: uniform} 
is commonly used in the literature of SGD methods for queueing models to ensure that the steady-state mean waiting time $\EE[W_\infty(\mu,p)]$ is differentiable with respect to model parameters. See \cite{Chong1993}, \cite{Fu1990}, \cite{LEcuyer94b}, \cite{LEcuyer1994}, and also Theorem 3.2 of \cite{Glasserman_1992}.
In Section \ref{sec: LiQUAR-b}, we present an initial attempt to relax Assumption \ref{assmpt: uniform}(c).

We do not require full knowledge of service and inter-arrival time distributions. But in order to bound  the estimation error of the queueing data, we require the individual workload to be light-tailed. Specifically, 
we make the following assumptions on $V_n$.
\begin{assumption}\label{assmpt: light tail} \textbf{$($Light-tailed individual workload$)$}
%There exists a constant $\eta>0$ such that the moment-generating functions
%	\begin{align*}
%		\EE[\exp(\eta V_n)] <\infty \quad \text{and}\quad \EE[\exp(\eta U_n)]<\infty.
%	\end{align*}
There exists a sufficiently small constant $\eta>0$ such that
$$\EE[\exp(\eta V_n)]<\infty.$$ 
In addition, there exist constants $0<\theta <\eta/2\bar{\mu}$ and $\gamma_0 >0$ such that
\begin{equation}\label{eq: gamma0}
\phi_V(\theta) < \log\left(1+\underline{\mu}\theta/\bar{\lambda}\right)-\gamma_0,
\end{equation}
where $\phi_V(\theta)\equiv \log\EE[\exp(\theta V_n)]$ is the cumulant generating functions of $V_n$.
\end{assumption}

Note that $\phi'_V(0)=1$ as $\EE[V_n]=1$. 
Suppose  $\phi_V$ is smooth around 0, then we have  $\phi_V(\theta) = \theta +o(\theta)$ by Taylor's expansion. On the other hand, as $\underline{\mu}>\bar{\lambda}$ under Assumption \ref{assmpt: uniform}, there exists $a>0$ such that $\log\left(1+\underline{\mu}\theta/\bar{\lambda}\right) = (1+a)\theta+o(\theta)$. This implies that, we can choose $\theta$ small enough such that $\log\left(1+\underline{\mu}\theta/\bar{\lambda}\right) -\phi_V(\theta)>\frac{a\theta}{2}$ and then we set $\gamma_0 = \frac{a\theta}{2}$.
Hence, a sufficient condition that warrants \eqref{eq: gamma0} is to require that $\phi_V$ be smooth around 0, which is true for many distributions of  $V$ considered in common queueing models. Assumption \ref{assmpt: light tail} will be used in our proofs to establish ergodicity result.

\section{Our Online Learning Algorithm}\label{sec: algorithm}
We first explain the main ideas in the design of LiQUAR and provide the algorithm outline in Section \ref{subsec: algorithm outline}. The key step in our algorithm design is to construct a data-based gradient estimator, which is explained with details in Section \ref{subsec: system dynamics}. As a unique feature of service systems,  there is a delay in data observation of individual workloads, i.e., they are revealed only after service completion. We also explain how to deal with this issue in Section \ref{subsec: system dynamics}. The design of algorithm hyperparameters in LiQUAR will be specified later in Section \ref{sec: regret upper bound} based on the regret analysis results. In the rest of the paper, we use bold symbols for vectors and matrices.

\subsection{Algorithm outline}\label{subsec: algorithm outline}
The basic structure of LiQUAR follows the online learning scheme as  illustrated in Figure \ref{fig:flow}. It interacts with the queueing system in continuous time and improves pricing and staffing policies iteratively. In each iteration $k\in\{1,2,...\}$, LiQUAR operates the queueing system according to control parameters $\bar{\bx}_k\equiv (\bar{\mu}_k,\bar{p}_k)$ for a certain time period, and collects data generated by the queueing system during the period.
%As the queueing system evolves in continuous time, LiQUAR organizesthe time horizon into successive operational cycles and updates pricing and staffing policies at the beginning of each cycle using data collected in previous cycles. In particular, the time horizon is divided into different cycles of $T_k$ units of time for $k=1,2,...$ and $T_k$ is called the \text{cycle length}. To mitigate the bias due to the transient performance of the queueingprocess, LiQUAR sends the cycle length $T_k$ increasing to infinity  as $k\to\infty$. 
At the end of an iteration, LiQUAR estimates the gradient of the objective function $\nabla f(\bar{\bx}_k)$ based on the collected data and accordingly updates the control parameters. The updated control parameters will be used in the next iteration.

We use the \textit{finite difference} (FD) method \citep{BoradieZeevi2011} to construct our gradient estimator. Our main purpose is to make LiQUAR model-free and applicable to the settings where the demand function $\lambda(p)$ is unknown.  %Next, we shall explain the gradient estimator used in LiQUAR which can be computed from the observed queueing data including customers arrival and service times. 
%In detail, for each iteration $k$,  denote the current control parameter as $\bar{x}_k\equiv (\bar{\mu}_k,\bar{p}_k)$.   % Let   The finite-difference approximation of the gradient $\nabla f(\bar{x}_k)$ yields, 
%$$\nabla f(\bar{x}_k)\approx \frac{~\EE[f(\bar{x}_k+\delta_k\cdot Z_k/2) - f(\bar{x}_k-\delta_k\cdot Z_k/2)]~}{\delta_k}.$$
%To explain the FD gradient estimator, we index the iteration by $k$ and denote the control parameter in iteration $k$ as $\bar{x}_k\equiv (\bar{\mu}_k,\bar{p}_k)$.  
To obtain the FD estimator of $\nabla f(\bar{\bx}_k)$, LiQUAR splits total time of iteration $k$ into two equally divided intervals (i.e., cycles) each with $T_k$ time units. We index the two cycles by $2k-1$ and $2k$, in which the system is respectively operated under control parameters
\begin{equation}\label{eq: cycle policies}
\bx_{2k-1} \equiv \bar{\bx}_k-\delta_k\cdot \bZ_k/2\equiv (\mu_{2k-1},p_{2k-1})\quad \text{ and }\quad \bx_{2k} \equiv \bar{\bx}_k+\delta_k\cdot \bZ_k/2\equiv(\mu_{2k},p_{2k}), 
\end{equation}
where $\delta_k$ is a positive and small number and $\bZ_k\in\mathbb{R}^2$ is a random vector independent of system dynamics such that $\EE[\bZ_k]=(1,1)^\top$. 
Using data collected in the two cycles, LiQUAR obtains  estimates of the system performance $\hat{f}(\bx_{2k})$ and $\hat{f}(\bx_{2k-1})$, which in turn yield the FD approximation for the gradient  $\nabla f(\bar{\bx}_k)$: 
$$\bH_k\equiv\frac{\hat{f}(\bx_{2k})-\hat{f}(\bx_{2k-1})}{\delta_k}.$$ Then, LiQUAR updates the control parameter according to a SGD recursion as $\bar{\bx}_{k+1}=\Pi_\mathcal{B}(\bar{\bx}_k-\eta_k\bH_k)$, where $\Pi_\mathcal{B}$ is the operator that projects $\bar{\bx}_k-\eta_k\bH_k$ to $\mathcal{B}$.
We give the outline of LiQUAR below.

\noindent\textbf{Outline of LiQUAR:}
\begin{enumerate}
\item[0.] Input: hyper-parameters $\{T_k, \eta_k, \delta_k\}$ for $k=1,2,...$, initial policy $\bar{\bx}_1=(\bar{\mu}_1,\bar{p}_1)$. \\
For $k=1,2,...,L$,
\item Obtain $\bx_l$ according to \eqref{eq: cycle policies} for $l=2k-1$ and $2k$. In cycle $l$, operate the system with policy $\bx_l$ for $T_k$ units of time.
\item Compute $\hat{f}(\bx_{2k-1})$ and $\hat{f}(\bx_{2k})$ from the queueing data to build an estimator $\bH_k$ for $\nabla f(\mu_k,p_k)$.
\item Update $\bar{\bx}_{k+1}=\Pi_\mathcal{B}(\bar{\bx}_k-\eta_k \bH_k)$.
\end{enumerate}

Next, we explain in details how the gradient estimator $\bH_k$, along with $\hat{f}(\bx_{2k-1})$ and $\hat{f}(\bx_{2k})$, are computed from the queueing data in Step 2.

\subsection{Computing Gradient Estimator from Queueing Data}\label{subsec: system dynamics}

We first introduce some notation to describe the system dynamics under LiQUAR and the queueing data generated by LiQUAR. For $l\in\{2k-1, 2k\}$, let $W_l(t)$ be the present workload at time $t\in[0,T_k]$ in cycle $l$.  By definition, we have $W_{l+1}(0) = W_{l}(T_k)$ for all $l\geq 1$.  We assume that the system starts empty, i.e., $W_1(0)=0$. At the beginning of each cycle $l$, the control parameter is updated to $(\mu_l,p_l)$. The customers arrive in cycle $l$ according to a Poisson process $N_l(t)$ with rate $\lambda(p_l)$, $0\leq t\leq T_k$. 
%Let's denote by $\{N_l(t):0\leq t\leq T_k\}$ the arrival process of customers in cycle $l$. 
Let $\{V_i^l: i=1,2,...,N_l\}$ be a sequence of I.I.D. random variables denoting customers' individual workloads, where $N_l=N_l(T_k)$ is the total number of customer arrival in cycle $l$. Then, the dynamics of the workload process $W_l(t)$ is described by the SDE:
\begin{equation}\label{eq: W SDE}
W_l(t) = W_l(0) +\sum_{i=1}^{N_l(t)}V^l_i -\mu_l\int_0^t {\bf 1}(W_l(s)>0)ds.
\end{equation}
If the system dynamics is available continuously in time (i.e. $W_l(t)$ was known for all $t\in [0, T_k]$ and $l=2k-1, 2k$), then a natural estimator for $f(\mu_l,p_l)$ would be
\begin{align*}
\hat{f}(\mu_l,p_l)&=\frac{-pN_l}{T_{k}}+\frac{h_0}{T_{k}}\int_{0}^{T_{k}}W_l(t)dt+c(\mu_l).
\end{align*}

\begin{figure}
\centering
\includegraphics[width=1.02\linewidth]{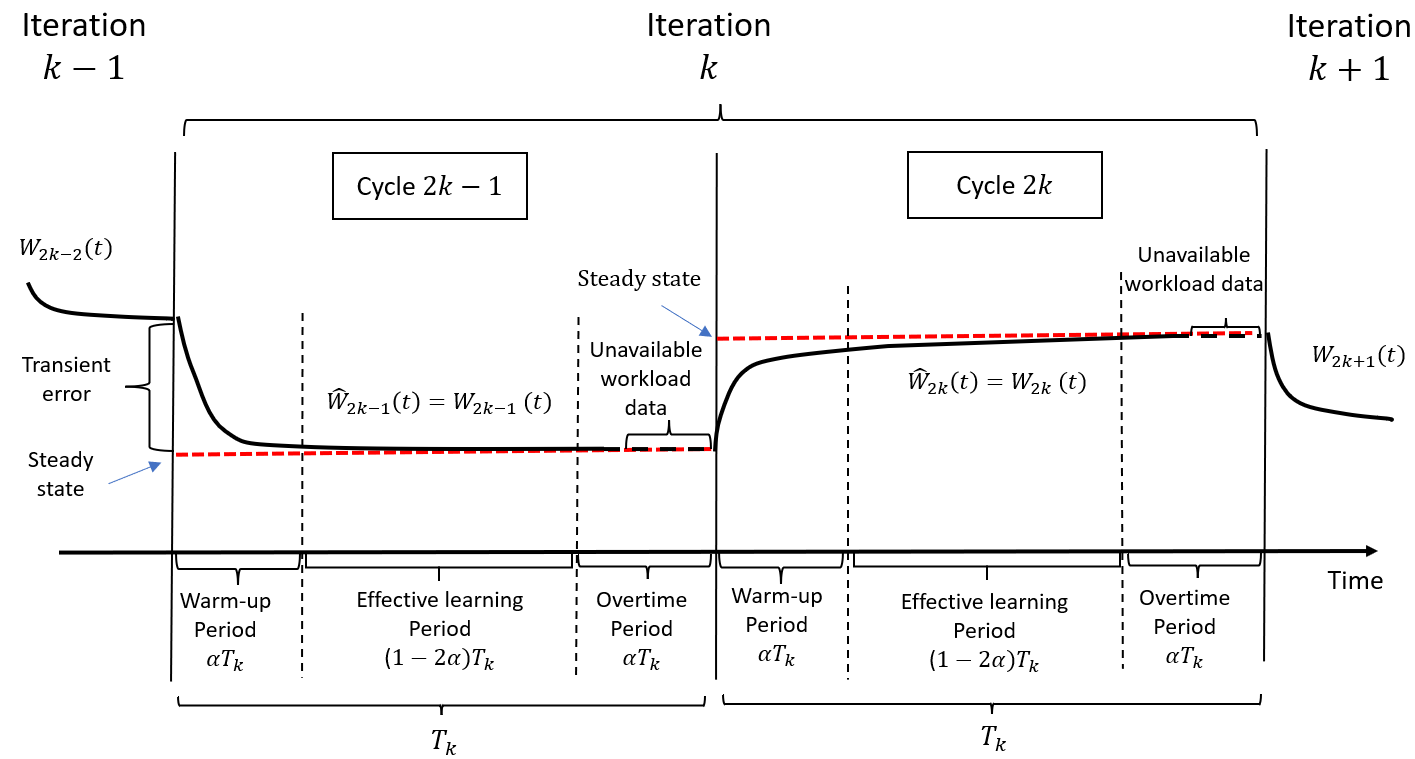}
\vspace{-0.35in}
\caption{The system dynamics under LiQUAR.}
\label{fig: algor flow}
\end{figure}

\subsubsection{\textbf{Retrieving workload data from service and arrival times.}} We assume that LiQUAR can observe each customer's arrivals in real time, but can only recover the individual workload at the service completion time. This assumption is consistent with real practice in many service systems. For example, in call center, hospital, etc., customer's individual workload is realized only after the service is completed. Hence, the workload process  $W_l(t)$ is not immediately observable at $t$. 

In LiQUAR, we approximate $W_l(t)$ by $\hat{W}_l(t)$ which we elaborate below. For given $l\geq 1$ and $t\in [0,T_k]$, if all customers arriving by time $t$ can finish service by the end of cycle $l$, then all of their service times are  realized, so we can recover $W_l(t)$ from the arrival times and service times of these customers using \eqref{eq: W SDE}. Since customers are served under FIFO, it is straightforward to see that this happens if and only if $W_l(t)\leq \mu_l(T_k-t)$,  {i.e., the  workload at time $t$ is completely processed by $T_k$.} Hence, we define the approximate workload as
%, we need to approximate $W_l(t)$ by the process $\hat{W}_l(t)$.% which is calculated using arrival times of all arriving customers and service times of customers who have finished service in cycle $l$.
%In our setting, we assume that, for each customer, we can observe the arrival time immediately upon her arrival, and the service time when the service completes. Then, the individual workload can be calculated by multiplying the service time with service rate. For given $l\geq 1$ and $t\in [0,T_k]$, if all customers that have arrived before time $t$ can finish service by the end of cycle $l$, their service times are all realized, and as a consequence, we can recover $W_l(t)$ from the arrival times and service times of these customers using \eqref{eq: W SDE}. As customers are served FIFO, one can check that this happens if and only if $W_l(t)\leq \mu_l(T_k-t)$, i.e. the presented workload at time $t$ is served by $T_k$. Following the above analysis, we define the approximated workload process as
\begin{equation}\label{eq: hat W}
\hat{W}_l(t) = \begin{cases}
W_l(t),& \text{ if }W_l(t)\leq \mu_l(T_k-t)\\
0,&\text{ otherwise.}
\end{cases}
\end{equation}
As illustrated in Figure \ref{fig: algor flow}, to reduce approximation error incurred by delayed observations  of service times, we discard the  $\hat{W}_l(t)$ data for $t\in ((1-\alpha) T_k,T_k]$; we call the subinterval $((1-\alpha)T_k,T_k]$ the overtime period in cycle $l$. The following Proposition \ref{prop: delay observation} ensures that the approximation error $|\hat{W}_l(t)-W_l(t)|$ incurred by delayed observation vanishes exponentially fast as length of the overtime period increases. This result will be used in Section \ref{sec: regret upper bound} to bound the estimation errors of  the FD gradient estimator $H_k$.
\begin{proposition}[Bound on Error of Delayed Observation]
\label{prop: delay observation}
Under Assumptions \ref{assmpt: uniform} and \ref{assmpt: light tail}, there exist some constants $M$ and $\theta_0>0$ such that, for all $l\geq 1$ and $0\leq t\leq T_k$,
$$\EE[|\hat{W}_l(t)-W_l(t)|]\leq \exp(-\theta_0\underline{\mu}/2\cdot(T_k-t))M.$$
%\red{Check: $\text{ with }\theta_1=\theta_0\underline{\mu}/2$}
\end{proposition}

{Roughly speaking, $M$ is the moment bound under the busiest traffic intensity, and $\theta_0$ is a small number depending on $\theta$ in Assumption \ref{assmpt: light tail}. %\textcolor{red}{and $\theta_0$ is in the domain of the moment generating function of certain queueing functions around the origin.} 
{The existence of them are shown in the Lemma \ref{lmm: uniform bound} in Section \ref{appx: lemma proofs}}}. %Their values are specified in Lemma \ref{lmm: uniform bound} in Appendix \ref{appx: uniform bound}. %The proof of Proposition \ref{prop: delay observation} is given in Appendix  \ref{appx: delay observation}.

\subsubsection{\textbf{Computing the gradient estimator.}} 
As illustrated in Figure \ref{fig: algor flow},  we also discard the data at the beginning of each cycle (i.e.,  $\hat{W}_l(t)$ for  $t\in[0,\alpha T_k]$ in cycle $l$) in order to reduce the bias due to transient  queueing dynamics incurred by the changes of the control parameters. We call $[0, \alpha T_k]$ the warm-up period of cycle $l$. Thus, we give the following system performance estimator under control $x_l$, $l\in\{2k-1, 2k\}$:
\begin{equation}\label{eq: f G hat}
\hat{f}^G(\mu_l,p_l)=\frac{-pN_l}{T_{k}}+\frac{h_0}{(1-2\alpha)T_{k}}\int_{\alpha T_k}^{(1-\alpha)T_{k}}\hat{W}_l(t)dt+c(\mu_l),
\end{equation}
and the corresponding FD gradient estimator
\begin{equation}\label{eq: FD estimator}
\bH_k = \frac{\bZ_k\cdot(\hat{f}^G(\mu_{2k},p_{2k})-\hat{f}^G(\mu_{2k-1},p_{2k-1}))}{\delta_k}.
\end{equation}
Unlike standard zero-order methods used in offline optimization problems, our data is generated through online interactions with the real system. Consequently, we must carefully tune our algorithm parameters to control the variance of 
$H_k$, as techniques like common random numbers can not be used to achieve variance reduction  for $\hat{f}^G(\mu_{2k},p_{2k})-\hat{f}^G(\mu_{2k-1},p_{2k-1})$. % by using  common random numbers. As a consequence,  %along with the bias, relying on theoretic regret analysis using statistical properties of queueing data.  

The psuedo code of LiQUAR is given in Algorithm \ref{alg: GOLiQ-UDS}. To complete the design of LiQUAR algorithm, we still need to specify the hyperparameters $T_k, \eta_k, \delta_k$ for $k\geq 1$.  We seek to optimize these hyperparameters in Section \ref{sec: regret upper bound} to achieve minimized regret bound.

\begin{algorithm}[h]
\SetAlgoLined
\KwIn{number of iterations $L$\;
parameters $0<\alpha<1$, and $T_k$, $\eta_k$, $\delta_k$ for $k=1,2,..,L$\;
initial value $\bar{x}_1=(\bar{\mu}_1,\bar{p}_1)$, $W_1(0)=0$\;}
\For{$k=1,2,...,L$}{
\vskip 1ex
Randomly draw $\bZ_k\in\{(0,2),(2,0)\}$\;
\vskip 1ex
\textbf{Run Cycle $2k-1$:} Run the system for $T_k$ units of time under control parameter $$\bx_{2k-1}= \bar{\bx}_k-\delta_k \bZ_k/2=(\bar{\mu}_k,\bar{p}_k)-\delta_k \bZ_k/2 ,$$ 
%Calculate $(\hat{W}_{2k-1}(t))$ according to \red{Procedure ??}\;
\vskip 1ex
\textbf{Run Cycle $2k$: } Run the system for $T_k$ units of time under control parameter $$\bx_{2k}=\bar{\bx}_k+\delta_k \bZ_k/2=(\bar{\mu}_k,\bar{p}_k)+\delta_k \bZ_k/2 ,$$ 
%Calculate $(\hat{W}_{2k}(t))$ according to \red{Procedure ??}\;
\vskip 2ex
\textbf{Compute FD gradient estimator:} 
{\small
	\begin{equation*}
		\begin{aligned}
			\bH_k=\frac{\bZ_k}{\delta_k}\Bigg[\frac{h_0}{(1-2\alpha)T_k}&\int_{\alpha T_k}^{(1-\alpha) T_k}\left(\hat{W}_{2k}(t)-\hat{W}_{2k-1}(t)\right)dt-\frac{p_{2k}N_{2k}-p_{2k-1}N_{2k-1}}{T_k}+c(\mu_{2k})-c(\mu_{2k-1})\Bigg]
		\end{aligned}
	\end{equation*}
	where $\hat{W}_l(\cdot)$ is an approximate  of $W_l(\cdot)$ as specified in \eqref{eq: hat W}.
}
%where the integral term is calculated according to \eqref{eq: integral}\;
\vskip 2ex
\textbf{Update} $\bar{\bx}_{k+1}= \Pi_\mathcal{B} (\bar{\bx}_k- \eta_k \bH_k)$.	
\vskip 1ex
}
\caption{LiQUAR}
\label{alg: GOLiQ-UDS}

\end{algorithm}
\begin{remark}[Integrating queueing features into LiQUAR]\label{remark: queue structure utilization}
%\blue{The objective function \eqref{eqObjmin} does not only involves the queueing function, it also contains a demand function $\lambda(\cdot)$ with unknown structure. Thus, the current form of our algorithm intends to ``minimize" the assumptions on the demand structure while keeping the algorithm simple and flexible to implement. 
For LiQUAR to effectively address our queueing control problem, its design should be well informed by the features and structures inherent to the underlying queueing system. 
%On the other hand, the full design of Algorithm \ref{alg: GOLiQ-UDS} actually implicitly needs to utilize some queueing structure to achieve a good performance.  Specifically, 
{First, we utilize workload data to construct the gradient estimator, moving beyond traditional reliance on arrival and service time data. This shift is particularly advantageous in heavy-traffic scenarios where queueing formulas are highly sensitive to input estimates such as demand and service distributions. To calculate the gradient estimator which involves integrating the workload process with delayed observations, we leverage the fact that the workload process is almost surely piecewise linear. This helps simplify the computation of the gradient estimator.}  
Second, to mitigate transient biases and errors caused by delayed observations in the queueing data, we exclude data from a designated warm-up interval at the start of each cycle and an over-time interval at its end. The lengths of these intervals are determined by the exponential ergodic rate of the $M/GI/1$ queue. 
Third, the efficiency of learning algorithms, particularly gradient-based methods, critically depends on the choice of hyperparameters. Leveraging the regret analysis in Section \ref{sec: regret upper bound}, which extensively utilizes queueing properties such as transient bias and autocorrelation in the data, we carefully determine optimal hyperparameters to minimize regret. 
\end{remark}

	\section{Convergence Rate and Regret Analysis} \label{sec: regret upper bound}
	%In this section, we develop the theoretic guarantees on  the performance of LiQUAR. We first establish asymptotic optimality of LiQUAR in Section \ref{subsec: convergence to optimality}, showing that $\bar{x}_k$ will converges to the optimal control $x^*$ as time goes by, and provide an expression of the convergence rate in terms of the algorithm hyperparameters. In Section \ref{subsec: regret upper bound}, we investigate the performance of dynamic pricing and capacity sizing decisions according to LiQUAR, by analyzing the total regret, which is the gap between the amount of cost produced by LiQUAR and that by the optimal control $x^*$.
	In Section \ref{subsec: convergence to optimality}, we establish the rate of convergence for our decision variables $(\mu_k,p_k)$ under LiQUAR (Theorem \ref{thm: convergence}). Besides, our analysis illustrate	how the estimation errors in the queueing data will propagate to the iteration of $(\mu_k,p_k)$ and thus affect the quality of decision making. We follow three steps:
	First, we quantify the bias and mean square error of the estimated system performance $\hat{f}^G(\mu_l,p_l)$ computed     from the queueing data via \eqref{eq: f G hat} (Proposition \ref{prop: estimation error}). To bound the estimation errors, we need to deal with the transient bias and stochastic variability in the queueing data. Next, using these estimation error bounds, we can determine the accuracy of the FD gradient estimator $H_k$ in terms of the algorithm hyperparameters (Proposition \ref{prop: bound Bk V_k}). Finally, following the convergence analysis framework of SGD algorithms, we obtain the convergence rate of LiQUAR in terms of the algorithm hyperparameters (Theorem \ref{thm: convergence}). The above three steps together form a quantitative explanation of how the errors are passed on from the queueing data to the learned decisions (whereas there is no such steps in PTO so its performance is much more sensitive to the errors in the data).
	In addition, the convergence result enables us to obtain the optimal choice of hyperparameters if the goal is to approximate $\bx^*=(\mu^*,p^*)$ accurately with minimum number of iterations, which is often preferred  in simulation-based offline learning settings.
	
	In Section \ref{subsec: regret upper bound}, we investigate the cost performance of dynamic pricing and capacity sizing decisions made by LiQUAR, via analyzing the total regret, which is the gap between the amount of cost produced by LiQUAR and that by the optimal control $\bx^*$. Utilizing the convergence rate established by Theorem \ref{thm: convergence} and a separate analysis on the transient behavior of the system dynamics under LiQUAR (Proposition \ref{prop: regret of nonstationarity}), we obtain a theoretic bound for the total regret of LiQUAR in terms of the algorithm hyperparameters. By simple optimization, we obtain an optimal choice of hyperparameters which leads to a total regret bound of order $O(\sqrt{T}\log(T))$ (Theorem \ref{thm: upper bound}), where $T$ is the total amount of time in which the system is operated by LiQUAR .

	\subsection{Convergence Rate of Decision Variables}\label{subsec: convergence to optimality}
	As $\bar{\bx}_k$ evolves according to an SGD iteration in  Algorithm \ref{alg: GOLiQ-UDS}, its convergence  depends largely on how accurate the gradient is approximated by the FD estimator $H_k$. %The accuracy of $H_k$ essentially depends on the estimation errors of  the system performance calculated from the queueing data in cycles $l=2k-1$ and $2k$, i.e. $\{\hat{W}_l(t): \alpha T_k\leq t\leq (1-\alpha)T_k\}$. % In Section \ref{subsec: estimation error}, we first establish bounds for the estimation errors of  queueing data in terms of the algorithm hyperparameters. In Section \ref{subsec: convergence to optimal solution}, we connect the estimation errors to the convergence rate of $\bar{x}_k$, i.e. the price and service capacity learned by Algorithm \ref{alg: GOLiQ-UDS}, to the optimal solution. Based on these building blocks, in Section \ref{subsec: regret upper bound}, we are able to express the regret bound in terms of the hyperparamters and obtain a ``near optimal" regret bound of order $\tilde{O}(\sqrt{T})$ for Algorithm \ref{alg: GOLiQ-UDS} under  the properly chosen  hyperparamters.
	In the theoretical analysis, the accuracy of $\bH_k$ is measured by the following two quantities: 
	\begin{equation*}
		B_k\equiv\EE\left[\|\EE[\bH_k-\nabla f(\bar{\bx}_k)| \mathcal{F}_k]\|^2\right]^{1/2}\quad \text{and}\quad \mathcal{V}_k \equiv \EE[\|\bH_k\|^2],
	\end{equation*}
	where $\mathcal{F}_k$ is the $\sigma$-algebra including all events in the first $2(k-2)$ cycles and $\|\cdot\|$ is Euclidean norm in $\mathbb{R}^2$. Intuitively, $B_k$ measures the bias of the gradient estimator $\bH_k$ and $\mathcal{V}_k$ measures its variability.
	
	According to \eqref{eq: FD estimator}, the gradient estimator $\bH_k$ is computed using the estimated system performance $\hat{f}^G(\mu_{2k},p_{2k})$ and $\hat{f}^G(\mu_{2k-1},p_{2k-1})$. So, the accuracy of $\bH_k$ essentially depends on the estimation errors of  the system performance, i.e., how close is $\hat{f}^G(\mu_l,p_l)$ to $f(\mu_l,p_l)$. Note that the control parameters  $(\mu_l,p_l)$ for $l\in\{2k-1,2k\}$ are random and dependent on the events in the first $2(k-2)$ cycles. Accordingly, we need to analyze the estimation error of $\hat{f}^G(\mu_l,p_l)$ conditional on the past events, which is also consistent with our definition of $B_k$.  %Recall that,
	%in each cycle $l=2k-1, 2k$,  the estimated system performance is computed from the queueing data as follows
	%$$
	%\hat{f}^G(\mu_l,p_l)=\frac{-pN_l}{T_{k}}+\frac{h_0}{(1-2\alpha)T_{k}}\int_{\alpha T_k}^{(1-\alpha)T_{k}}\hat{W}_l(t)dt+c(\mu_l).
	%$$
	%The following Proposition \ref{prop: estimation error} establish bounds on the conditional bias and mean square error of $\hat{f}^G(\mu_l,p_l)$, for any given control parameter $(\mu_l, p_l)$ and initial workload $W_l(0)$, in terms of the algorithm hyperparameter $T_k$.
	For this purpose, we denote by $\mathcal{G}_l$  the $\sigma$-algebra including all events in the first $l-1$ cycles and write $\EE_l[\cdot]\equiv\EE[\cdot\vert \mathcal{G}_l]$. The following Proposition \ref{prop: estimation error} establishes bounds on the conditional bias and mean square error of $\hat{f}^G(\mu_l,p_l)$, in terms of %the control parameter $(\mu_l, p_l)$, 
	the initial workload $W_l(0)$ and the  hyperparameter $T_k$.
	
	\begin{proposition}[Estimation Errors of System Performance]\label{prop: estimation error} Under Assumptions \ref{assmpt: uniform} and \ref{assmpt: light tail}, for any $T_k>0$, the bias and mean square error of  $\hat{f}^G(\mu_l,p_l)$, conditional on $\mathcal{G}_l$, have the following bounds: 
		\begin{enumerate}
			\item Bias$$\left|\EE_l\left[\hat{f}^G(\mu_{l},p_{l})-f(\mu_{l},p_{l}) \right]\right|\leq \frac{2\exp(-\theta_1\alpha T_k)}{(1-2\alpha)\theta_1 T_k}\cdot M(M+W_l(0))(\exp(\theta_0 W_l(0))+M).$$
			\item Mean square error
			$$\EE_l[(\hat{f}^G(\mu_l,p_l)-f(\mu_l,p_l))^2] \leq K_MT_k^{-1}(W_l^2(0)+1)\exp(\theta_0 W_l(0)),$$
		\end{enumerate}
		where $\theta_1 \equiv \min(\gamma, \theta_0\underline{\mu}/2)%\theta_0\underline{\mu}/2)
		$ , and $\gamma$ and $K_M$ are two positive constants that are independent of $l, T_k, W_l(0),  \mu_l$ and $p_l$. 
	\end{proposition}
	The proof of Proposition \ref{prop: estimation error} is given in Section \ref{sec: proofs}, where the specification of the two constants $\gamma$ and $K_M$ are given in \eqref{eq: gamma choice} and \eqref{eq: KM} respectively. The key step in the proof is to bound the transient bias (from the steady-state distribution) and auto-correlation of the workload process $\{W_l(t):0\leq t\leq T_k\}$, utilizing an ergodicity analysis. This approach can be applied to other queueing models which share similar ergodicity properties, e.g., GI/GI/1 queue and  stochastic networks \citep{blanchet2020rates}.  %to derive bounds on the estimation errors of system performance.
	
	Based on Proposition \ref{prop: estimation error}, we establish the following bounds on $B_k$ and $\mathcal{V}_k$ in terms of the algorithm hyperparameters $T_k$ and $\delta_k$.
	\begin{proposition}[Bounds for $B_k$ and $\mathcal{V}_k$]\label{prop: bound Bk V_k}
		Under Assumptions \ref{assmpt: uniform} and \ref{assmpt: light tail}, the bias and variance of the gradient estimator satisfy
		\begin{align}\label{eq: B_k V_k bound}
			B_k=O\left(\delta_k^2+\delta_k^{-1}\exp(-\theta_1 \alpha T_k)\right),\quad \mathcal{V}_k= O\left(\delta_k^{-2}T_k^{-1}\vee 1\right).
		\end{align}
		
	\end{proposition}

	{Assumption \ref{assmpt: uniform} guarantees that the objective function  $f(\mu,p)$ in \eqref{eqObjmin} has desired convex structure (see Lemma 
		\ref{lmm: convexity} in Section \ref{sec: proofs} for details).} Hence, the SGD iteration is guaranteed to converge to its optimal solution $x^*$ as long as the gradient bias $B_k$ and variance $\mathcal{V}_k$ are properly bounded.  Utilizing the bounds on $B_k$ and $\mathcal{V}_k$ as given in  Proposition \ref{prop: bound Bk V_k}, we are able to prove the convergence of LiQUAR and obtain an explicit expression of the convergence rate in terms of algorithm hyperparameters. 
	\begin{thm}[Convergence rate of decision variables]\label{thm: convergence}
		Suppose Assumption \ref{assmpt: uniform} holds. If there exists a constant $\beta\in (0,1]$ such that the following inequalities  hold for all $k$ large enough:
		\begin{equation}\label{eq: beta conditions}
			\left(1+\frac{1}{k}\right)^{\beta}\leq 1+\frac{K_0}{2}\eta_k, \quad B_k\leq \frac{K_0}{8}k^{-\beta}, \quad \eta_k\mathcal{V}_k = O(k^{-\beta}).
		\end{equation}
		Then,  we have 
		\begin{equation}\label{eq: convergence beta}
			\EE\left[\|\bar{\bx}_k-\bx^*\|^2\right]=O(k^{-\beta}).
		\end{equation}
		If, in further, Assumption \ref{assmpt: light tail} holds and the algorithm hyperparameters are set as $\eta_k=O(k^{-a})$, $T_k = O(k^b)$, and $\delta_k=O(k^{-c})$ for some constants $a,b,c\in(0,1]$. We have
		\begin{equation}\label{eq: convergence rate}
			\EE\left[\|\bar{\bx}_k-\bx^*\|^2\right]=O\left(k^{\max(-a,-a-b+2c,-2c)}\right).
		\end{equation}

		%			\begin{enumerate}
			%				\item[$(a)$]  $\left(1+\frac{1}{k}\right)^{\beta}\leq 1+\frac{K_0}{2}\eta_k$,
			%				\item[$(b)$] $B_k\leq \frac{K_0}{8}k^{-\beta}$,
			%				\item[$(c)$] $\eta_k\mathcal{V}_k = O(k^{-\beta})$.
			%			\end{enumerate}
		
	\end{thm}
	\begin{remark}[optimal convergence rate]\label{rmk: convergence} According to the bound \eqref{eq: convergence rate}, by minimizing the term $\max(-a,a-b+2c, -2c)$, one can obtain an optimal choice of hyperparameters $\eta_k = O(k^{-1}), T_k=O(k)$ and $\delta_k=O(k^{-1/2})$ under which the decision parameter $\bx_k$ converges to $\bx^*$ at a fastest rate of $O(L^{-1})$, in terms of the total number of iterations $L$. 
		Of course, the above convergence rate analysis does not focus on reducing the total system cost generated through the learning process, which is what we will do in Section \ref{subsec: regret upper bound}.
	\end{remark}

	\subsection{Regret Analysis}\label{subsec: regret upper bound}
	Having established the convergence of control parameters 
	under Assumption \ref{assmpt: uniform}, we next investigate
	the efficacy of LiQUAR as measured
	by the cumulative regret which measures the gap between the cost under LiQUAR and that under the optimal control. 
	According to the system dynamics described in Section \ref{subsec: system dynamics}, under LiQUAR, the expected cost incurred in cycle $l$ is
	\begin{align}\label{regret_rhol}
		\rho_l \equiv \EE\left[h_0\int_0^{T_k}W_l(t)dt+c(\mu_l)T_k-p_lN_l\right],
	\end{align}
	where $k = \lceil l/2\rceil $. The total regret in the first $L$ iterations (each iteration contains two cycles) is
	$$R(L)=\sum_{k=1}^L\sum_{l=2k-1}^{2k}R_l=\sum_{l=1}^{2L}R_l,\quad \text{with } R_l\equiv \rho_{l}-T_kf(\mu^*,\rho^*).$$
	Our main idea  is to 
	separate the total regret $R(L)$ into three parts as
	\begin{align}\label{eq: regret decomposition}
		R(L) = &\sum_{k=1}^L \underbrace{\EE\left[2T_{k}(f(\bar{\bx}_k)-f(\bx^*))\right]}_{\equiv R_{1k}\text{: regret of suboptimality}}\notag\\
		&+\sum_{k=1}^L\underbrace{\EE\left[\left(\rho_{2k-1}  -T_{k}f(\bx_{2k-1})\right)+\left(\rho_{2k}  -T_{k}f(\bx_{2k})\right)\right] }_{\equiv R_{2k}\text{: regret of nonstationarity}}\\
		&+\sum_{k=1}^L\underbrace{\EE\left[T_{k}(f(\bx_{2k-1}) + f(\bx_{2k}) - 2f(\bar{\bx}_k))\right] }_{\equiv R_{3k}\text{: regret of {finite difference} }}\notag,
	\end{align}
	%Given the previous results on the convergence rate and estimation error bounds, we are able to develop an explicit expression for the asymptotic order of regret upper bound in terms of the algorithm hyperparameters. Then, we utilize this expression to optimize our choice of hyperparameters by minimizing the asymptotic order of the regret bound. Thus, we finalize the algorithm design and obtain the theoretic regret upper bound.
	{which arise from the errors due to the suboptimal decisions ($R_{1k}$),  the transient system dynamics ($R_{2k}$),  and  the estimation of gradient ($R_{3k}$), respectively.}
	Then we aim to minimize the orders of all three regret terms by selecting the ``optimal'' algorithm hyperparameters $T_k, \eta_k$ and $\delta_k$ for $k\geq 1$. %Then, based on such connection,  we can optimize the choice of algorithm hyperparameters and obtain the regret bound of LiQUAR. 
	
	\paragraph{\textbf{Treating $R_{1k},R_{2k},R_{3k}$ separately.} } Suppose the hyperparameters of LiQUAR are set in the form of $\eta_k=O(k^{-a})$, $T_k = O(k^b)$, and $\delta_k=O(k^{-c})$ for some constants $a,b,c\in(0,1]$. The first regret term  $R_{1k}$ is determined by the convergence rate of control parameter $\bar{x}_k$. By Taylor's expansion, $f(\bar{\bx}_k)-f(\bx^*)= O(\|\bar{\bx}_k-\bx^*\|_2^2 )$, and hence, $R_{1k} = O(T_k\|\bar{\bx}_k-\bx^*\|_2^2 )$. Following Theorem \ref{thm: convergence}, we have $R_{1k} = O(k^{\max(b-a,b-2c,-a+2c)})$. By the smoothness condition in Assumption \ref{assmpt: uniform}, we can check that $R_{3k} = O(T_k\delta_k^2)= O(k^{b-2c})$ (Lemma \ref{lmm: exploration cost} in Section \ref{sec: proofs}). 
	
	The remaining regret analysis will focus on  the regret of nonstationarity $R_{2k}$. Intuitively, it depends on the rate at which the (transient) queueing dynamics converges to its steady state. Applying the same ergodicity analysis as used in the analysis of estimation errors of system performance, we can find a proper bound on the transient bias after the warm-up period, i.e., for $W_l(t)$ with $t\geq \alpha T_k$. Derivation of a desirable bound on the transient bias in the warm-up period, i.e., for $W_l(t)$ with $t\in[0,\alpha T_k]$, is less straightforward. The main idea is based on the two facts that (1) $W_l(t)$, when $t$ is small, is close to the steady-state workload corresponding to $(\mu_{l-1},p_{l-1})$ and that (2) the steady-state workload corresponding to $(\mu_{l-1},p_{l-1})$ is close to that of  $(\mu_l,p_l)$. We formalize the bound on $R_{2k}$ in Proposition \ref{prop: regret of nonstationarity} below. The complete proof  is given in Section \ref{subsec: proof of regret bound}. 
	\begin{proposition}[Regret of Nonstationarity]\label{prop: regret of nonstationarity}
		Suppose Assumptions \ref{assmpt: uniform} and \ref{assmpt: light tail} hold. If $T_k>\log(k)/\gamma$
		and  there exists some constant $\xi\in(0,1]$ such that $\max(\eta_k\sqrt{\mathcal{V}_k},\delta_k)=O(k^{-\xi})$.
		Then,
		\begin{equation}\label{eq: R_1k}
			R_{2k} = O\left(k^{-\xi}\log(k)\right).
		\end{equation}
		If, in further, the algorithm hyperparameters are set as $\eta_k=O(k^{-a})$, $T_k = O(k^b)$, and $\delta_k=O(k^{-c})$ for some constants $a,b,c\in(0,1]$, we have
		\begin{equation*}
			R_{2k} = O\left(k^{\max(-a-b/2+c,-a, -c)}\log(k)\right).
		\end{equation*}
	\end{proposition}
	
	By summing up the three regret terms, we can conclude that
	$$R(L)\leq \sum_{k=1}^L C\left(k^{\max(-a-b/2+c,-a, -c)}\log(k) +k^{\max (b-a,b-2c, -a+2c)}+k^{b-2c}\right),$$
	for some positive constant $C$ that is large enough. The order of the upper bound on the right hand side reaches its minimum at $(a,b,c)=(1,1/3,1/3)$. The corresponding total regret and time elapsed in the first $L$ iterations are, respectively,
	$$R(L)=O(L^{2/3}\log(L)) \quad \text{ and }\quad T(L)=O(L^{4/3}).$$
	As a consequence, we have $R(T) = O(\sqrt{T}\log(T))$. 
	
	\begin{thm}[Regret Upper Bound]\label{thm: upper bound}
		Suppose Assumptions \ref{assmpt: uniform} and \ref{assmpt: light tail} hold. If we choose $\eta_k = c_\eta k^{-1}$ for some $c_\eta> 2/K_0$, $T_k = c_Tk^{1/3}$ for {some $ c_T> 0$} and $\delta_k = c_\delta k^{1/3}$ for {some $0<c_\delta<\sqrt{K_0/32c}$}, where $c$ is a smoothness constant given in {Lemma \ref{lmm: fd approximation},} then the total regret accumulated in the first $L$ rounds by LiQUAR
		$$R(L)=O(L^{2/3}\log(L))=O(\sqrt{T(L)}\log(T(L))).$$ Here $T(L)$ is the total units of time elapsed in $L$ cycles. 
	\end{thm}
	\begin{remark}[On the $O(\sqrt{T}\log(T))$ Regret Bound]
		Consider a hypothetical setting in which we are no longer concerned with the transient behavior of the queueing system, i.e., somehow we can directly observe  an unbiased and independent sample of the objective function with uniform bounded variance in each iteration. In this case, we know that the Kiefer-Wolfowitz algorithm and its variate provide an {effective approach for model-free stochastic optimization}  \citep{BoradieZeevi2011}. %Besides, it is equivalent to measure the regret by the number of iterations or by the units of time elapsed. 
		According to  \cite{BoradieZeevi2011},  the convergence rate of Kiefer-Wolfowitz algorithm is $\|\bar{\bx}_k-\bx^*\|^2=O(\eta_k/\delta_k^2)$. In addition, the regret of {finite difference} is $f(\bx_{2k-1})+f(\bx_{2k})-2f(\bar{\bx}_k) = O(\delta_k^2)$. Since $\eta_k/\delta_k^2+ \delta_k^2\geq 2\sqrt{\eta_k}\geq k^{-1/2}$, we can conclude that the optimal convergence rate in such a hypothetical setting is $O(k^{-1/2})$. This accounts for the $\sqrt{T}$ part of our regret in Theorem \ref{thm: upper bound}. Unfortunately, unlike the hypothetical setting, our queueing samples are biased and correlated. Such a complication is due to the nonstationary error at the beginning of cycles which gives rise to the extra $\log(T)$ term in the regret bound; see Proposition \ref{prop: regret of nonstationarity} for additional discussion of the $\log(T)$ term in our regret.
	\end{remark}

	\subsection{LIQUAR in Heavy Traffic}\label{sec: heavy traffic}
	We now evaluate LiQUAR's performance under heavy traffic conditions. To do this, we construct a series of queueing models with traffic intensities that approach the critical threshold of 1, and define the associated profit optimization problems. Our goal is to derive an explicit regret expression as a function of traffic intensity. To achieve this, we will need to consider a simplified model in order for reduced technicalities in our regret analysis.
	We consider two simplifications: first, we focus on the case of exponential service times; second, we focus on a pricing problem, treating the service rate $\mu$ as a constant. We stress that the simplified model still preserve the challenge of dealing with unknown demand (the core aspect of the learning problem). 
	
	Consider a sequence of $M/M/1$ systems indexed by the parameter $h>0$, which is the queueing congestion cost per time unit. They share a common demand curve $\lambda(p)$ and service rate $\mu=1$. For the $h^{\rm th}$ model, we aim to  minimize the objective function below
	\begin{equation}\label{eq: objct heavy traffic}
		\min_{p\in\mathcal{B}_h}\quad f_h(p)\equiv -p\lambda(p)+\frac{h\lambda(p)}{1-\lambda(p)},
	\end{equation}
	where %$h$ is the queueing congestion cost per time unit, and
	$\mathcal{B}_h$ will be specified later.
	When the holding cost $h$ in \eqref{eq: objct heavy traffic} decreases, the service provider is incentivized to increase service utilization to maximize profit which places the system under the heavy-traffic regime. Below we will formally show that the traffic intensity under the optimal price $\rho^*_h$ converges to $1$ as $h\rightarrow0$, specifically, below we will show $1-\rho^*_h=O(\sqrt{h})$.

	%In a moment, we will show that under moderate assumptions on demand function, the optimal service excess 
	We denote by $p^*_h$ as the optimal solution to \eqref{eq: objct heavy traffic}. To explicitly show the relationship between $\rho^*_h$ and $h$, we impose the following assumptions on the demand curve.
	
	\begin{assumption}[Demand Curve] \label{ass: demand heavy traffic}We assume that the arrival rate function $\lambda(p)$ satisfies the following conditions:
		\begin{enumerate}
			\item The demand function is non-increasing and twice-differentiable.
			\item The function $r(p)=p\lambda(p)$ is strictly concave. %, i.e., $2\lambda'(p)^2+p\lambda''(p)>0$.
			\item {The demand function is elastic in the feasible regions: $-\frac{\lambda'(p)}{\lambda(p)}\cdot p>1$ for $p\in\mathcal{B}_h$. }
		\end{enumerate}
	\end{assumption}
	{The third technical condition is commonly used in the literature of revenue management in queues; see for example, Assumption 1 in \cite{Maglaras2003}. Essentially, this condition assumes that customers are price sensitive in the feasible region.}
	
	% \red{
		% 	\begin{remark}
			% 		Our heavy traffic setting, including the objective function and Assumption 3, is consistent with the literature on pricing the capacity sizing for queueing systems with large market size. In particular, letting $h\to 0$ in our setting is equivalent to multiplying the demand function $\lambda(p)$ and service rate $\mu$ by $1/h$, i.e., sending to market scale to infinity. Assumption 3, especially condition (3), has also been used in previous works, see for example \cite{Maglaras2003}, and drives the economically optimal system to heavy traffic when market size grows up.
			% 	\end{remark}
		% }
	Denote $p_0$ as the price that makes the system critically loaded, i.e., $\lambda(p_0)=1$.
	Before giving our regret bound, we first characterize the optimal pricing decision $p_h^*$ as a function of $h$ and relate the traffic intensity $\rho_h^*$ to $h$.
	\begin{proposition}\label{prop: heavy traffic problem rho}
		Under Assumption \ref{ass: demand heavy traffic}, we have the optimal price 
		$$p^*_h \equiv \arg\min f_h(p)= p_0 +\sqrt{h}\cdot\sqrt{\frac{1}{(1+p_0\lambda'(p_0))\lambda'(p_0)}}+o(\sqrt{h}),$$
		and the corresponding optimal service excess
		$$1-\rho_h^*=\sqrt{h}\cdot \sqrt{\frac{-\lambda'(p_0)}{1+p_0\lambda'(p_0)}}+o(\sqrt{h})=O(\sqrt{h}).$$
		%\textcolor{red}{with $1+p_0\lambda'(p_0)<0$.}		
	\end{proposition}
	%	Using Taylor's expansion, we can conclude that, as $h\to 0$,
	%	$$p^*_h \equiv \arg\min f_h(p)= p^* +\sqrt{h}\cdot\sqrt{\frac{1}{(1+p^*\lambda'(p^*))\lambda'(p^*)}}+o(\sqrt{h}),$$
	%	where $p^*$ and $\lambda(p)$ satisfies
	%	\begin{equation}\label{eq: condition on lambda in heavy traffic}
		%		\lambda(p^*)=1,\quad p^*\lambda'(p^*)<-1.
		%	\end{equation}
	
	To investigate LiQUAR's performance in heavy traffic, we consider  $\mathcal{B}_h$ which asymptotically operates the system in heavy traffic as $h\rightarrow 0$. Following Proposition \ref{prop: heavy traffic problem rho}, we let
	$$p^*_h\in \left[~p_0+ c_1\sqrt{h}, p_0+c_2\sqrt{h}~\right]\equiv \mathcal{B}_h,$$
	where the constants $c_1$ and $c_2$ satisfy $0<c_1<c_0\equiv1/\sqrt{\lambda'(p_0)(1+p_0\lambda'(p_0))}<c_2.$ 
	%Then, we apply LiQUAR to solve the problem sequences \eqref{eq: objct heavy traffic} with $h\rightarrow0$. Before giving the formal theorem for LiQUAR under heavy-traffic, we introduce a heavy-traffic learning scheme. 
	Note that as $\rho_h^*\rightarrow1$ (or equivalently, $h\rightarrow0$), {the queueing system takes a longer time to converge to its steady-state}. Therefore, as $h\rightarrow0$, we need to increase the length of the learning cycle. %To explicitly show the influence due to heavy-traffic, we consider $T_0$ as time coefficient independent with $h$ and scale time by $T^h\equiv T_0/h$. 
	We operate the $h^{\rm th}$ model under LiQUAR with a total time duration of $T^h\equiv T_0/h$ units where $T_0$ is a positive constant independent of $h$. We evaluate our performance using the regret
	%, we operate LiQUAR to solve \eqref{eq: objct heavy traffic} for the time $T^h$ and the performance metric of the $h$-th system is measured by regret 
	$$R^h(T_0)\equiv R(T^h).$$
	%	to show that this regret is measured in our heavy-traffic scheme indexed by $h$ with time $T_0$ independent of traffic intensity. 
	In what follows, we report the theoretical regret bound as a function of the traffic intensity. %Then, by properly choosing the hyperparameters of LiQUAR according to the holding cost $h$, we could give a theoretical upper bound of LiQUAR in heavy-traffic as below:
	%Then, we intend to apply LiQUAR to solve the following problem 
	%	\begin{equation*}
		%		\min_{p\in\mathcal{B}_h} f_h(p)=-p\lambda(p)+h\frac{\lambda(p)}{\mu-\lambda(p)},
		%	\end{equation*}
	%	with $\mu=1$ fixed. Following the analysis in our main paper, we have the following Theorem \ref{thm: heavy traffic} showing the performance of LiQUAR in heavy-traffic (as $h\rightarrow0$).
	\begin{thm}[Regret Bound in Heavy Traffic]\label{thm: heavy traffic}	
		For the $h^{\rm th}$ system operated under LiQUAR to minimize \eqref{eq: objct heavy traffic}, when the hyper-parameters are chosen as
		\begin{equation}
			T^h_k = c_Th^{-1}k^{1/3},\quad\delta^h_k=c_\delta\sqrt{h}k^{-1/3},\quad \eta^h_k=c_\eta\sqrt{h}k^{-1},
		\end{equation}
		where the constants $c_T, c_\delta$ and $c_\eta$ are independent of $h$, then the  regret is given by 
		\begin{equation}
			R^h(T_0) \leq C\sqrt{h^{-1}T_0\log T_0}=O(\sqrt{T_0\log T_0}/(1-\rho_h^*)),
		\end{equation} 
		where $C$ is a positive constant independent of $h$. 
	\end{thm}
	{The result in Theorem \ref{thm: heavy traffic} refines that in Theorem \ref{thm: upper bound} by emphasizing how the regret depends on the system's traffic intensity.} In Section 7, we also conduct heavy-traffic analysis for a nonparametric PTO framework; and we compare the performance under both methods in our heavy-traffic regime both theoretically and numerically.

	%%%%%%%%%%%%%%%%%%%%%%%%%%%%%%%%%%%%%%%%%%%%%
	%%%%%%%%%%%%%%%%%%%%%%%%%%%%%%%%%%%%%%%%%%%%%
	\section{Numerical Experiments} \label{sec: Numerical}
	We provide engineering confirmations of the effectiveness of LiQUAR by conducting a series of numerical experiments. We will use simulated data to visualize the convergence of LiQUAR, estimate the regret curves and benchmark them with our theoretical bounds. 
	In Section \ref{subsec: basic exp}, we evaluate the performance of LiQUAR using an $M/M/1$ base example with logit demand functions. 
	In Section \ref{subsec: hyperparameters}, we discuss how to fine-tune the algorithm's hyperparameters including $T_k$ and $\eta_k$. 
	%Then, we compare the performance of LiQUAR and PTO and investigate the impact of traffic intensity $\rho^*$ on both methods in Section \ref{subsec: PTO}.
	In Section \ref{subsec: GG1}, we generalize LiQUAR to $GI/GI/1$ queues with non-Poisson arrivals and evaluate its performance.

	\subsection{An $M/M/1$ base example}\label{subsec: basic exp}
	Our base model is an $M/M/1$ queues having Poisson arrivals with rate $\lambda(p)$ and exponential service times with rate $\mu$. We consider a logistic demand function %\eqref{eq: logit demand} 
	\citep{BesbesZeevi2015}
	\begin{equation}\label{eq: logit demand}
		\lambda(p)=M_0\cdot\frac{\exp(a-bp)}{1+\exp(a-bp)},
	\end{equation}
	with $M_0=10,a=4.1,b=1$ %\citep{BesbesZeevi2015}, %\eqref{eq: logit demand} with $M=10,a=4.1$ and $b=1$ and
	and a linear staffing cost function 
	\begin{equation}\label{eq: cost function}
		c(\mu)=c_0\mu.
	\end{equation}
	The demand function is shown in the top left panel in Figure \ref{fig: basic plot}. Then, the service provider's profit optimization problem \eqref{objective} reduces to
	\begin{equation}\label{eq: base obj}
		\max_{\mu,p}\left\{p\lambda(p)-h_0\frac{\lambda(p)/\mu}{1-\lambda(p)/\mu}-c_0\mu\right\}.
	\end{equation}
	
	\subsubsection{Performance sensitivity to parameter errors without learning}\label{sec:PTOSensitivity}
	We first illustrate %the aforementioned curse of inaccuracy in heavy traffic 
	how the parameter estimation error impacts the performance. 
	%we consider an $M/M/1$ example with a logit demand function  \citep{BesbesZeevi2015}
	%\begin{equation}\label{eq: logit demand}
	%	\lambda(p)=M_0\cdot\frac{\exp(a-bp)}{1+\exp(a-bp)}, 
	%\end{equation}
	%and a linear cost $c(\mu)=\mu$. 
	Here we assume the service provider does not know the true value of $\lambda(p)$ but rather make decisions based on an estimated  arrival rate $\hat{\lambda}_{\epsilon}(p) \equiv (1-\epsilon\%)\lambda(p)$, where $\epsilon$ is the percentage estimator error. Let $(\hat{\mu}_{\epsilon},\hat{p}_{\epsilon})$ and $(\mu^*,p^*)$ be the solutions under the estimated $\hat{\lambda}_{\epsilon}$ and the true value of $\lambda$. We next compute the relative profit loss due to the misspecification of the demand function $\left(\mathcal{P}(\mu^*,p^*)-\mathcal{P}(\hat{\mu}_{\epsilon},\hat{p}_{\epsilon})\right)/\mathcal{P}(\mu^*,p^*)$, which is the relative difference between profit under the miscalculated solutions using the believed $\hat{\lambda}_{\epsilon}$ and the true optimal profit under $\lambda$. %, normalized by the true optimal profit. 

\begin{figure}
		\centering
		\includegraphics[width=0.9\linewidth]{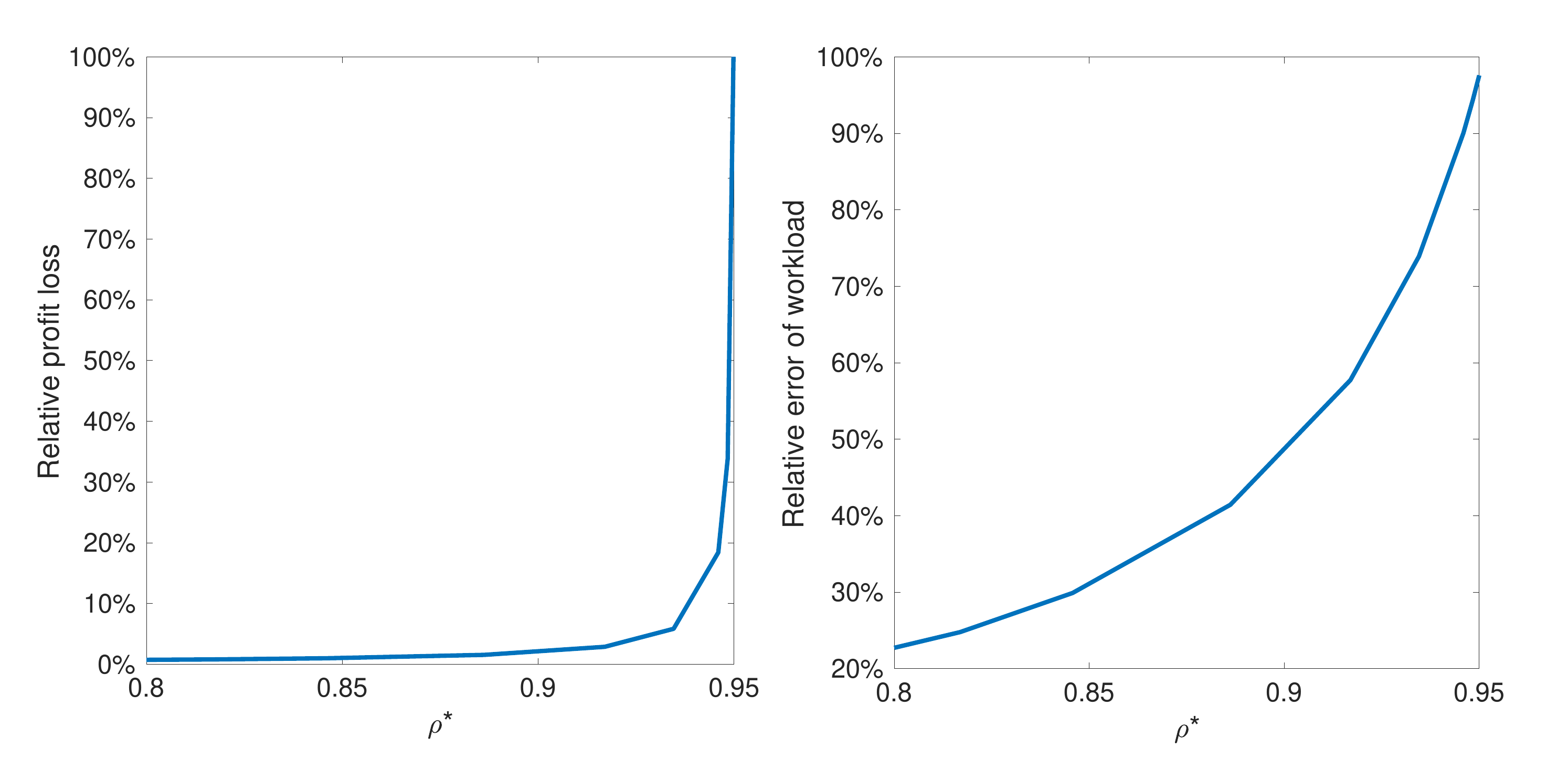}
		\caption{Relative profit loss (left) and workload error (right) for the $M/M/1$ example with $M_0=10, a=4.1$, and $b=1$ and linear staffing cost $c(\mu)=\mu$. }
		\label{fig: intro relative loss}
	\end{figure}

	Let $\rho^*\equiv \lambda(p^*)/\mu^*$ be the traffic intensity under the true optimal solution. 
	We are able to impact the value of $\rho^*$ by varying the queueing penalty {coefficient} $h_0$. We provide an illustration in Figure \ref{fig: intro relative loss} with $\epsilon = 5$. From the left panel of Figure \ref{fig: intro relative loss}, we can see that as $\rho^*$ increases, the model fidelity becomes more sensitive to the misspecification error in the demand rate and the relative loss of profit grows dramatically as $\rho^*$ goes closer to 1. 
	This effect arises from the fact that the error predicted workload is extremely sensitive to that in the arrival rate and is disproportionally amplified by the PK formula when the system is in heavy traffic (see panel (b) for the relative error of the workload). Later in Section \ref{sec:PTO_compare}, we will conduct a careful comparison to the PTO method where we will compute the PTO regret including profit losses in both the prediction and optimization steps.
	%because the estimation errors of the input parameters may be disproportionally amplified by this formula.
	%From the right panel of Figure \ref{fig: intro relative loss}, we could see that the drastic loss is originating from the underestimation of the workload. As the $\rho^*\rightarrow1$, the negative impact of the underestimation becomes more obvious. Therefore, an adaptive learning method based directly on workload is needed to avoid this sensitivity issue, which is the point of departure for our article.

	\subsubsection{Performance of LiQUAR}
\begin{figure}[t]
		\centering
		\includegraphics[width=.98\linewidth]{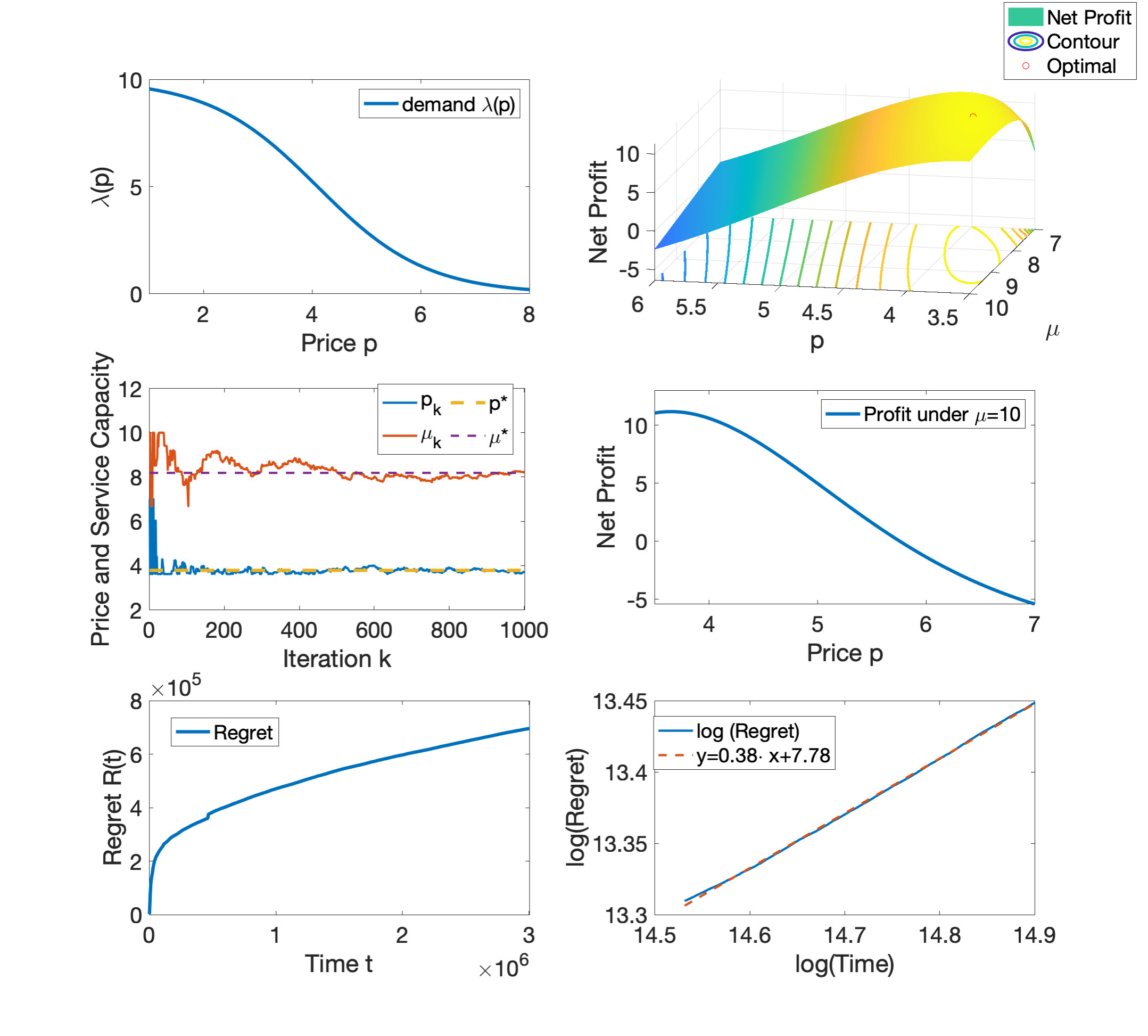}
		\vspace{-0.3in}
		\caption{Joint pricing and staffing in the $M/M/1$ logistic demand base example with  $\eta_k=4k^{-1},\delta_k=\min(0.1,0.5k^{-1/3})$, $T_k=200k^{1/3}$, $p_0=5$, $\mu_0=10$ and $\alpha=0.1$: (i) Demand function $\lambda(p)$ (top left panel); (ii) net profit function (top right panel); (iii) sample trajectories of decision parameters (middle left); (iv) One dimensional net profit function when $\mu=10$;  (v) average regret curve estimated by 100 independent runs (bottom left); (vi) a linear fit to the regret curve in logarithm scale. }
		\label{fig: basic plot}
	\end{figure}

	Using the explicit forms of  \eqref{eq: base obj}, we first numerically obtain the exact optimal solution $(\mu^*,p^*)$ and the maximum profit $\mathcal{P}(\mu^*,p^*)$ which will serve as benchmarks for LiQUAR. Taking $h_0=1$ and $c_0=1$ yields $(\mu^*,p^*)=(8.18,3.79)$, and the corresponding profit plot is shown in the top right panel of Figure \ref{fig: basic plot}. To test the criticality of condition (b) in Assumption \ref{assmpt: uniform}, we implement LiQUAR when condition (b) does not hold. For this purpose, we set $\mathcal{B}=[6.5,10]\times[3.5,7]$, in which the objective \eqref{eqObjmin} is not always convex, let alone the condition (b) of Assumption \ref{assmpt: uniform} (top right and middle right panel of Figure \ref{fig: basic plot}).
	
	Then we implement LiQUAR without exploiting the  specific knowledge of the exponential service distribution or the form of $\lambda(p)$. 
	In light of Theorem \ref{thm: upper bound}, we set the hyperparameters  $\eta_k=4k^{-1},\delta_k=\min(0.1,0.5k^{-1/3})$, $T_k=200k^{1/3}$ and $\alpha=0.1$. From Figure \ref{fig: basic plot}, we observe that the pair $(\mu_k,p_k)$, despite some stochastic fluctuations, converges to the optimal decision rapidly. The regret is estimated by averaging 100 sample paths and showed in the bottom left panel of Figure \ref{fig: basic plot}. To better relate the regret curve to its theoretical bounds as established in Theorem \ref{thm: upper bound}, 
	we also draw the logarithm of regret as a function of the logarithm of the total time; we fit the log-log curve to a straight line (bottom right panel of Figure \ref{fig: basic plot}) so that the slope of the line may be used to quantify the theoretic order of regret: the fitted slope ($0.38$) is less than its theoretical upper bound ($0.5$). Such ``overperformance" is not too surprising because the theoretic regret bound is established based on a worst-case analysis. In summary, our numerical experiment shows that the technical condition (b) in Assumption \ref{assmpt: uniform} does not seem to be too restrictive.
	\subsection{Tuning the hyperparameters for LiQUAR} \label{subsec: hyperparameters}
	Next, we test the performance of LiQUAR on the base $M/M/1$ example under different hyperparameters. 
	%The numerical results indicate that the performance of GOLiQ-UDS is overall robust with respect to the hyperparameters. In addition, w
	We also provide some general guidelines on the choices of hyperparameters when applying LiQUAR in practice.  %Specifically, we consider the joint pricing and staffing problem for an $M/M/1$ queue in Subsection \ref{subsec: basic exp}. For this problem, GOLiQ-UDS is performed under different hyperparameters and we observe the regret curves to gain some empirical insight into hyperparameters' choices. In Section \ref{subsubsec: step length insight}, we investigate the influence of $\eta_k$ and $\delta_k$. The comparison of the regret under different $T_k$ is reported in \ref{subsubsec: test T}.
	\subsubsection{Step lengths $\eta_k$ and $\delta_k$.}\label{subsubsec: step length insight}
\begin{figure}[t]
		\centering
		\includegraphics[width=\linewidth]{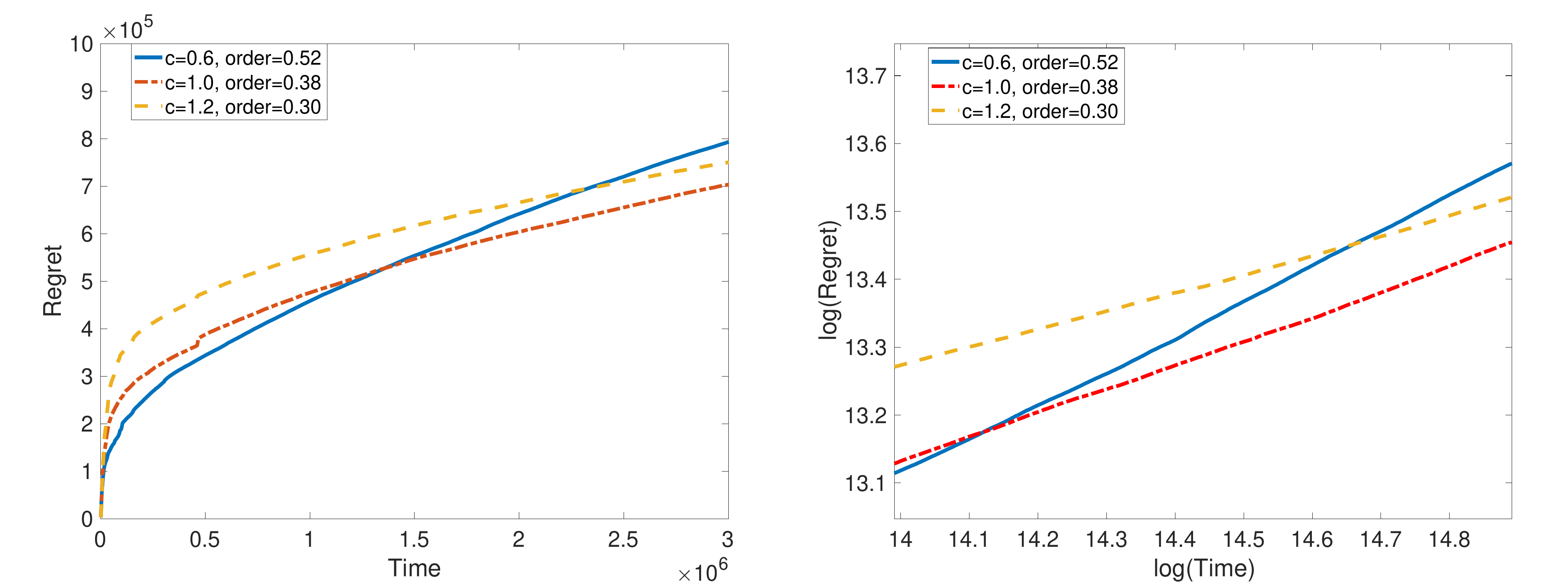}
		\caption{Regret under different $c\in\{0.6,1.0,1.2\}$: (i) average regret from 100 independent runs (left panel); (ii) regret curve in logarithm scale, with $T_k=200k^{1/3}$, $\eta_k=c\cdot 4k^{-1}, \delta_k=\min(0.1,c\cdot0.5k^{-1/3})$ and $\alpha=0.1$.}
		\label{fig: c_robust}
	\end{figure}

	In the first experiment, we tune the step length $\eta_k$ and $\delta_k$ jointly within the following form:
	\begin{equation}
		\eta_k=c\cdot 4k^{-1},\quad \text{and} \quad \delta_k=\min(0.1,c\cdot 0.5k^{-1/3}).
		\label{eq: special form of tuning}
	\end{equation}
	To understand the rationale of this form, note that these parameters give critical control to the variance of the gradient estimator. We aim to keep the variance of the term $\eta_k H_k$ at the same level in the gradient descent update
	\begin{align*}
		\bx_{k+1}=\Pi_\mathcal{B}(\bx_k-\eta_k \bH_k), \qquad \text{ with }\quad 
		\eta_k \bH_k=\eta_k\frac{\hat{f}(\bx_k+\delta_k/2\cdot \bZ_k)-\hat{f}(\bx_k-\delta_k/2\cdot\bZ_k)}{\delta_k}.
	\end{align*}
	In this experiment, we let $c\in\{0.6,1.0,1.2\}$ and fix $T_k=200k^{1/3}$ and $\alpha=0.1$. For each case, the regret curve is estimated by 100 independent runs for $L=1000$ iterations. 
	The regret and its linear fit are reported in  Figure \ref{fig: c_robust}. As shown in the right panel of Figure \ref{fig: c_robust}, the regret of LiQUAR has slopes of the linear regret fit close to 0.5 in all three cases.
	%Moreover, because the asymptotic order only reflects the long-term profit level, we are also interested in the short-term performance under different $c$. Therefore, we calculate the average unit time cost in different stages. In specific, we decompose $L=1000$ cycles into 4 quarters, with each quarter contains 250 cycles. Then the average regret incremental per unit time in each quarter is calculated by 
	%$$\text{Average incremental}=\frac{\text{Regret in each quarter}}{\text{Total time in each quarter}}.$$
	%This metric roughly measures the growth speed in different stages. The results are shown in the right panel of Figure \ref{fig: c_robust}.
	Comparing the two curves with $c=0.6$ and $c=1.2$ (left panel of Figure \ref{fig: c_robust}), we find that the larger value of $c$ immediately accumulates a large regret in the early stages but performs better in the later iterations. This observation may be explained by the trade-off between the level of exploration and learning rate of LiQUAR. In particular, a larger value of $c$ leads to larger values of $\eta_k$ and $\delta_k$, which allows more aggressive exploration and higher learning rate. %As a consequence, we observe a higher growth rate in regret at the early iterations due to more aggressive exploration, but a slower growth rate at the later iterations due to faster learning rate.    
	
	Although the tuning of $c$ will not affect the convergence of asymptotic regret of the algorithm, it may be critical to decision making in a finite-time period. %The above analysis provides some qualitative guidance on hyperparameter choice when applying GOLiQ-UDS in practice. Roughly speaking, 
	For example, a myopic decision maker who prefers good  system performance in a short term should consider small values of $c$, while a far-sighted decision maker who values more the long-term performance should adopt a larger $c$.% could be more helpful as it leads to better convergence to the optimal decision.
	
	%Therefore, to obtain a smaller regret in finite time, the service provider needs to balance between the regret order and the initialization cost. Smaller initialization cost may leads to larger regret order, which means more revenue loss in the long run. This also relates to the exploration and exploitation trade-off in learning problems. See the following Remark \ref{rmk: exp-exp trade-off}.
	%Discussion
	%\begin{remark}[The management insight of length $c$]
	%	\label{rmk: exp-exp trade-off}
	%	To some extent, different $c$ reflects different level of exploration. A myopic service provider, who focuses more on exploitation and less on exploration, would prefer smaller $c$, for example, $c=0.6$; on the contrary, larger $c$ means that GOLiQ-UDS tends to explore more with large $\delta_k$ and $\eta_k$. In this sense, the service provider needs to balance the exploration and exploitation to earn more profit. If GOLiQ-UDS explores too much, it is more likely to find the optimal policy rapidly but more exploration cost needs to be paid, which is reflected by the high regret in the beginning stage for the curve $c=1.2$ and smaller asymptotic order of the regret. However, GOLiQ-UDS needs sufficient exploration to learn well as many learning algorithms do. So a too conservative service provider who chooses too small $c$ incurs high regret order. As long as the exploration is sufficient, i.e. $c$ is not too small, GOLiQ-UDS performs well and converges to the optimal decision rapidly as our theory suggests.
	%\end{remark}
	\subsubsection{Cycle length $T_k$.}\label{subsubsec: test T} 
\begin{figure}[h]
		\centering
		\includegraphics[width=\linewidth]{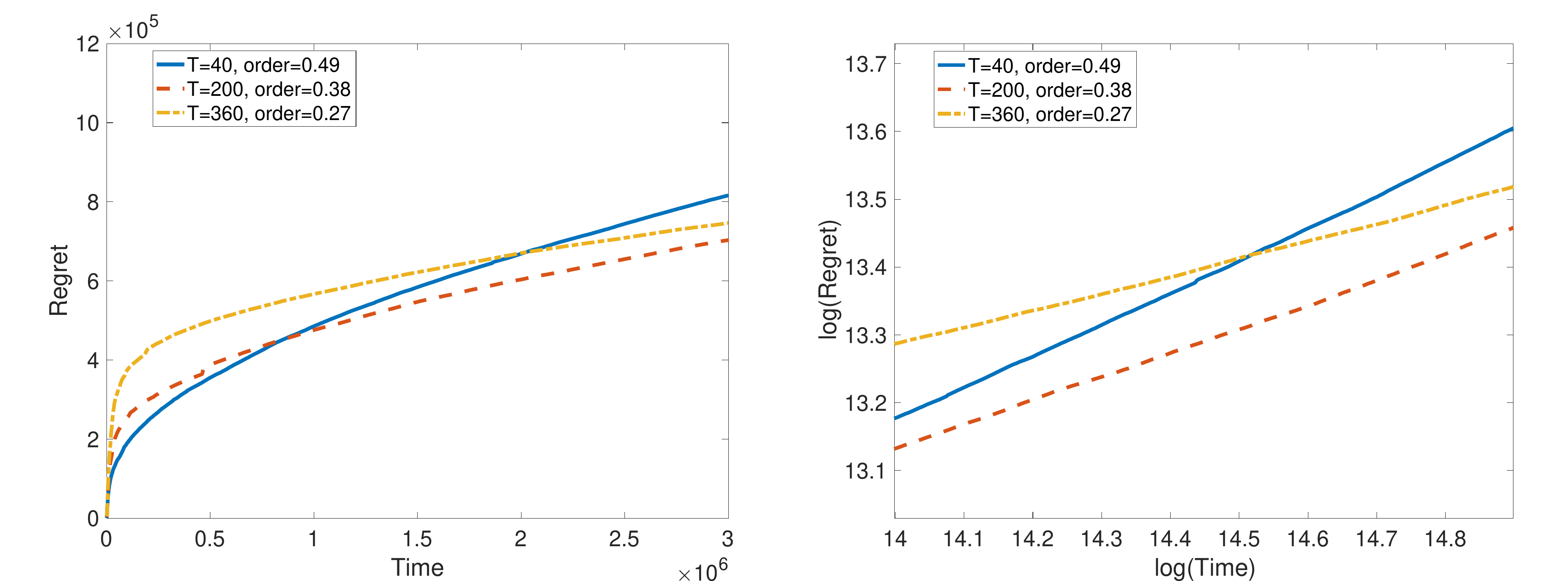}
		\caption{Regret under different $T\in\{40,200,360\}$: (i) average regret from 100 independent runs (left panel); (ii) regret curve in logarithm scale, with $T_k=T\cdot k^{1/3}$, $\eta_k=4k^{-1}, \delta_k=\min(0.1,0.5k^{-1/3})$ and $\alpha=0.1$.}
		\label{fig: test T}
	\end{figure}

In this experiment, we test the impact of $T_k$ on the performance of LiQUAR. We again use the $M/M/1$ base example. The step-length hyperparameters are set to
	$\eta_k=4k^{-1}$ and $\delta_k=\min(0.1, 0.5k^{-1/3})$.
	We choose different values of $T_k$ in the form of
	$$T_k=T\cdot k^{1/3},~ T\in\{40,200,360\}.$$
	For different values of $T$, iteration numbers $L_T$ are chosen to maintain equal total running times for LiQUAR. In particular, we choose 
	$L_T=\left\lceil1000\cdot \left(200/T\right)^{3/4}\right\rceil.$
	%For each $T$, the corresponding regret curve is estimated by 100 independent runs of GOLiQ-UDS, and the regret order is estimated by a linear fit of the log-log curve. 
	Results of all above-mentioned cases are reported in Figure \ref{fig: test T}.
	
	The right panel of Figure \ref{fig: test T} shows that the slope of the linear fits all below 0.5.
	According to the three regret curves in the left panel, we can see how different values of $T$ impact the exploration-exploitation trade-off: a larger value of $T$, e.g., $T=360$, yields a higher regret in the early iterations but ensures a flatter curve in the later iterations.
	This is essentially due to the trade-off between the learning cost and the quality of the gradient estimator. A larger cycle length $T_k$ guarantees a high-quality gradient estimators as more data are generated and used in each iteration which help reduce the gradient estimator's transient bias and variance. On the other hand, it demands that the system be operated for a longer time  under suboptimal control policies, especially in the early iterations. 
	The above analysis provides the following guidance for choosing $T$ in practice: A smaller $T$ is preferred if the service provider's goal is to make the most efficient use of the data in order to make timely adjustment on the control policy. This guarantees good performance in short term (the philosophy here is similar to that of the temporal-difference method with a small updating cycle). But if the decision maker is more patient and aims for good long-term performance, he/she should select a larger  $T$ which ensures that the decision update is indeed meaningful with sufficient data (this idea is similar to the Monte-Carlo method with batch updates).
	\subsection{Queues with non-Poisson arrivals}
	%\subsection{Queues with non-Poisson arrivals}
	\label{subsec: GG1}
	In this section, we consider the more general $GI/GI/1$ model having arrivals according to a renewal process.
	Similar to the service times, we model the interarrival times using scaled random variables $U_1/\lambda(p),U_2/\lambda(p), \ldots$ for a given $p$, with $U_1,U_2,\ldots$ being a sequence of I.I.D. random variables with $\mathbb{E}[U_n]=1$. 
	
	The PTO framework is not applicable here because $\EE[W_\infty]$ does not have a closed-form solution in the $GI/GI/1$ setting. This provides additional motivations for our online learning approach. 
	On the other hand, generalizing the theoretical regret analysis rigorously from $M/GI/1$ to  $GI/GI/1$ is by no means a straightforward extension. 
	A key step in our analysis is to give a proper bound for the bias of the gradient estimator. When the arrival process is Poisson, the memoryless property ensures that $N_l/T_k$ in \eqref{eq: f G hat} is an unbiased estimator for the arrival rate. 
	%in the $GI$ arrival case, this estimator has a bias which will largely amplify the bias of the gradient estimator. 
	For renewal arrivals, the arrival rate bias has an order $O(1/T_k)=O(k^{-1/3})$ (see for example Lorden's inequality \citep[Section V, Proposition 6.2]{AsmussenBook}), which contributes to the bias of the FD with an order of $O(1/T_k\delta_k)=O(1)$. This contradicts  Theorem \ref{thm: convergence} which requires $B_k=O(k^{-1})$. 
	This part of the analysis requires additional investigations (in order to establish a more delicate bias bound). %Unfortunately, such inference on the arrival bias is not straightforward and usually needs careful construction of  couplings, e.g., \cite{lindvall1986coupling}, \cite{hermann2000coupling}}. 
We leave the careful regret analysis of $GI/GI/1$ to future research.

{Nevertheless, from the engineering perspective, the increased bias due to the $GI$ arrival process may not be too significant (note that the theoretical bias bound is obtained from a worst-case analysis). %Therefore, LiQUAR is still applicable when the arrival is $GI$ because LiQUAR is model-free and does not draw the Poisson arrival assumption heavily. 
	We next conduct some preliminary numerical experiments to test the performance of LiQUAR under $GI$ arrivals.} 	
We consider an $E_2/M/1$ queue example having Erlang-2 interarrival times with mean $1/\lambda(p)$ and exponential service times with rate $\mu$ to illustrate the performance of LiQUAR in $GI/GI/1$'s case. We continue to consider the logit demand function \eqref{eq: logit demand} with $M=10,a=4.1,b=1$ and linear staffing cost function \eqref{eq: cost function}. Unlike the $M/GI/1$ case where the PK formula provides a closed-form formula for the steady-state waiting time, here we numerically compute the optimal solution $(\mu^*,p^*)$ by using matrix geometric method (note that the state process of $E_2/M/1$ is quasi-birth-and-death process). Letting $h_0=c_0=1$ yields the optimal decision  $(\mu^*,p^*)=(7.78,3.75)$.
%the stationary distribution of queue length for an $E_2/M/1$ queue is a geometric distribution with failure rate $\sigma$, where $\sigma$ is the solution of 
%$$\log(\sigma)=\phi_A(\mu(\sigma-1))$$
%and $\phi_A$ is the c.g.f. of the Erlang-2 distribution with mean $1/\lambda(p)$. In this way, we could numerically obtain the exact optimal solution $(\mu^*,p^*)$ and the maximum profit $\mathcal{P}(\mu^*,p^*)$ as well. 
\begin{figure}[h]
	\centering
	\vspace{-0.2in}
	\includegraphics[width=.98\linewidth]{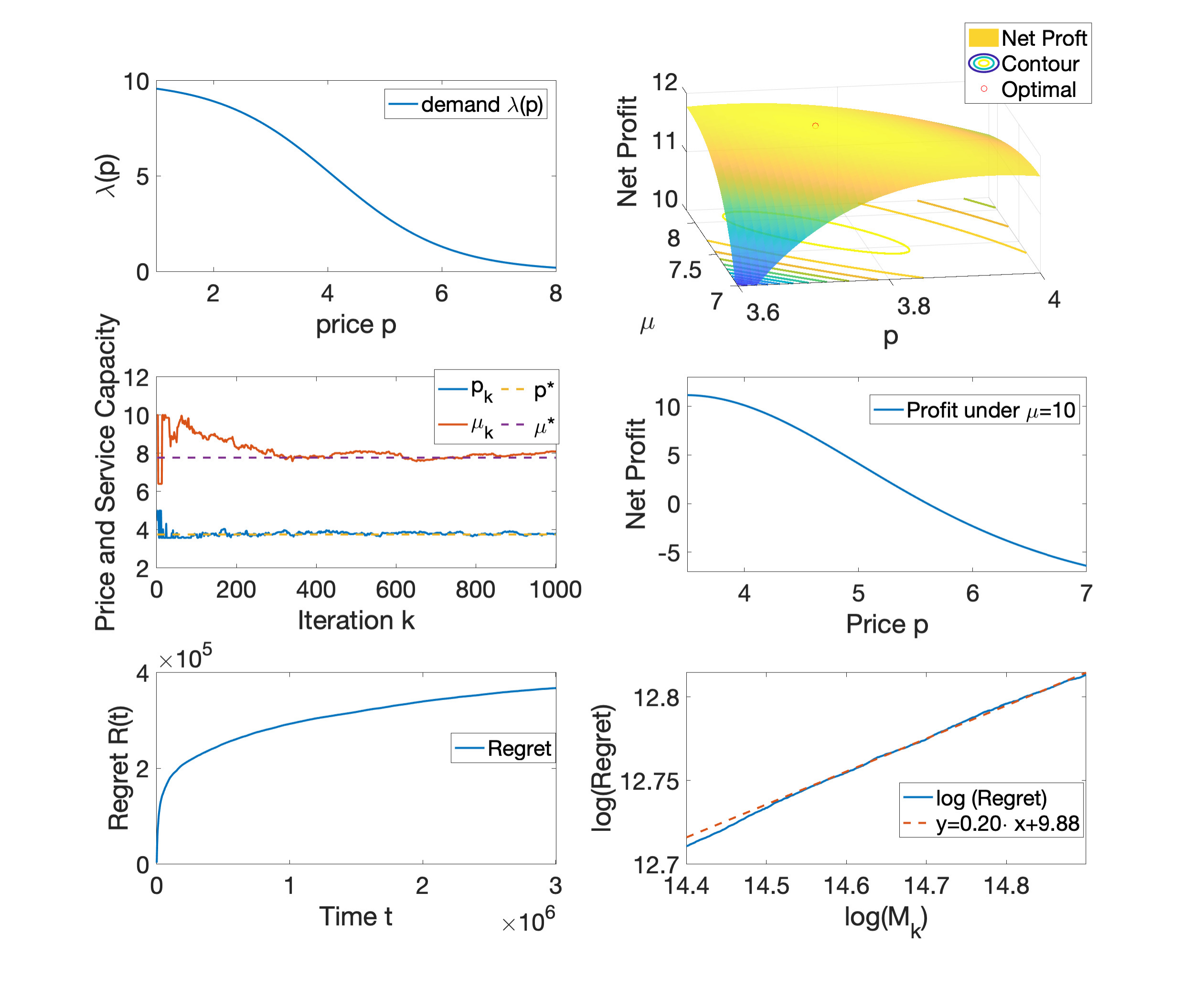}
	\vspace{-0.2in}
	\caption{Joint pricing and staffing in the $E_2/M/1$ queue with $\eta_k=4k^{-1},\delta_k=\min(0.1,0.5k^{-1/3})$, $T_k=200k^{1/3}$, $\alpha=0.1$, $p_0=5$ and $\mu_0=10$.% (i) Demand function $\lambda(p)$ (top left panel); (ii) net profit function and corresponding contours (top right panel); (iii) sample path of decision parameters (middle left); (iv) sample path of gradient (middle right); (iv) estimated regret curve (bottom left); (iv) logarithm of regret versus logarithm of total time $M_k$ plot with linear fit. The regret is estimated by 100 independent Monte-Carlos runs.
	}
	\label{fig: E2M1}
\end{figure}
We implement LiQUAR with hyperparameters $\eta_k=4k^{-1},\delta_k=\min(0.1,0.5k^{-1/3})$ , $T_k=200k^{1/3}$, and $\alpha=0.1$. Figure \ref{fig: E2M1}, as an analog to Figure \ref{fig: basic plot}, shows that the refined LiQUAR continues to be effective, exhibiting a rapid converge to the optimal decision and a slowly growing regret curve (bottom left panel of Figure \ref{fig: E2M1}). %with asymptotic order $0.18$ (bottom right panel of Figure \ref{fig: E2M1}). %In this sense, GOLiQ-UDS continues to perform well in non-Markovian arrival settings. We believe that the effectiveness would preserve in cases of more $GI/GI/1$ queues due to the model-free design GOLiQ-UDS.
Despite of the good performance of the above $E_2/M/1$ example, we acknowledge that this is only a preliminary step, and the full investigation of the $GI/GI/1$ case requires careful theoretical analysis and comprehensive numerical studies. 	
%%%%%%%%%%%%%%%%%%%%%%%%%%%%%%%%%%%%%%%%%%%%%%%%
%%%%%%%%%%%%%%%%%%%%%%%%%%%%%%%%%%%%%%%%%%%%%%%%

\section{LiQUAR vs. PTO}\label{sec:PTO_compare}
In this section, we contrast the performance of LiQUAR to that of the conventional PTO method. In principle, a PTO algorithm undergoes two phases: (i) ``prediction" of the model (e.g., estimation of the demand function and service distribution) and (ii) ``optimization" of the decision variables (e.g., setting the optimal service price and capacity). 
Taking the demand function $\lambda(\cdot)$ as an example, PTO relies on the ``prediction" phase to provide a good estimate $\widehat{\lambda}(p)$, which will next be fed to the ``optimization'' phase for generating desired control policies. In case no historical data is available so that the ``prediction'' completely relies on the newly generated data, one needs to learn the unknown demand curve $\lambda(p)$ by significantly experimenting the decision parameters in real time in order to generate sufficient demand data that can be used to obtain an accurate $\widehat{\lambda}(p)$. 

%We compare LiQUAR with PTO in this section both theoretically and numerically. 
We begin by establishing theoretical results to compare the performance of LiQUAR and PTO within a heavy-traffic framework. These findings are then supplemented by numerical experiments to give engineering confirmations. %; We also consider an advanced PTO aided by the  downstream optimization structure.

\subsection{LiQUAR vs. PTO in heavy traffic}\label{sec: PTO_compare_HT}
\begin{figure}
	\centering
	\includegraphics[width=\linewidth]{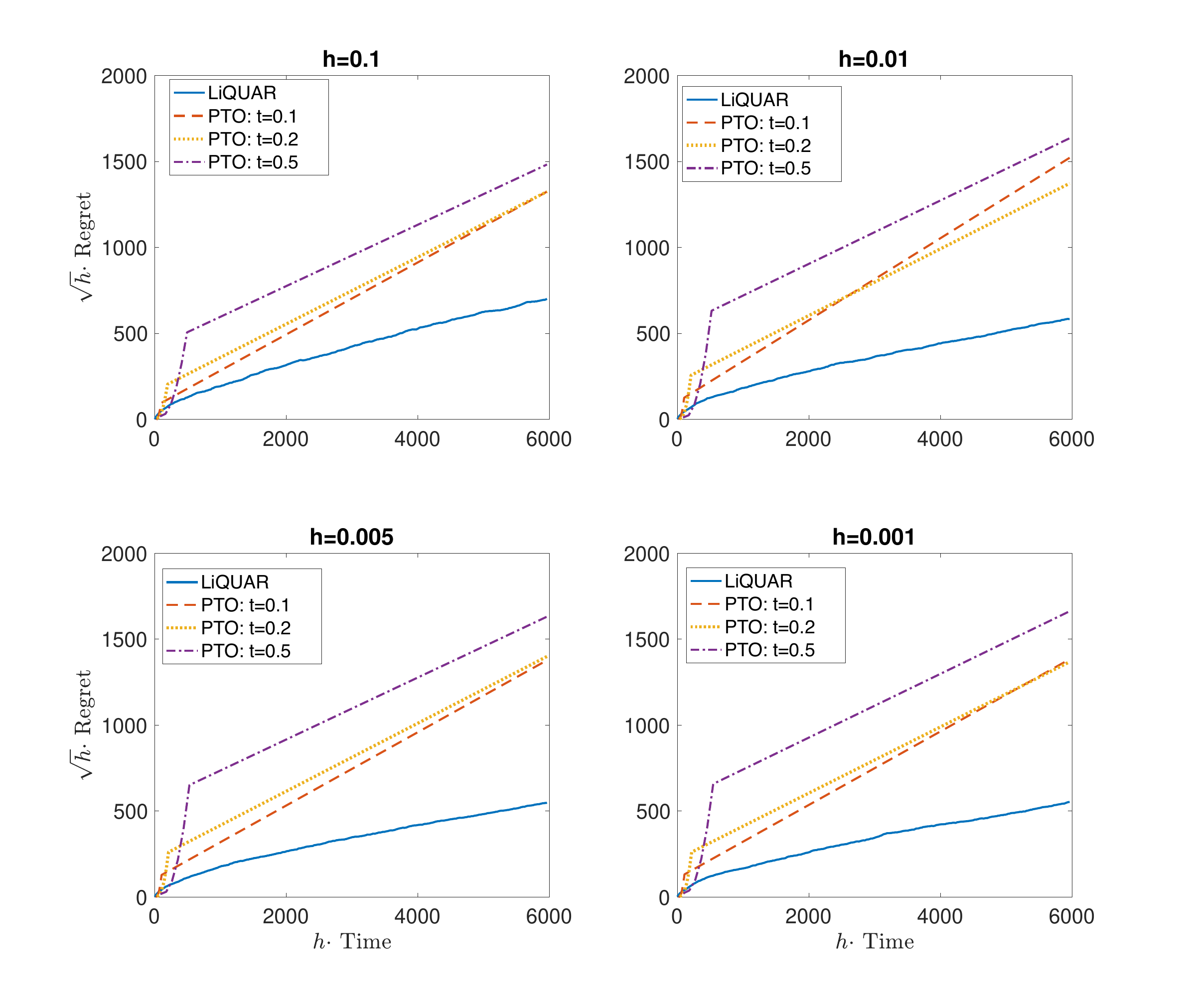}
	\caption{PTO vs. LiQUAR; $\eta^h_k=4\sqrt{h}k^{-1},\delta^h_k=2\sqrt{h}k^{-1/3}$ with $T^h_k=h^{-1}k^{1/3}$; $t_0^h=t\cdot \frac{T_0^{5/7}\log (T_0/h)^{2/7}}{\sqrt{h}}$ and $\kappa= t_0^{1/5}\cdot h^{1/5}$ with $t\in\{0.1,0.2,0.5\}$}
	\label{fig: NPPTO}
\end{figure}

In this section, we compare the regret bounds of LiQUAR and PTO when the system is in heavy traffic. Consider the sequence of $h$-indexed systems described in Section 5.3, we now intend to use PTO to find the optimal value of \eqref{eq: objct heavy traffic} within $\mathcal{B}_h$ and measure its performance by computing the regret $R^h(T_0)$ at time $T^h=T_0/h$. We consider the following PTO algorithm \citep{besbes2009dynamic}:
\begin{itemize}
	\item \textbf{Input:} Total running time $T^h$, number of testing points $\kappa^h$, time of prediction $t_0^h$
	\item \textbf{Step 1. Prediction:}
	\begin{enumerate}
		\item[a.] Organize $\mathcal{B}_h$ into $\kappa^h$ evenly spaced grids and distribute the testing points in all grids.
		\item[b.] For each $i=1,\cdots,\kappa^h$, operate the system at $i^{\rm th}$ under the testing point $p_i$ for $t_0^h/\kappa^h$ units of time, and approximate the demand curve $\hat{\lambda}(p_i)$ by the time-averaged arrival rate.
		%\item Calculate $\hat{f}(p_i)$ by plugging in $\hat{\lambda}(p_i)$
	\end{enumerate}
	\item \textbf{Step 2. Optimization:} 
	\begin{enumerate}
		\item[a.] Calculate $\hat{f}(p_i)$ using $\hat{\lambda}(p_i)$ in the PK formula.
		\item[b.] Operate the system under $\hat{p}^*=\arg\max_{i} \hat{f}(p_i)$ for the rest of time horizon.
	\end{enumerate}
	
\end{itemize}
%The essential error of this method comes from three sources:
%\begin{enumerate}
%	\item Price experimentation cost $O(\Delta_0^h t_0^h)=O(\sqrt{h} t_0^h)$
%	\item Stochastic Fluctuation Error $|f(p^*)-f(p_G^*)|=O(1/\sqrt{t_0^h/\kappa ^h})=O(\sqrt{\kappa^h/t_0^h})$
%	\item Deterministic Grid Error $|p^*-p_G^*|=O(|\mathcal{B}_h|/\kappa^h)$,
%\end{enumerate}
%where $\Delta_0^h=\max_{p\in\mathcal{B}_h} f^h(p) - f^h(p^*)$, $p_G^*$ is the optimal points in the testing pricing grid. We neglect the superscribe $h$ in the following analysis to ease the burden of notation.
%Following the analysis in \cite{besbes2009dynamic}, we have
%\begin{align*}
%	R(t_0,\kappa,h,T)&\leq \Delta_0 t_0+ (T-t_0) \EE[f(\hat{p}^*)-f(p^*)]\\
%	&\leq \Delta_0 t_0+T\cdot [f(\hat{p}^*)-f(p^*_G)+f(p^*_G)-f(p^*)]\\
%	&\leq \Delta_0 t_0 + T\cdot \sqrt{\frac{\kappa \log T}{t_0}} + CT\cdot \frac{\nabla^2 f(\xi) }{2} |p_G^*-p^*|^2\\
%	&\leq C\sqrt{h} t_0 + CT\cdot \sqrt{\frac{\kappa \log T}{t_0}}  +C\frac{T}{\sqrt{h}}\cdot \left( \frac{h}{\kappa}\right)^2
%\end{align*}
%For this term, we find that $t_0=O(\frac{T^{5/7}\log (T) ^{2/7}}{h^{2/7}})=O(\frac{T_0^{5/7}\log (T)^{2/7}}{h})$ and $\kappa=O(h\cdot (t_0/\log T)^{1/5})=O\left(\frac{hT}{\log T}\right)^{1/7}=O(\frac{T_0^{1/7}}{\log T})$ are optimal. Therefore, we would have the optimal regret order that
%$$R(t_0^*,\kappa^*,h,T)=O\left(\frac{T_0^{5/7}\log (T_0/h)^{2/7}}{\sqrt{h}}\right).$$
%Following the analysis in \cite{besbes2009dynamic}, 
We next give a regret bound for the above PTO method under the heavy-traffic learning scheme. %and we summarize it in the following Proposition \ref{prop: PTO heavy-traffic}.
\begin{proposition}[PTO in heavy traffic]\label{prop: PTO heavy-traffic}
	Under Assumption \ref{ass: demand heavy traffic}, in the $h^{\rm th}$ system, PTO with hyperparameters $\kappa^h,t_0^h$ yields the regret bound: 
	\begin{align}\label{PTO_HT}
		R^h(T_0)\leq C\sqrt{h}t_0^h+CT_0\cdot\frac{\sqrt{\kappa^h\log T_0/h}}{h\sqrt{t_0^h}}+C\frac{\sqrt{h}T_0}{(\kappa^h)^2}.
	\end{align}
	In addition, if we select $t_0^h=O\left(\frac{T_0^{5/7}(\log T^h)^{2/7}}{h}\right)$ and $\kappa =O\left(\frac{T_0^{1/7}}{\log T^h}\right)$, then the PTO regret bound can be minimized as below: 
	$$R^h(T_0)\leq C\cdot \frac{T_0^{5/7}\log (T_0/h)^{2/7}}{\sqrt{h}}=\tilde{O}\left(\frac{T_0^{5/7}}{1-\rho_h^*}\right),$$
	with $C$ being some constant independent with $h$ and $1-\rho_h^*$.
\end{proposition}
\begin{remark}[LiQUAR vs. PTO in heavy traffic]\label{rem:LiQUARvsPTO}
	Compared to the original regret bounds presented in \cite{besbes2009dynamic}, we have re-optimized the hyperparameters and derived reduced regret bounds under Assumption \ref{assmpt: uniform}, enabling a fair comparison between LiQUAR and PTO. According to Theorem \ref{thm: heavy traffic} and Proposition \ref{prop: PTO heavy-traffic}, the regret bounds for both LiQUAR and PTO share a dependence on the traffic intensity in the order of $1/(1-\rho_h^*)$. However, LiQUAR exhibits a slower growth rate with respect to the time horizon, scaling as $\sqrt{T_0}$. This implies that over the long run, LiQUAR achieves a smaller regret bound than PTO.
	Furthermore, as the two terms involving $T_0$ and $\rho^*$ interact multiplicatively in the regret bound, the factor $1/(1-\rho^*)$ amplifies LiQUAR's advantage over PTO in heavy-traffic conditions, i.e., as $\rho^* \to 1$. This trend is further validated by the numerical results illustrated in Figure \ref{fig: NPPTO}.
\end{remark}

Next, we numerically investigate the performance of LiQUAR and PTO for a one-dimensional pricing problem under heavy traffic. %Following Assumption \ref{ass: demand heavy traffic}, w
We consider an $M/M/1$ system having exponential demand function
$$\lambda(p)=\exp(a-bp),$$
with $a=1+\log 2$, $b=1$ and exponential service time distributions. Following the settings in Theorem \ref{thm: heavy traffic}, we consider a sequence of objective functions in \eqref{eq: base obj} indexed by $h$.
We keep $\mu=1$ held fixed and allow $h\in\{0.1,0.01,0.005,0.001\}$ to account for different values of  the traffic intensity. As $h\rightarrow0$, the feasible region $\mathcal{B}_h$ takes the form  
$$ \mathcal{B}_h=\left[~p_0+ c_1\sqrt{h}, p_0+c_2\sqrt{h}~\right],$$
with $c_1=0.6\cdot c_0,c_2=5c_0$ and $c_0=1/\sqrt{\lambda'(p_0)(1+p^*\lambda'(p_0))}=1.20$.

Then, we apply LiQUAR and PTO in all the instances with different $h$. % without using specific knowledge of exponential demand function $\lambda(p)$ beyond $\mathcal{B}_h$. 
For LiQUAR, following  Theorem \ref{thm: heavy traffic}, we choose $\eta^h_k=4\sqrt{h}k^{-1},\delta^h_k=2\sqrt{h}k^{-1/3}$ with $T^h_k=h^{-1}k^{1/3}$ for 500 iterations. To make a fair comparison, we pick an equal runtime for LiQUAR and PTO  with $T_0^h=\sum_{k=1}^{500} T_k^h$ and $T_0=h\cdot T_0^h$ for all $h$. %Following the analysis of PTO, the 
PTO's hyperparameters are chosen as $t_0^h=t\cdot \frac{T_0^{5/7}\log (T_0/h)^{2/7}}{\sqrt{h}}$ and $\kappa^h= h\cdot (t_0/\log T^h)^{1/5} $ with $t\in\{0.1,0.2,0.5\}$.

%\paragraph{Experiment Results}  
In Figure \ref{fig: NPPTO}, we report the scaled regret curves of both methods under different holding costs $h$ where each regret curve is estimated by the average of 100 independent replications. To understand how the regret is influenced by the heavy-traffic scaling factor $h$, we scale time by $h\cdot T^h=T_0$ and scale the regret by $\sqrt{h}\cdot R^h(T_0)$ (as in Theorem \ref{thm: heavy traffic} and Proposition \ref{prop: PTO heavy-traffic}). 
From Figure \ref{fig: NPPTO}, we confirm that LiQUAR significantly outperforms PTO in all heavy-traffic scenarios. %In addition, for each method, the scaled regret curves have roughly the same regret, which implies that for both PTO and LiQUAR, they have a $1/\sqrt{h}$ dependence on $h$. This result is consistent with our analysis in Theorem \ref{thm: heavy traffic} and the analysis on PTO.
%\begin{figure}
%	\centering
%	\includegraphics[width=\linewidth]{Heavy-traf-eps-converted-to.pdf}
%	\caption{Comparison between NPPTO and LiQUAR; $\eta^h_k=4\sqrt{h}k^{-1},\delta^h_k=2\sqrt{h}k^{-1/3}$ with $T^h_k=h^{-1}k^{1/3}$; $t_0^h=t\cdot \frac{T_0^{5/7}\log (T_0/h)^{2/7}}{\sqrt{h}}$ and $\kappa= t_0^{1/5}\cdot h^{1/5}$ with $t\in\{0.1,0.2,0.5\}$}
%\end{figure}

% \subsection{Additional Numerical Experiments}
% In this section, we provide simulation experiments that compare LiQUAR with other learning methods in various settings numerically. In detail, we first extend our numerical experiment in Section 6.3 by comparing LiQUAR with pPTO and Smart-pPTO method in Section \ref{subsec: SPTO}. Then, we focus on the numerical comparison between LiQUAR and non-parameteric PTO (NPPTO)  under heavy-traffic conditions in Section xxx. Finally, we compare LiQUAR with Reinforcement Learning methods, specifically, Policy Gradient (PG), numerically and investigate how our analysis could potentially accelerate PG in our queueing setting in Section xxx.
\subsection{LiQUAR vs. objective-informed PTO} \label{subsec: SPTO}

In this section, we compare LiQUAR to an advanced parametric PTO framework, where the prediction step incorporates information from the downstream objective function. We refer to this approach as \textit{objective-informed PTO} (oiPTO). In oiPTO, it is assumed that the decision-maker knows the parametric form of the demand function, $\lambda(\cdot; \boldsymbol{\beta})$, with parameters $\boldsymbol{\beta}$ that are initially unknown.
During the prediction phase, oiPTO estimates $\boldsymbol{\beta}$ by leveraging the structure of the downstream objective function (see \eqref{eq: SpPTO obj}). In the optimization phase, oiPTO uses the demand function with the estimated parameters, denoted as $\boldsymbol{\beta}_{oi}$, to compute the optimal decisions.

Specifically, let $\theta \in (0,1)$ represent the exploration ratio and $T$ denote the total time (or learning budget). The oiPTO approach divides the total time horizon $T$ into two phases. In the first phase, corresponding to the interval $[0, \theta T]$, the focus is on learning the parameters of the demand function. In the second phase, covering the interval $[\theta T, T]$, the system operates using decisions optimized based on the estimated parameters. The details of these two steps are provided below:

%In this section, we extend our numerical analysis in Section 6.3 to contrast the performance of LiQUAR to that of the Smart-parametric-Predict-Then-Optimization (SpPTO). In principle, the SpPTO is similar to pPTO with prediction and optimization phases. The critical example between SpPTO and pPTO lies in its prediction objective. We give the details as below:
\begin{itemize}
	\item \textbf{Prediction:} Suppose $m$ parameters of the demand function are to estimated, and we uniformly select $(p_1,\mu_1),\cdots,(p_m,\mu_m)\in\mathcal{B}$ as experimentation decisions. We sequentially operate the system under each of the experimentation decision for $\theta T/m$ units of time. Then, based on the arrival and workload data, the estimated parameter is given by
	%		\begin{equation}\label{eq: pPTO obj}
		%		\boldsymbol{\beta}^*=\argmin_{\boldsymbol{\beta}}\sum_{k=1}^{m}(\lambda(p_k;\boldsymbol{\beta})-N_km/(\theta T))^2,\quad \quad \quad \text{for pPTO}
		%	\end{equation}
	\begin{equation}\label{eq: SpPTO obj}
		\boldsymbol{\beta}_{oi}=\arg\min_{\boldsymbol{\beta}} \sum_{k=1}^{m}(f(p_k,\mu_k;\boldsymbol{\beta})-\hat{f}(p_k,\mu_k))^2,
	\end{equation}
	where   $f(p,\mu;\boldsymbol{\beta})=-p\lambda(p;\boldsymbol{\beta})+h\cdot \frac{\lambda(p;\boldsymbol{\beta})}{\mu-\lambda(p;\boldsymbol{\beta})}+c(\mu)$  and $\hat{f}(p,\mu)$ is the time average cost estimation of $f(p,\mu)$.
	\item \textbf{Optimization:} Next, we obtain the oiPTO-optimal policy $\hat{x}^*$ by maximizing our objective function with $\lambda(\cdot;\boldsymbol{\beta})$ replaced by $\lambda(\cdot;\boldsymbol{\beta}_{oi})$. Then we implement this policy for the rest of time horizon.
\end{itemize}
\paragraph{Experiment settings and results.} We consider our base logit example in Section \ref{subsec: basic exp} having demand function \eqref{eq: logit demand} with $M_0=10,a=4.1,b=1$ and exponential service times. Throughout this experiment, we fix the staffing cost $c(\mu)=\mu$. 
To understand the impact of the system's congestion level on performance of oiPTO and LiQUAR, we consider two scenarios specified by the optimal traffic intensity $\rho^*$: (i) A light-traffic case with $\rho^* = 0.709$ ($h_0=1$) and  (ii) A heavy-traffic case with $\rho^* = 0.987$ ($h=0.001$).  %In both scenarios, we compare pPTO and LiQUAR with the same uniformly stable region $\mathcal{B}=[6.18,10]\times[3.62,7]$. 

%As a comparison, we also consider a base parametric PTO method dubbed pPTO.The only difference of pPTO from oiPTO is that in prediction, the parameter is estimated only based on the empirical arrival rates, i.e., 
%$$\boldsymbol{\beta}^*_{p}=\arg\min_{\boldsymbol{\beta}}\sum_{k=1}^{m}(m\cdot N_k/\theta T-\lambda(p_k;\boldsymbol{\beta}))^2.$$

For LiQUAR, we consistently select the hyperparameters $\eta_k=4k^{-1},\delta_k=\min(0.1,0.5k^{-1/3})$, initial values $(\mu_0,p_0)=(10,7)$ and $T_k=200k^{1/3}$ for $L=1000$ iterations with a total running time $T=2\sum_{k=1}^{L}200k^{1/3}$.		
For oiPTO, we use the same total time $T$ 
and consider several values of the exploration ratio $\theta\in\{0.3\%,0.9\%,1.5\%,6\%,15\%\}$ to account for different levels of exploration efforts. 

\begin{figure}[htp]
	
	\subfloat[Light traffic case with $\rho^*=0.706$]{%
		\includegraphics[clip,width=\columnwidth]{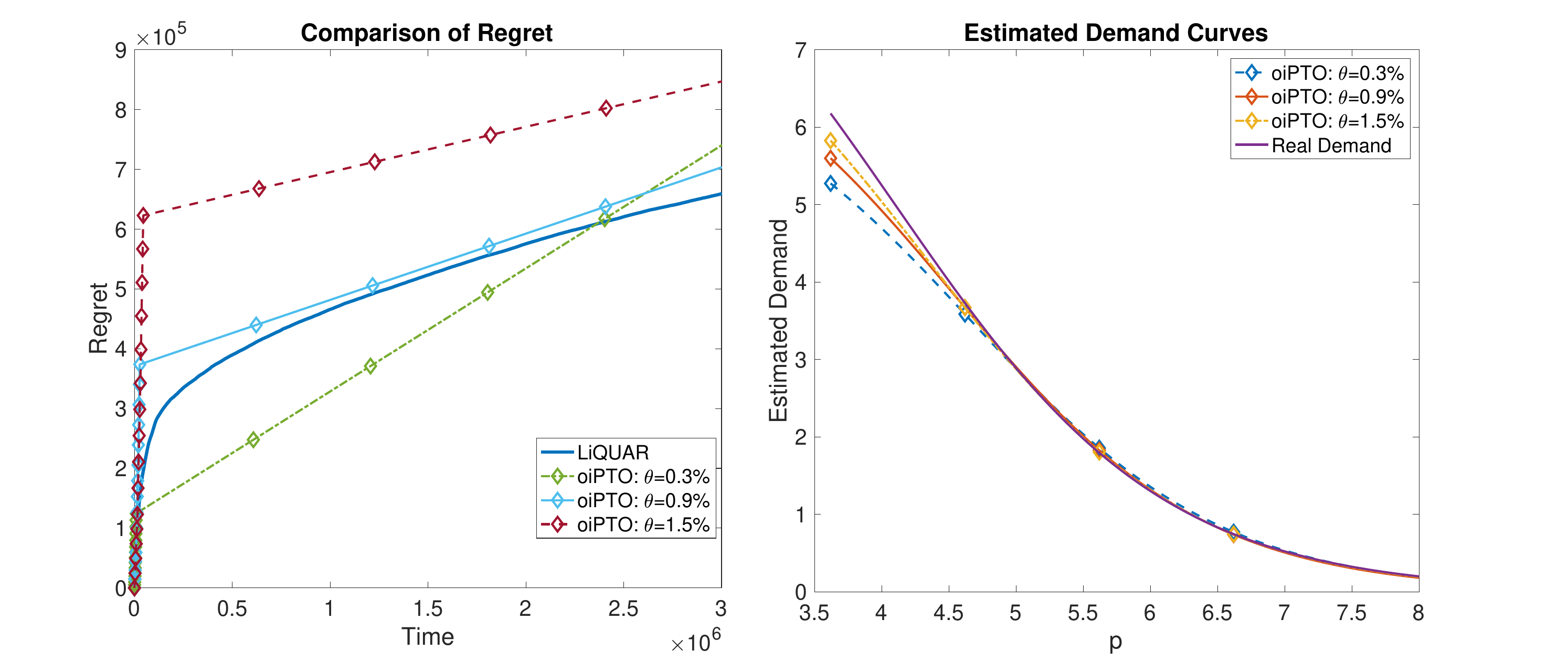}%
	}
	
	\subfloat[Heavy traffic case with $\rho^*=0.987$]{%
		\includegraphics[clip,width=\columnwidth]{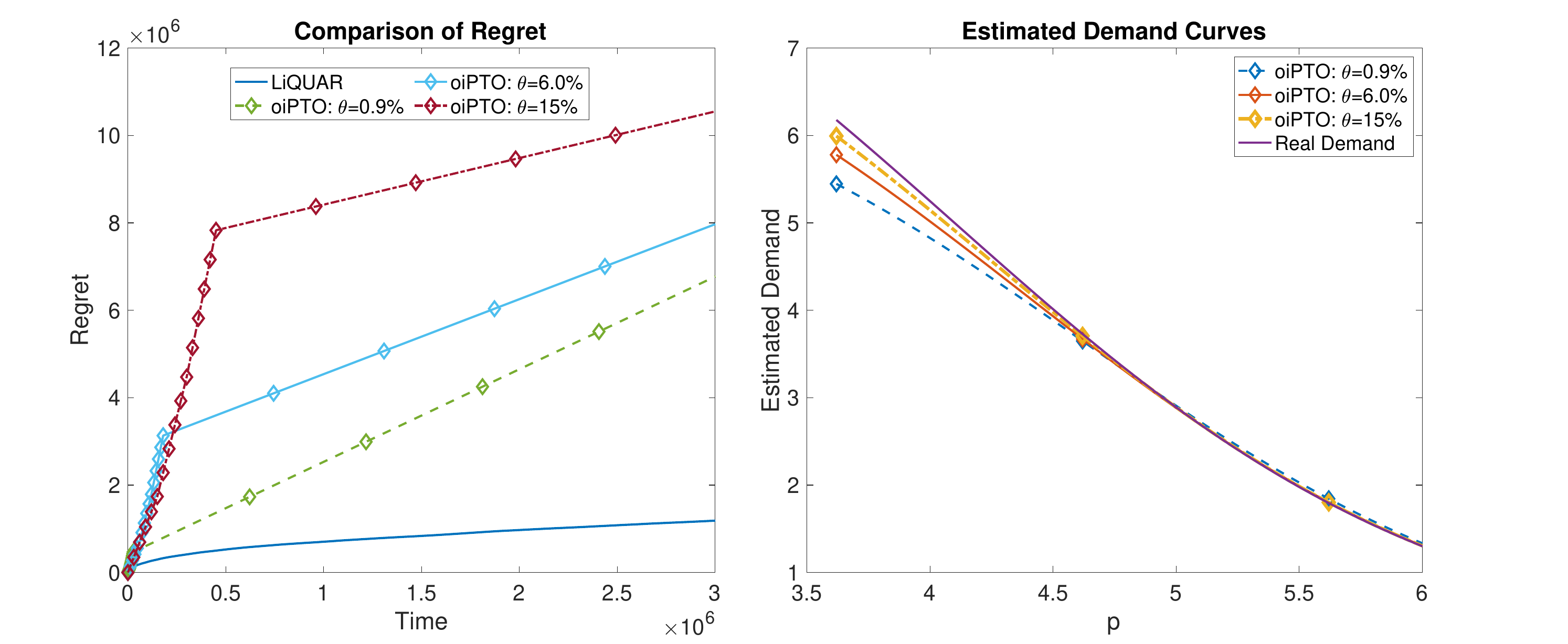}%
	}
	\vspace{0.1in}
	\caption{oiPTO vs. LiQUAR: (i) low traffic scenario $\rho^*=0.705$; (ii) high traffic scenario $\rho^*=0.987$.  Hyperparameters for LiQUAR are $\eta_k=4k^{-1},\delta_k=\min(0.1,0.5k^{-1/3}),T_k=200k^{1/3}$ in both scenarios. All regret curves are estimated by averaging 1,000 independent simulation runs.}
	\label{fig: SpPTO}
\end{figure}

%\paragraph{Experiment Results.}
In Figure \ref{fig: SpPTO}, we present the regret results for LiQUAR and oiPTO, showcasing the three oiPTO curves with the lowest regrets. The left-hand panels illustrate that the exploration ratio $\theta$ has a significant impact on oiPTO's performance. The regret for oiPTO exhibits a piecewise linear pattern: during the prediction phase, regret grows rapidly due to periodic exploration across all experimentation variables; in the optimization phase, regret continues to increase linearly, albeit at a slower rate, as the system operates based on the oiPTO-optimized solution, which remains suboptimal.
A larger (smaller) $\theta$ leads to higher (lower) regret during the prediction phase but results in a more (less) accurate model. This improved accuracy generates decisions that are closer to optimal, resulting in a slower (faster) regret growth during the optimization phase.

Furthermore, comparing case (a) to case (b), we observe that while oiPTO incorporates downstream objective information, LiQUAR consistently outperforms oiPTO by achieving a lower regret. This advantage is especially pronounced in heavy-traffic scenarios. The primary reason is that LiQUAR employs an integrated learning approach, continuously refining its decision-making through direct interaction with the environment. In contrast, oiPTO follows a static learning strategy, investing fixed efforts into parameter prediction and relying entirely on these predictions in the optimization phase. The limitations of oiPTO become more apparent in heavy traffic, where the nonlinear structure of workload amplifies the cost of suboptimal decisions.%With a short exploration time, the prediction is inaccurate and the cost of the inaccuracy of the final decision prevails, leading to a regret slope and a high relative cost (see the regret and estimated demand curve of $\theta=0.3\%$ in panel (a)). On the other hand, if pPTO spends too much time on the ``prediction" stage, the initial exploration cost overshadows the cost of the final decision (see the $\theta=1.5\%$ line in panel (a)). The service provider needs to take a careful balance between the intrinsic exploration-exploitation trade-off in pPTO. 

\section{LiQUAR vs. Reinforcement Learning}
\label{sec: RL}
In this section, we compare LiQUAR with reinforcement learning (RL) methods. While the machine learning literature offers a wide range of RL approaches, we focus on the policy gradient (PG) method for comparison due to the following reasons: (i) both LiQUAR and PG rely on gradient-based optimization, making them conceptually aligned; and (ii) our problem involves an infinite state space (e.g., queue length or workload) and a continuous action space, where the PG method demonstrates particular advantages over other RL techniques. 
%We first introduce the problem setting in MDP and then introduce the benchmark policy gradient method in our specific problem.

\paragraph{Problem Settings and Algorithms.} Because an RL method is underpinned by its corresponding Markov decision process (MDP), and setting up an MDP requires the model to be Markovian, we now restrict our attention to the $M/M/1$ queue. Specifically, we consider the following discrete-time MDP with the objective of maximizing its long-run average reward. Our MDP has
\begin{itemize}
	\item \textit{Time steps:} $t=1,2,\dots$.
	\item \textit{States:} Queue length at beginning of period $t$, denoted by $S_t$.
	\item \textit{Actions:} Choices of price and service rate at each time step $A_t=(p,\mu)$.
	%\item %\textit{Transition probabilities:} Transitions of the queue length are determined by the arrival rate $\lambda(p)$ and service rate $\mu$.
	\item \textit{Rewards:} The net profit gained in time slot $t$, denoted by $R_t$.
\end{itemize} 
Under the above setting, we write the Bellman equation as below:
\begin{align*}
	q_\pi(s,a)+h_{\pi}=\EE_{s'\sim P(s,a),a'\sim\pi}[R(s,a)+q(s',a')], 
\end{align*}
with $q_{\pi}(s,a)$ is the  $q$-function of policy $\pi$ and $h_\pi$ is the long-run average revenue under $\pi$.

Following \citep[Section 13.6]{SB18}, we apply the Gaussian parameterization for our actions. Specifically, we draw 
$p\sim N(\bar{p},\sigma_p^2)$ and $\mu\sim N(\bar{\mu},\sigma_\mu^2)$ independently. Let
$\theta\equiv(\bar{p},\bar{\mu},\sigma_p^2,\sigma_\mu^2)$ and we denote $\pi_\theta$ as the Gaussian density with parameter $\theta$. %Here we implicitly consider static policies independent of the states for two purposes: (i) our objective intends to search for the optimal static policy and we make this point consistent;(ii) refraining the policy class into the static policy class would ease the burden of learning.
%Based on the normal representation, we next describe the policy gradient method and the full PG algorithm in our problem is given in Algorithm \ref{alg: PG}. To better describe the policy gradient method, we further introduce the following notations:
%\begin{itemize}
%	\item Long-run Average revenue: $J(\theta)=r(\pi)=\lim\limits_{T\rightarrow\infty} \frac{1}{T}\sum_{t=1}^{T}\EE[R(t)|S_t,A_t\sim \pi_\theta]$
%	\item Differential Total Rewards: $G_t=R_t-r(\pi)+R_{t+1}-r(\pi)+\cdots$, with $R_t$ being the reward within time slot $t$.
%	%\item Differential Value Functions: $v^\pi(s)=\EE[G_t| S_t=s]$
%	\item Differential q-function: $q^\pi(s,a)=\EE[G_t|S_t=s,A_t=a]$
%\end{itemize} 
According to the policy gradient theorem \cite[p.339]{SB18}, the gradient on the policy function can be represented as 
$\nabla_\theta h_{\pi_\theta}=\EE[\nabla \log \pi_\theta(A_t|S_t)\cdot q(S_t,A_t)]$.

We next quickly explain how the PG algorithm works. We organize the time into successive cycles each of which contains several episodes. In each episode, the PG algorithm operates the system under policy $\pi_\theta$ and generates a sample of gradient estimator. The averaged value of these samples from different episode gives the PG estimator in each cycle. See  Algorithm \ref{alg: PG} in the appendix for the detailed description of the PG algorithm.

%\begin{algorithm}[h]
%	\SetAlgoLined
%	\KwIn{normal parameterization $\pi(a|\theta)$, step size $\eta>0$, initial policy parameter $\theta:(\bar{p}_1,\bar{\mu}_1,\sigma_{p,1}^2,\sigma_{\mu,1}^2)$, cycle length $L$ (how many episodes in each episode), episode length $T$ (how many time slots in each episode)\;}
%	\For{each cycle}{
	%		\For{episode $i=1:L$}
	%		{
		%			Generate an episode $Q_1,(p_1,\mu_1),R_1,\cdots,Q_{T-1},(p_{T},\mu_{T}),R_T$ following $\pi^\theta$\;
		%			$\bar{R}=\frac{1}{T}\sum_{t=1}^{T}R_T$\;
		%			\For{$t=1,\cdots,T$}{
			%				$G=\sum_{k=t}^{T}R_k-\bar{R}$\;
			%				$\hat{\nabla}_{i,t}\leftarrow G\cdot  \begin{pmatrix}
				%					(p_t-\bar{p})/\sigma_p^2\\
				%					(\mu_t-\bar{\mu})/\sigma_\mu^2\\
				%					\left[(p_t-\bar{p})^2-\sigma_p^2\right]/\sigma_p^3\\
				%					\left[(\mu_t-\bar{p})^2-\sigma_\mu^2\right]/\sigma_\mu^3
				%				\end{pmatrix}$\;
			%			}
		%			$\hat{\nabla}_i=\frac{1}{T}\sum_{t=1}^{T}\hat{\nabla}_{i,t}$\;
		%		}
	%		$\theta\leftarrow \theta +\eta \cdot \frac{1}{L}\sum_{i=1}^{L}\hat{\nabla}_i$\;
	%	}
%	\caption{Episodic Policy Gradient method}
%	\label{alg: PG}
%\end{algorithm}

\paragraph{Experiment settings and results.} We now compare PG with LiQUAR using our base example as described in Section 6.1. Specifically, we consider an $M/M/1$ queue having logit demand function \eqref{eq: logit demand} with $M_0=10$, $a=4.1$, $b=1$ and exponential service times with holding cost $c(\mu)=\mu$. 
For LiQUAR, the hyperparameter are $\eta_k=4k^{-1}$, $T_k=200k^{1/3}$, and $\delta_k=\min(0.1,0.5k^{-1/3})$. For PG, we consider pick the step size $\eta\in\{0.1,0.01,0.001,0.0001\}$ and the cycle length $T\in\{10,100,300\}$. In addition, we keep the total running time of LiQUAR and PG equal in order for a fair comparison.

\begin{figure}
	\centering
	\includegraphics[width=0.8\linewidth]{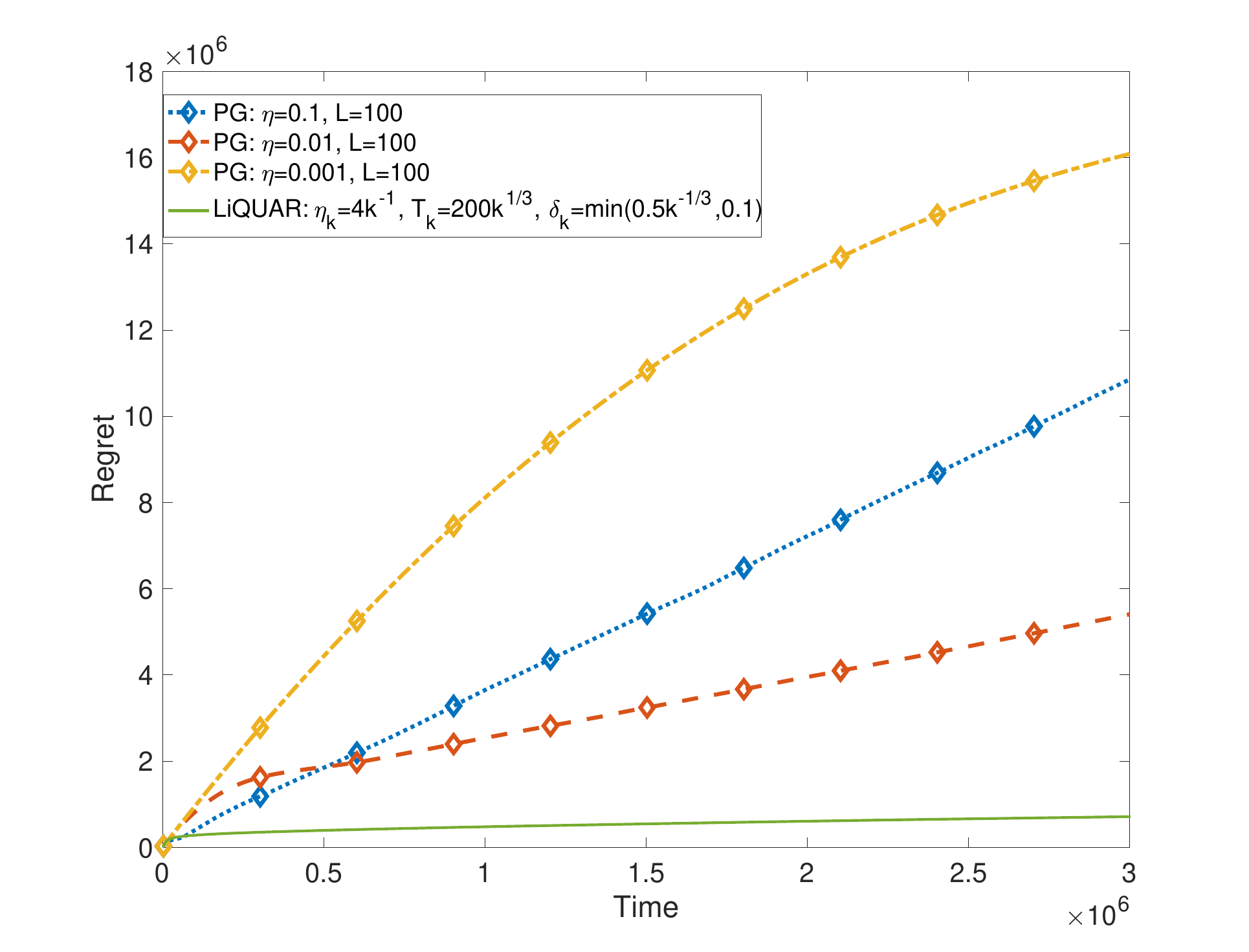}
	\caption{Comparison of the regret of LiQUAR with PG in our base example. The hyperparameter choices are: (i) LiQUAR : $\eta_k=4k^{-1},T_k=200k^{1/3},\delta_k=\min (0.5k^{1/3},0.1)$; (ii) PG: $\eta\in\{0.1,0.01,0.001\}, L\in\{1,10,100,300\}$, Episode length $T=3000/L$. All regret curves are estimated from 1,000 independent simulation replications.}
	\label{fig: PG}
\end{figure}

We report the regret curves in Figure \ref{fig: PG}. For the clarity of the figure, we report the curves with the lowest regret for each step length choice $\eta$. From Figure \ref{fig: PG}, we find that  LiQUAR is more effective than PG in a wide range of hyper-parameter choices.
\begin{remark}[LiQUAR vs. PG]
	In the PG algorithm, the gradient estimator relies on accurately learning the $q_\pi(s,a)$ function (as outlined in the policy gradient theorem), which is a two-dimensional function. Inaccuracies in this estimation can result in significant variance in the gradient calculations. In contrast, LiQUAR only requires learning the values of individual actions, substantially reducing the complexity and effort required for learning. Furthermore, by the design of LiQUAR, the tuning of its hyperparameters can leverage domain-specific knowledge of the queueing system, e.g., the transient bias and auto-correlation between queueing data; also see Remark \ref{remark: queue structure utilization}. In comparison, tuning hyperparameters in the PG algorithm is considerably more challenging, as RL methods are generally considered black-box approaches.
\end{remark}

\section{Conclusions}\label{sec: Conclusions}
In this paper we develop an online learning framework, dubbed LiQUAR, designed for dynamic pricing and staffing in an $M/GI/1$ queue with unknown arrival rate function and service distribution. LiQUAR's main appeal is its ``model-free" attribute. 
Unlike 
%Instead of using the precise knowledge of the demand function and the service-time distribution as required in 
the conventional ``predict-then-optimize" approach where precise estimations of the demand function and service distribution must be conducted (as a separate step) before the decisions may be optimized, LiQUAR is an integrated method that recursively evolves the control policy to optimality by effectively using the newly generated queueing data (e.g., arrival times and service times). 
LiQUAR's main advantage is its solution robustness; its algorithm design is able to automatically relate the parameter estimation errors to the fidelity of the optimized solutions. Comparing to the conventional method, this advantage becomes more significant when the system is in  heavy traffic.

Effectiveness
of LiQUAR is substantiated by (i) theoretical results including the
algorithm convergence and regret analysis, and (ii) engineering confirmation via simulation
experiments of a variety of representative  queueing models. Theoretical analysis of the
regret bound in the present paper may shed lights on the design of efficient online learning algorithms (e.g., bounding gradient estimation error and controlling proper learning rate) for more general queueing systems. In addition, the analysis on the statistical properties for our gradient estimator has independent interests and may contribute to the general literature of stochastic gradient decent.
We also extend LiQUAR to the more general $GI/GI/1$ model and confirm its good performance by conducting numerical studies.

There are several venues for future research. One dimension is to extend the method to queueing models under more general settings such as non-Poisson arrivals, customer abandonment and multiple servers, which will make the framework more practical for service systems such as call centers and healthcare. Another interesting direction is
to theoretically relax the assumption of uniform stability by developing a ``smarter" algorithm that automatically explore and then stick to control policies that guarantee a stable system performance. %otherwise the decision maker will face a low profit (high cost).

\bibliography{references}
\bibliographystyle{chicago}
%%%%%%%%%%%%%%%%%%%%%%%%%%%%%%%%%%%%%%%%%%%%%%%%%%%%%%%%%%%%%%%%%%%%%%%%%%%%%%%%%%%%%%%%%%%%%%%%%%%%%%%%%%%%%%%%%%
\newpage
\appendix
\begin{center}
{\large\bf SUPPLEMENTARY MATERIAL}
\end{center}
This e-companion provides supplementary materials to the main paper. 
In Section \ref{sec: proofs}, we provide the proofs for our main results in the main paper. 
In Section \ref{appx: proofs}, we supplement Section \ref{sec: proofs} to give additional proofs. 
In Section \ref{sec: lower bound}, we give a regret lower bound for the $M/M/1$ queue.
In Section \ref{appx: demand function}, we verify that the Condition (a) of Assumption \ref{assmpt: uniform} holds for some commonly used demand functions. %including: (i) linear demand functions; (ii) quadratic demand functions; (iii) exponential demand functions; (iv) logit demand functions. 
In Section \ref{appx: additional numeric}, we conduct additional numerical studies. 
In Section \ref{sec: AlgRL}, we provide the detailed description for the PG algorithm in Section \ref{sec: RL}. To facilitate readability, all notations are summarized in Table \ref{tab: notations} including all model parameters and functions, algorithmic parameters and variables, and constants in the regret analysis.
%\begin{APPENDICES}

\section{Proofs of Main Results}\label{sec: proofs}
In this section, we provide the proofs of the main theorems and propositions. Proofs of technical lemmas are given in the Section \ref{appx: proofs}.

\subsection{Proof of Proposition \ref{prop: delay observation}}\label{subsec: proof of delay error}
First, we introduce a technical lemma to uniformly bound the moments of workload under arbitrary control policies.
\begin{lemma}[Uniform Moment Bounds]\label{lmm: uniform bound short} Under Assumptions \ref{assmpt: uniform} and \ref{assmpt: light tail}, there exist some constants $\theta_0>0$ and  $M>1$ such that, for any sequence of control parameters $\{(\mu_l,p_l):l\geq 1\}$, 
	$$\EE[W_l(t)^m]\leq M, \quad\EE[W_l(t)^m\exp(2\theta_0 W_l(t))]\leq M,$$
	for all $m\in\{0, 1,2\}$, $l\geq 1$ and $0\leq t\leq T_{k}$ with $k=\lceil l/2 \rceil$.
\end{lemma}
Then, following \eqref{eq: hat W},
%$$\EE_l\left[|\hat{W}_l(t)-W_l(t)|\right]=\EE_l\left[W_l(t)\cdot\ind{t>A_{N^D(l)}^{C^D(l)}}\right].$$
%By definition, the event $\{t\leq A_{N^D(l)}^{C^D(l)}\}$ happens as long as the workload $W_l(t)$ are all served by $T_k$, i.e. $W_l(t)\leq \mu_l(T_k-t)$. Consequently,
\begin{align*}	%\EE_l\left[|\hat{W}_l(t)-W_l(t)|\right]&=\EE_l\left[W_l(t)\cdot\ind{t>A_{N^D(l)}^{C^D(l)}}\right]=
	\EE&\left[|\hat{W}_l(t)-W_l(t)|\right]=\EE\left[W_l(t)\cdot\ind{W_l(t)>\mu_l(T_k-t)}\right]\leq \EE\left[W_l(t)^2\right]^{1/2}\PP\left(W_l(t)>\mu_l(T_k-t)\right)^{1/2}\\
	&\leq \EE\left[W_l(t)^2\right]^{1/2}\cdot\exp\left(-\frac{1}{2}\theta_0\mu_l(T_k-t)\right)\EE\left[\exp(\theta_0W_l(t))\right]^{1/2}
	\leq \exp\left(-\frac{1}{2}\theta_0\underline{\mu}(T_k-t)\right)M,
\end{align*}
where  the last inequality follows from  Lemma \ref{lmm: uniform bound short}.
\hfill $\Box$
\subsection{Proof of Proposition \ref{prop: estimation error}}\label{subsec: proof of estimation error}
For each cycle $l$, the difference between the estimated system performance $\hat{f}^G(\mu_l,p_l)$ and its true value is
\begin{align*}
	\hat{f}^G(\mu_l,p_l)-f(\mu_l,p_l)=\frac{-p_l(N_l-\lambda(p_l)T_k)}{T_k}+\frac{1}{(1-2\alpha)T_k}\int_{\alpha T_k}^{(1-\alpha)T_k}[\underbrace{\hat{W}_l(t)-W_l(t)}_{\text{delayed observation}} + \underbrace{W_l(t)-w_l}_{\text{transient error}}]~dt,
\end{align*}
where $w_l=\EE[W_\infty(\mu_l,p_l)]$ is the steady-state mean workload.
To bound the moments of this difference, which correspond to the bias and MSE of $\hat{f}^G(\mu_l,p_l)$, we construct a stationary workload process $\bar{W}_l(t)$ for $0\leq t\leq T_{k}$. At $t=0$, the initial value $\bar{W}^l(0)$ is independently drawn from the stationary distribution $W_\infty(\mu_l,p_l)$ and $\bar{W}_l(t)$ is \textit{synchronously coupled} with $W_l(t)$ in the sense that they share the same sequence of arrivals and individual workload on $[0,T_k]$. 

\noindent\textbf{Bound on the Bias.} The \textit{bias} of $\hat{f}^G(\mu_l,p_l)$ can be decomposed as
\begin{align}\label{eq: bias decomposition}
	&\EE_l\left[\hat{f}^G(\mu_l,p_l)-f(\mu_l,p_l)\right]\notag\\
	=~&\frac{1}{(1-2\alpha)T_{k}}\int_{\alpha T_{k}}^{(1-\alpha)T_{k}}\left(\EE_l\left[\hat{W}_l(t) \right]-\EE_l[\bar{W}_l(t)]\right)dt
	\leq \frac{1}{(1-2\alpha)T_k}\int_{\alpha T_k}^{(1-\alpha)T_k}\EE_l\left[|\hat{W}_l(t)-\bar{W}_l(t)|\right]dt.\notag\\
	\leq~& \frac{1}{(1-2\alpha)T_k}\left(\int_{\alpha T_k}^{(1-\alpha)T_k}\EE_l[|\hat{W}_l(t)-W_l(t)|]dt+\int_{\alpha T_k}^{(1-\alpha)T_k}\EE_l[|W_l(t)-\bar{W}_l(t)|]dt\right).
\end{align}
The first term in \eqref{eq: bias decomposition} is the error caused by delayed observation. Following the same analysis as in Section \ref{subsec: proof of delay error}, 
\begin{align*}
	\EE_l\left[|\hat{W}_l(t) - W_l(t)|\right] %&=\EE_l\left[ W_l(t)\cdot1\left(t>A_{n(l)}^{c(l)}\right)\right]\leq \EE_l\left[W_l(t)\cdot1\left(W_l(t)> \mu_l(T_k-t)\right)\right]\\
	%&\leq \EE_l[W_l(t) ^2]^{1/2}\PP_l\left(W_l(t)> \mu_l(T_k-t)\right)^{1/2}\\
	\leq \EE_l[W_l(t) ^2]^{1/2}\cdot \exp(-a\mu_l(T_k-t))\EE_l[\exp(2aW_l(t))]^{1/2},
\end{align*}
for $a=\theta_0/2$. It is easy to check that $W_l(t) \leq  W_l(0)+\bar{W}_l(t)$. Conditional on $\mathcal{G}_l$, for all $0\leq t\leq T_k$, $\bar{W}_l(t)$ is the stationary workload with parameter $(\mu_l,p_l)$. Following the proof of Lemma \ref{lmm: uniform bound short},  $\bar{W}_l(t)$ is stochastic bounded by the stationary workload with parameter $(\underline{\mu},\underline{p})$. Therefore,
\begin{align}\label{eq: bias bound 1}
	&\EE_l\left[|\hat{W}_l(t) - W_l(t)|\right] \leq \EE_l[W_l(t) ^2]^{1/2}\cdot \exp(-\theta_0\mu_l(T_k-t)/2)\EE_l[\exp(\theta_0 W_l(t))]^{1/2}\notag\\
	\leq ~&\exp(-\theta_0\mu_l(T_k-t)/2)(W_l(0)^2 + 2W_l(0)\EE_l[\bar{W}_l(t)]+\EE_l[\bar{W}_l(t)^2])^{1/2} \exp(\theta_0 W_l(0))\EE_l[\exp(\theta_0\bar{W}_l(t))]^{1/2}\notag\\
	\leq~& \exp(-\theta_0\mu_l(T_k-t)/2)(W_l(0)^2 + 2MW_l(0)+M^2)^{1/2} \exp(\theta_0 W_l(0))M^{1/2}\notag\\
	\leq~& \exp(-\theta_0 \mu_l(T_k-t)/2)M(M+W_l(0))\exp(\theta_0 W_l(0)).
\end{align}
The last inequality holds as  $M\geq1$. {The second term in \eqref{eq: bias decomposition} will be bounded using the following lemma on convergence rate of  two synchronously coupled workload processes. }

\begin{lemma}[Ergodicity Convergence]\label{lmm: ergodicity}
	Suppose Assumptions \ref{assmpt: uniform} and \ref{assmpt: light tail} hold. Two workload processes $W(t)$ and $\bar{W}(t)$ with equal control parameters $(\mu,p)\in\mathcal{B}$ are synchronously coupled with initial states $(W(0), \bar{W}(0))$. Then, there exists $\gamma>0$ independent of $(\mu,p)$, such that
	$$\EE\left[|W(t)-\bar{W}(t)|^m~|~W(0), \bar{W}(0) \right]\leq e^{-\gamma t}(e^{\theta_0 W(0)}+e^{\theta_0 \bar{W}(0)})|W(0)-\bar{W}(0)|^m.$$
\end{lemma}
%for \red{$a= \eta$} according to Lemma \ref{lmm: uniform bound}.\\
{Using this lemma, we can compute}
\begin{align}\label{eq: bias bound 2}
	\EE_l[|W_l(t)-\bar{W}_l(t)|]&\leq \exp(-\gamma t)\EE_l\left[|W_l(0)-\bar{W}_l(0)|(\exp(\theta_0 W_l(0))+\exp(\theta_0 \bar{W}_l(0)))\right]\notag\\
	&\leq \exp(-\gamma t) \left(W_l(0)\exp(\theta_0 W_l(0)) + MW_l(0)+M\exp(\theta_0 W_l(0)) + M\right)\notag\\
	&\leq \exp(-\gamma t)(M+W_l(0))(\exp(\theta_0 W_l(0))+M).
\end{align}
Let $\theta_1 = \min(\gamma, \theta_0\underline{\mu}/2)$.  Plugging inequalities \eqref{eq: bias bound 1} and \eqref{eq: bias bound 2} into \eqref{eq: bias decomposition}, we obtain the following bound for the  bias
\begin{align*}
	\left|\EE_l\left[\hat{f}^G(\mu_l,p_l)-f(\mu_l,p_l)\right]\right|\leq \frac{1}{(1-2\alpha)T_k}\cdot\frac{2\exp(-\theta_1\alpha T_k)}{\theta_1}\cdot M(M+W_l(0))(\exp(\theta_0 W_l(0))+M).
\end{align*}

\noindent\textbf{Bound on the Mean Square Error.} The mean square error (MSE) of $\hat{f}^G(\mu_l,p_l)$
$$\EE_l[(\hat{f}^G(\mu_l,p_l)-f(\mu_l,p_l))^2]\leq 2\EE_l[E_1^2]+2\EE_l[E_2^2],$$
with
$$\hat{f}^G(\mu_l,p_l) -f(\mu_l,p_l)= \underbrace{\frac{-p_l(N_l-\lambda(p_l)T_{k})}{T_{k}}}_{E_1}+\underbrace{\frac{1}{(1-2\alpha)T_{k}}\int_{\alpha T_{k}}^{(1-\alpha)T_{k}}(\hat{W}_l(t)-w_l)dt}_{E_2}.$$
%Let  $\mathcal{G}_l$ be $\sigma$-algebra generated by $\{W_m(t): 1\leq m\leq l-1, 0\leq t\leq T_{k(m)}\}$ and $Z_k$ for $l=2k-1, 2k$. Then $(\mu_l,p_l), W_l(0)\in\mathcal{G}_l$. We write $\EE_l[\cdot] \equiv \EE[\cdot|\mathcal{G}_l]$. 
Conditional on $\mathcal{G}_l$, the observed number of arrivals $N_l$ is a Poisson r.v. with mean $\lambda(p_l)T_k$. So, $
\EE_l[E_1^2] = p_l^2\lambda(p_l)T_k^{-1}\leq \bar{p}^2\bar{\lambda}T_k^{-1}$.

For $E_2$,  we have
\begin{align*}
	\EE_l[E_2^2] = \frac{1}{(1-2\alpha)^2T_k^2}\int_{\alpha T_k}^{(1-\alpha)T_k}\int_{\alpha T_k}^{(1-\alpha)T_k}\EE_l\left[(\hat{W}_l(t)-w_l)(\hat{W}_l(s)-w_l)\right] dt ds.
\end{align*}
According to \eqref{eq: hat W}, $\hat{W}_l(\cdot)\leq W_l(\cdot)$ and therefore, for any $0\leq s\leq  t\leq T_k$, 
\begin{align*}
	&\EE_l[(\hat{W}_l(t)-w_l)(\hat{W}_l(s)-w_l)] =\EE_l[\hat{W}_l(t)\hat{W}_l(s)-w_l(\hat{W}_l(s)+\hat{W}_l(t))+w^2_l]\\
	\leq~& \EE_l[W_l(t){W}_l(s)-w_l(\hat{W}_l(s)+\hat{W}_l(t))+w^2_l]\\
	\leq ~&\EE_l[({W}_l(t)-w_l)({W}_l(s)-w_l)] + \left(\EE_l[w_l|W_l(s)-\hat{W}_l(s)|]+\EE_l[w_l|W_l(t)-\hat{W}_l(t)|]\right)\\
	\leq ~&\underbrace{\EE_l[({W}_l(t)-w_l)({W}_l(s)-w_l)]}_{\text{auto-covariance}} + M\underbrace{\left(\EE_l[|W_l(s)-\hat{W}_l(s)|]+\EE_l[|W_l(t)-\hat{W}_l(t)|]\right)}_{\text{error caused by delayed observations}}
\end{align*}
{To bound the auto-covariance term, we introduce the following lemma.}
\begin{lemma}[Auto-covariance of $W_l(t)$]\label{lmm: autocovariance W}
	There exists a constant $K_V> 0$ independent of $T_k,l,p_l,\mu_l$ such that, for any $l\geq 1$ and $0\leq  s\leq t\leq T_k$,
	\begin{equation}\label{eq: covariance bound}
		\EE_l[({W}_l(t)-w_l)({W}_l(s)-w_l)] \leq K_V\left(\exp(-\gamma(t-s))+\exp(-\gamma s)\right)(W_l(0)^2+1)\exp(\theta_0W_l(0)).
	\end{equation}
\end{lemma}
%\begin{remark}
%	The upper bound contains two terms of $\exp(-\gamma(t-s))$ and $\exp(-\gamma s)$. Intuitively, the first term comes from the fact that the dependence of workload process $W_l(t)$ on the initial state $W_l(s)$ decays exponentially fast, as proved in Lemma \ref{lmm: ergodicity}, and the second term comes from the transient bias $\EE[W_l(s)-w_l]$.
%\end{remark}

\noindent {Following \eqref{eq: covariance bound},} we write
\begin{align*}
	&\frac{1}{(1-2\alpha)^2T_k^2}\int_{\alpha T_k}^{(1-\alpha) T_k}\int_{\alpha T_k}^{(1-\alpha)T_k}\EE_l[(W_l(t)-w_l)(W_l(s)-w_l)]dtds\\
	\leq ~&\frac{2K_V(W_l(0)^2+1)\exp(\theta_0W_l(0))}{(1-2\alpha)^2T_k^2}\int_{\alpha T_k}^{(1-\alpha) T_k}\int_{\alpha T_k}^{t} (\exp(-\gamma(t-s))+\exp(-\gamma s))dsdt\\
	=~&\frac{2K_V(W_l(0)^2+1)\exp(\theta_0W_l(0))}{(1-2\alpha)^2T_k^2}\int_{\alpha T_k}^{(1-\alpha) T_k}\gamma^{-1}(1-\exp(-\gamma(t-\alpha T_k))+\exp(-\gamma \alpha T_k)-\exp(-\gamma t))dt\\
	\leq~&\frac{2K_V(W_l(0)^2+1)\exp(\theta_0W_l(0))}{(1-2\alpha)^2T_k^2}\int_{\alpha T_k}^{(1-\alpha) T_k}2\gamma^{-1}dt\leq\frac{4K_V(W_l(0)^2+1)\exp(\theta_0W_l(0))}{\gamma(1-2\alpha)T_k}.%=O(T_k^{-1}).
\end{align*}
For the error of delayed observation, by Proposition \ref{prop: delay observation},  we have
\begin{align*}
	&\frac{1}{(1-2\alpha)^2T_k^2}\int_{\alpha T_k}^{(1-\alpha) T_k}\int_{\alpha T_k}^{(1-\alpha) T_k}\left(\EE_l[|W_l(s)-\hat{W}_l(s)|]+\EE_l[|W_l(t)-\hat{W}_l(t)|]\right)dsdt\\
	\leq~&\frac{M(M+W_l(0))\exp(\theta_0 W_l(0))}{(1-2\alpha)^2T_k^2}\int_{\alpha T_k}^{(1-\alpha) T_k}\int_{\alpha T_k}^{(1-\alpha) T_k}\left(\exp(-\frac{\theta_0\mu_l}{2}(T_k-t))+\exp(-\frac{\theta_0\mu_l}{2}(T_k-s))\right)dsdt\\
	=~&\frac{4M(M+W_l(0))\exp(\theta_0 W_l(0))}{\theta_0\mu_l(1-2\alpha)T_k}\left(\exp(-\frac{\theta_0\mu_l}{2}\alpha T_k)-\exp(-\frac{\theta_0\mu_l}{2}(1-\alpha)T_k)\right)\\
	\leq~& \frac{4M(M+W_l(0))\exp(\theta_0 W_l(0))}{\theta_0\mu_l(1-2\alpha)T_k}.%=O(T_k^{-1}).
\end{align*}
%Recall that $\theta_2 = \min(\gamma, \theta_0\underline{\mu}/2)$. As $W_l(0)\leq (W_l(0)^2+1)/2$ 
As $W_l(0)\leq (W_l(0)^2+1)/2$ and $M\geq 1$, we have $M+W_l(0)\leq (M+1)(1+W_l(0)^2)$. %\red{In addition,  $\theta_2=\min(\gamma,\theta_0\underline{\mu}/2)$. ???} 、
Then, if we choose 
\begin{equation}
	\label{eq: KM}
	K_M=\frac{8(K_V+M^3+M^2)}{(1-2\alpha)\min(\gamma,\theta_0\underline{\mu})}+2\bar{p}^2\bar{\lambda}
\end{equation}
then, we have
$$\EE_l[(\hat{f}^G(\mu_l,p_l)-f(\mu_l,p_l))^2]\leq 2\EE_l[E_1^2]+2\EE_l[E_2^2]\leq K_MT_k^{-1}(W_l(0)^2+1)\exp(\theta_0W_l(0)).$$\hfill $\Box$
\subsection{Proof of Proposition \ref{prop: bound Bk V_k}}\label{subsec: proof of boound Bk Vk}
%The gradient estimation error
%		\begin{align}\label{eq: H_k decomposition}
	%			H_k& = \delta_k^{-1}\sum_{l=2k-1}^{2k}\left(\hat{f}^G(\mu_l,p_l)-f(\mu_l,p_l)\right) +\frac{f(\mu_{2k},p_{2k})-f(\mu_{2k-1},p_{2k-1})}{\delta_k}.
	%		\end{align}
According to the following lemma, the FD approximation error is of order $O(\delta_k^2)$.
\begin{lemma}\label{lmm: fd approximation}
	Under Assumption \ref{assmpt: uniform}, there exists a smoothness constant $c>0$ such that for any $\mu_1,\mu_2,\mu\in[\underline{\mu},\bar{\mu}]$ and $p_1,p_2,p\in[\underline{p},\bar{p}]$,
	\begin{align*}
		\left\vert\frac{f(\mu_1,p)-f(\mu_2,p)}{\mu_1-\mu_2}-\partial_\mu f\left(\frac{\mu_1+\mu_2}{2},p\right)\right\vert\leq c(\mu_1-\mu_2)^2\\
		\left\vert\frac{f(\mu,p_1)-f(\mu,p_2)}{p_1-p_2}-\partial_p f\left(\mu, \frac{p_1+p_2}{2}\right)\right\vert\leq c(p_1-p_2)^2.
	\end{align*}
\end{lemma}
%First, we show by Lemma \ref{lmm: fd approximation} in Section \ref{appx: proofs of convergence} that there exists a universal constant $c>0$ such that the finite-difference approximation error
%\begin{align}\label{eq: FD approximation error}
%		\left|\frac{f(\mu_{2k},p_{2k})-f(\mu_{2k-1},p_{2k-1})}{\delta_k} - \partial_\mu f(\bar{\mu}_k,\bar{p}_k)\right|%\leq c\delta_k^2,
%\leq c\delta_k^2.
%		\end{align}
So, to bound $B_k$, it remains to show that $$\EE[\EE[\hat{f}^G(\mu_l,p_l)-f(\mu_l,p_l)|\mathcal{F}_k]^2]^{1/2}=O(\exp(-\theta_1\alpha T_k)).$$
Recall that $\mathcal{F}_k$ is the $\sigma$-algebra including all events in the first $2(k-2)$ cycles, so $\mathcal{F}_k\subseteq\mathcal{G}_l$ for $l=2k-1,2k$. By Jensen's inequality,
\begin{align*}
\EE\left[\hat{f}^G(\mu_l,p_l)-f(\mu_l,p_l)|\mathcal{F}_k\right]^2&=\EE\left[\EE_l\left[\hat{f}^G(\mu_l,p_l)-f(\mu_l,p_l)\right]\Big|\mathcal{F}_k\right]^2\\
&\leq \EE\left[\EE_l\left[\hat{f}^G(\mu_l,p_l)-f(\mu_l,p_l)\right]^2\Big|\mathcal{F}_k\right].
\end{align*}
Therefore, $$\EE\left[\EE\left[\hat{f}^G(\mu_l,p_l)-f(\mu_l,p_l)|\mathcal{F}_k\right]^2\right]\leq\EE\left[\EE_l\left[\hat{f}^G(\mu_l,p_l)-f(\mu_l,p_l)\right]^2\right].$$
%To analyze the approximation error $\hat{f}^G-f$, we define another sequence of filtration $\mathcal{G}_l$ as the $\sigma$-algebra generated by $\{W_m(t): 1\leq m\leq l-1, 0\leq t\leq T_{k(m)}\}$ and $Z_k$ for $l=2k-1, 2k$. By definition, $\mathcal{F}_k\subseteq \mathcal{G}_l$ for $l=2k-1, 2k$, and as a consequence,
%$$\EE[\EE[\hat{f}^G(\mu_l,p_l)-f(\mu_l,p_l)|\mathcal{G}_l]^2]\geq \EE[\EE[\hat{f}^G(\mu_l,p_l)-f(\mu_l,p_l)|\mathcal{F}_k]^2].$$
%So it suffices to show that 
%$$\EE[\EE[\hat{f}^G(\mu_l,p_l)-f(\mu_l,p_l)|\mathcal{G}_l]^2]=O(\exp(-2\theta_2\alpha T_k)).$$
%To simplify the notation, we write $\EE[\cdot|\mathcal{G}_l] = \EE_l[\cdot]$ and $\PP(\cdot|\mathcal{G}_l) = \PP_l(\cdot)$. 
By Proposition \ref{prop: estimation error}, the bias of estimated system performance
$$\left|\EE_l\left[\hat{f}^G(\mu_{l},p_{l})-f(\mu_{l},p_{l}) \right]\right|\leq \frac{2\exp(-\theta_1\alpha T_k)}{(1-2\alpha)\theta_1 T_k}\cdot M(M+W_l(0))(\exp(\theta_0 W_l(0))+M).$$
As  $(x+y)^2\leq 2x^2+2y^2$, we have, by Lemma \ref{lmm: uniform bound short}, 
\begin{align*}
&\EE[\EE_l[\hat{f}^G(\mu_l,p_l)-f(\mu_l,p_l)]^2]\\
\leq& \frac{4\exp(-2\theta_1\alpha T_k)}{(1-2\alpha)^2\theta_1^2T_k^2}\left(4M^4\EE[\exp(2\theta_0 W_l(0))+W_l(0)^2]+4M^2\EE[W_l(0)^2\exp(2\theta_0 W_l(0))]+4M^6\right)\\
\leq& \frac{4\exp(-2\theta_1\alpha T_k)}{(1-2\alpha)^2\theta_1^2T_k^2}\cdot(8M^5+4M^3+4M^6)=O(\exp(-2\theta_1\alpha T_k)).
\end{align*}
Therefore, $B_k=O\left(\delta_k^2+\delta_k^{-1}\exp(-\theta_1 \alpha T_k)\right)$.
The variance
$$\EE[\|H_k\|^2]\leq 3\delta_k^{-2}\sum_{l=2k-1}^{2k}\EE[(\hat{f}^G(\mu_l,p_l)-f(\mu_l,p_l))^2] + 3\delta_k^{-2}\EE[(f(\mu_{2k},p_{2k})-f(\mu_{2k-1},p_{2k-1}))^2].$$
By the smoothness condition of the objective function $f(x)$ as given in Assumption \ref{assmpt: uniform}, 
$$3\delta_k^{-2}\EE[(f(\mu_{2k},p_{2k})-f(\mu_{2k-1},p_{2k-1}))^2]\leq \max_{(\mu,p)\in\mathcal{B}}\|\nabla f(\mu,p)\|^2=O(1).$$ 
Following Proposition \ref{prop: estimation error}, for $l=2k-1, 2k$,
\begin{align*}
\EE[(\hat{f}^G(\mu_l,p_l)-f(\mu_l,p_l))^2] \leq K_MT_k^{-1}\EE[(W_l(0)^2+1)\exp(\theta_0W_l(0))]=O(T_k^{-1}).
\end{align*}
Therefore, $\EE[\|H_k\|^2]=O(\delta_k^{-2}T_k^{-1}\vee1)$.
\hfill $\Box{}$
\subsection{Proof of Theorem \ref{thm: convergence}}\label{subsec: proof of convergence}
To obtain convergence of the SGD iteration, we first need to establish a desirable convex structure of the objective function \eqref{eqObjmin}.
\begin{lemma}[Convexity and Smoothness of $f(\mu,p)$]\label{lmm: convexity}
Suppose Assumption \ref{assmpt: uniform} holds. Then, there exist finite positive constants $0<K_0\leq 1$ and $K_1>K_0$ such that for all $x=(\mu,p)\in\mathcal{B}$,
\begin{enumerate}
	%\item[$(a)$] $\partial^2_\mu f(\mu,p), \partial^2_p f(\mu,p)\geq K_0$;
	\item[$(a)$] $(\bx-\bx^*)^T\nabla f(x)\geq K_0\|\bx-\bx^*\|^2$,
	%\item[$(b)$] $\|\nabla f(x)\|\leq K_1\|x-x^*\|$;
	\item[$(b)$] $|\partial^3_\mu f(\bx)|,|\partial^3_p f(\bx)|\leq K_1$.
	%\textcolor{red}{\item[$(c)$] $f(x)$ has bounded third derivatives, i.e.  $|\frac{\partial^3}{\partial \mu^3}f(\mu,p)|, |\frac{\partial^3}{\partial p^3}f(\mu,p)|\leq K_1$.}
\end{enumerate}
\end{lemma} 
We only sketch the key ideas in the proof of the convergence result \eqref{eq: convergence beta} under the convexity structure here; the full proof is given in Appendix \ref{appx: proof of convergence theorem}. Let $b_k = \EE[\|\bar{\bx}_k-\bx^*\|^2] $. Then, following the SGD recursion and some algebra, we get the following recursion on $b_k$:
$$b_{k+1}\leq (1-2K_0\eta_k+\eta_k B_k)b_k+\eta_kB_k+\eta_k^2\mathcal{V}_k.$$
Under condition \eqref{eq: beta conditions}, we can show that the recursion coefficient  $1-2K_0\eta_k+\eta_k B_k<1$, so $b_k$ eventually converges to 0. With more careful calculation as given in Appendix \ref{appx: proof of convergence theorem}, we can obtain the convergence rate \eqref{eq: convergence beta} by induction using the above recursion.

Applying the convergence result \eqref{eq: convergence beta} to LiQUAR relies on knowing the bounds on $B_k$ and $\mathcal{V}_k$. Given Proposition \ref{prop: bound Bk V_k}, one can check that, if $\eta_k =O(k^{-a})$, $T_k=O(k^b)$ and $\delta_k=O(k^{-c})$, the bounds for $B_k$ and $\mathcal{V}_k$ as specified in condition \eqref{eq: beta conditions} holds with  $\beta = \max(-a, -a-b+2c, -2c)$. Then, \eqref{eq: convergence rate} follows immediately from \eqref{eq: convergence beta}.
\subsection{Proof of Proposition \ref{prop: regret of nonstationarity}}\label{subsec: proof of nonstationary regret}
The regret of nonstationarity
\begin{align*}
R_{2k}&=\sum_{l=2k-1}^{2k}\EE[\rho_l - T_k f(x_l)]
= \sum_{l=2k-1}^{2k}\EE\left[h_0\int_0^{T_k} (W_l(t)-w_l)dt -  p_l(N_l-T_k\lambda(p_l))\right],
\end{align*}
where $w_l = \EE_l[W_\infty(\mu_l,p_l)]$. Conditional on $p_l$, $N_l$ is a Poisson random variable with mean $T_k\lambda(p_l)$ and therefore,
\begin{align*}
R_{2k}=h_0\sum_{l=2k-1}^{2k}\EE\left[\int_0^{T_k} (W_l(t)-w_l) dt\right].
\end{align*}
Roughly speaking, $R_{2k}$ depends on how fast $W_l(t)$ converges to its steady state for given $(\mu_l,p_l)$. Given the ergodicity convergence result in Lemma \ref{lmm: ergodicity}, we can show that $W_l(t)$ becomes close to the steady-state distribution after a  warm-up period of length $t_k=O(\log(k))$.
\begin{lemma}[Nonstationary Error after Warm-up] \label{coro: regret of nonstationarity after warm-up period}Suppose $T_k>t_k \equiv \log(k)/\gamma$, then
$$\EE\left[\int_{t_k}^{T_k} (W_l(t)-w_l) dt\right]= O(k^{-1}).$$
\end{lemma}
To obtain a finer bound for small values of $t$, i.e., in the warm-up period, we  follow a similar idea as in  \cite{ChenLiuHong1} and decompose $\EE[W_l(t)-w_{l}] = \EE[W_l(t)-w_{l-1}]+\EE[w_{l-1}-w_{l}]$.  
\begin{lemma}[Nonstationary Error in Warm-up Period]\label{lmm: regret of nonstationarity in warm-up period}
Suppose $T_k>t_k \equiv \log(k)/\gamma$ for all $k\geq 1$. % and $\max(\eta_k\sqrt{\mathcal{V}_k},\delta_k)=O(k^{-\xi})$ for some $0<\xi\leq 1$. 
Then, there exists a constant $C_0$ such that for all $l=2k-1, 2k$,
\begin{enumerate}
	\item[$(a)$] $\EE[|w_l-w_{l-1}|]\leq C_0\EE[\|\bx_l-\bx_{l-1}\|];$
	\item[$(b)$] $\EE\left[\int_0^{t_k}W_l(t)-w_{l-1}dt\right]\leq C_0\EE[\|\bx_l-\bx_{l-1}\|^2]^{1/2}t_k.$
\end{enumerate}
As a consequence,
$$\EE\left[\int_{0}^{t_k} (W_l(t)-w_l) dt\right]= O\left(\max(\eta_k\sqrt{\mathcal{V}_k}, \delta_k)\log(k)\right).$$
\end{lemma}
%Then, the proof of Proposition \ref{prop: regret of nonstationarity} follows from Lemma \ref{coro: regret of nonstationarity after warm-up period} and \ref{lmm: regret of nonstationarity in warm-up period}. 
Following Lemma \ref{coro: regret of nonstationarity after warm-up period} and Lemma \ref{lmm: regret of nonstationarity in warm-up period}, we have
\begin{align*}
R_{2k}&=h_0\sum_{l=2k-1}^{2k}\EE\left[\int_{0}^{t_k}W_l(t)-w_ldt+\int_{t_k}^{T_k}W_l(t)-w_ldt\right]
=O(k^{-1})+O(\max(\eta_k\sqrt{\mathcal{V}_k},\delta_k)\log(k))\\
&=O(k^{-1})+O(k^{-\xi}\log(k))=O(k^{-\xi}\log(k)).
\end{align*}
Furthermore, if $\eta_k=O(k^{-a}),T_k=O(k^b)$ and $\delta_k=O(k^{-c})$, then by  Proposition \ref{prop: bound Bk V_k}, $\eta_k\sqrt{\mathcal{V}_k}=O(k^{\max(-a-b/2+c,-a)})$. As a result, $\max(\eta_k\sqrt{\mathcal{V}_k},\delta_k)=O(k^{\max(-a-b/2+c,-a,-c)})$. Therefore, setting $\xi=\max(-a-b/2+c,-a,-c)$ finishes the proof.\hfill$\Box{}$

\subsection{Proof of Theorem \ref{thm: upper bound}}\label{subsec: proof of regret bound}
As discussed in Section \ref{subsec: regret upper bound}, the bound for regret of suboptimality $R_{1k}$ follows immediately from Theorem \ref{thm: convergence}. The bound for $R_{2k}$ follows from Proposition \ref{prop: regret of nonstationarity}. The bound for $R_{3k}$ follows from the smooth condition in Assumption \ref{assmpt: uniform}. 
\begin{lemma}[Exploration Cost]
\label{lmm: exploration cost}
Under Assumption \ref{assmpt: uniform}, there exists a constant $K_4>0$ such that
\begin{equation}
	\label{eq: R_3k}
	R_{3k}\leq K_4T_k\delta_k^2.
\end{equation}
\end{lemma}

Now, given that  $\eta_k=c_\eta k^{-1}$ with $c_\eta>2/K_0$, $T_k =c_Tk^{1/2}$ with $c_T>0$ and $\delta_k=c_\delta k^{1/3}$ with $0<c_\delta<\sqrt{K_0/32c}$, by Proposition \ref{prop: bound Bk V_k},
$$B_k\leq 2c\delta_k^2 + O(\delta_k^{-1}\exp(-\theta_1\alpha T_k))= \frac{K_0}{16}k^{-2/3}+o(k^{-2/3})\leq \frac{K_0}{8}k^{-2/3},$$
for $k$ large enough, and $\mathcal{V}_k=O(k^{1/3})$. So condition \eqref{eq: beta conditions} is satisfied with $\beta=2/3$  and hence $R_{1k}= O(k^{-1/3})$. On the other hand, conditions in Proposition \ref{prop: regret of nonstationarity} hold with $\xi = 1/3$ and hence $R_{2k}= O(k^{-1/3}\log(k))$. Finally, $R_{3k} = O(T_k\delta_k^2)=O(k^{-1/3})$. So we can conclude that
$$R(L) = \sum_{k=1}^L(R_{1k}+R_{2k}+R_{3k})=\sum_{k=1}^LO(k^{-1/3}\log(k))=O(L^{2/3}\log(L)).$$
As $T_k =O(k^{1/3})$, we have $T(L)=O(L^{4/3})$, and therefore $R(L)=O(\sqrt{T(L)}\log(T(L))$.\hfill$\Box{}$
%%%%%%%%%%%%%%%%%%%%%%%%%%%%%%%%%%%%%%%%%%%%%%%%
%%%%%%%%%%%%%%%%%%%%%%%%%%%%%%%%%%%%%%%%%%%%%%%%
\section{Additional Proofs}\label{appx: proofs}
\subsection{Full Proof of Theorem \ref{thm: convergence}}\label{appx: proof of convergence theorem}
By the SGD recursion,
$\bar{\bx}_{k+1} = \Pi_{\mathcal{B}}(\bar{\bx}_k -\eta_k \bH_k)$.
Let $\mathcal{F}_k$ be the filtration up to iteration $k$, i.e. it includes all events in the first $2(k-1)$ cycles. By Lemma \ref{lmm: convexity}, we have
\begin{align*}
&\EE\left[\|\bar{\bx}_{k+1}-\bx^*\|^2]\leq\EE[\|\bar{\bx}_{k}-\bx^*-\eta_k \bH_k\|^2\right]\\
=~&\EE\left[\|\bar{\bx}_k-\bx^*\|^2 - 2\eta_k \bH_k\cdot( \bar{\bx}_k-x^*)+ \eta_k^2\|\bH_k\|^2\right]\\
=~&\EE\left[\|\bar{\bx}_k-\bx^*\|^2 -2\eta_k\nabla f(\bar{\bx}_k)\cdot (\bar{\bx}_k-\bx^*)\right]- \EE[2\eta_k(\bH_k-\nabla f(\bar{\bx}_k))\cdot( \bar{\bx}_k-\bx^*)] + \EE[\eta_k^2\|\bH_k\|^2]\\
\leq~&(1-2\eta_k K_0)\EE\left[\|\bar{\bx}_k-\bx^*\|^2\right] + \EE[2\eta_k(\bH_k-\nabla f(\bar{\bx}_k))\cdot(\bx^*-\bar{\bx}_k)] + \eta_k^2\EE[\|\bH_k\|^2].
\end{align*}
Note that
\begin{align*}
&\EE[2\eta_k(\bH_k-\nabla f(\bar{\bx}_k))\cdot( \bx^*-\bar{\bx}_k)] =~\EE[\EE[2\eta_k(\bH_k-\nabla f(\bar{\bx}_k))\cdot( \bx^*-\bar{\bx}_k)|\mathcal{F}_k]] \\
=~& 2\eta_k\EE[\EE[\bH_k-\nabla f(\bar{\bx}_k)|\mathcal{F}_k]\cdot( \bx^*-\bar{\bx}_k)] 
\leq~  2\eta_k \EE[\|\EE[\bH_k-\nabla f(\bar{\bx}_k)|\mathcal{F}_k]\|^2]^{1/2}\EE[\|\bx^*-\bar{\bx}_k\|^2]^{1/2}\\
\leq~& \eta_k \EE[\|\EE[\bH_k-\nabla f(\bar{\bx}_k)|\mathcal{F}_k]\|^2]^{1/2} (1+\EE[\|\bar{\bx}_k-\bx^*\|^2]).
\end{align*}
The second last inequality follows from 
%$ab+cd\leq \sqrt{a^2+c^2}\sqrt{b^2+d^2}$ and the 
H\"{o}lder's Inequality, and the last inequality follows from $2a\leq 1+a^2$.
Let  $b_k=\EE[\|\bar{\bx}_k-\bx^*\|^2]$ and recall that we have defined
$$B_k = \EE[\|\EE[\bH_k-\nabla f(\bar{\bx}_k)|\mathcal{F}_k]\|^2]^{1/2}, \quad \mathcal{V}_k = \EE[\|\bH_k\|^2].$$
Then, we obtain the recursion
\begin{equation}\label{eq: b_k recursion}
b_{k+1}\leq (1-2K_0\eta_k+\eta_kB_k)b_k + \eta_kB_k+\eta_k^2\mathcal{V}_k.
\end{equation}
Next, we prove by mathematical induction that there exists a large constant $K_2>0$ such that $b_k\leq K_2k^{-\beta}$ for all $k\geq 1$ using  recursion \eqref{eq: b_k recursion}. Given that $\eta_k\mathcal{V}_k=O(k^{-\beta})$, we can find a constant  $K_3>0$ large enough such that $\eta_k\mathcal{V}_k\leq K_3k^{-\beta}$ for all $k\geq 1$. Then, by the induction assumption that $b_k\leq K_2k^{-\beta}$, we have
$$
b_{k+1}\leq (1-2K_0\eta_k + \eta_kB_k)b_k +\eta_kB_k + \eta_k^2 \mathcal{V}_k\leq \left(1-2K_0\eta_k + \frac{K_0}{8}\eta_kk^{-\beta} \right)b_k +\frac{K_0}{8}\eta_kk^{-\beta} + K_3\eta_kk^{-\beta}.
$$
Note that $k^{-\beta}/(k+1)^{-\beta}= (1+\frac{1}{k})^\beta\leq 1+\frac{1}{k}\leq 1+\frac{K_0}{2}\eta_k$. So we have
%	$$b_{k+1}\leq (1-K_0\eta_k+K_2\eta_k\Delta_k)(1+A_k)C\Delta_{k-1}+K_2\eta_k\Delta_k + \eta_k^2K_2.$$
%	Since $\alpha\leq 0.5$ and we have chosen $\eta_k \geq 2/K_0k$, as a result, $1+A_k\leq 1+ K_0\eta_k/2$. Therefore,  we have
\begin{equation*}
\begin{aligned}
	b_{k+1}&\leq \left(1-2K_0\eta_k+\frac{K_0}{8}\eta_kk^{-\beta}\right)\left(1+\frac{K_0\eta_k}{2}\right)K_2(k+1)^{-\beta}+\frac{K_0}{8}\eta_kk^{-\beta} + K_3\eta_kk^{-\beta}\\
	&\leq K_2(k+1)^{-\beta} -\eta_kk^{-\beta}\left(\frac{3K_0K_2}{2}-\frac{K_0K_2}{8}k^{-\beta} - \frac{K_0^2K_2}{16}\eta_kk^{-\beta} -\frac{K_0}{8}-K_3\right).\\
	%& \leq  \Delta_k -\eta_k\Delta_k\left(\frac{K_0C}{2}-\frac{K_0}{4}-\frac{K_0C}{8}\Delta_k - K_0^2C\frac{\eta_k\Delta_k}{16} \right).
\end{aligned}
\end{equation*}
Then, we have $b_{k+1}\leq  K_2(k+1)^{-\beta}$
as long as
\begin{equation*}\label{eq: proof thm1}
\frac{3K_0K_2}{2}-\frac{K_0K_2}{8}k^{-\beta} - \frac{K_0^2K_2}{16}\eta_kk^{-\beta} -\frac{K_0}{8}-K_3 \geq 0.
\end{equation*}
As the step size $\eta_k\to 0$,  $\eta_kK_0\leq 1$ for $k$ large enough. Let $k_0=\max\{k\geq 1:\eta_kK_0> 1\}$. Then, if $K_2\geq 8K_3/K_0$, for all $k\geq k_0$,%To check \eqref{eq: proof thm1}, note that, $k^{-\beta}, K_0\leq 1$ and 
$$\frac{3K_0K_2}{2}-\frac{K_0K_2}{8}\Delta_k - \frac{K_0^2K_2}{16}\eta_k\Delta_k -\frac{K_0}{8}-K_3 \geq \frac{3K_0K_2}{2}-\frac{K_0K_2}{8}-\frac{K_0K_2}{16}-\frac{K_0K_2}{8}-\frac{K_0K_2}{8} = \frac{17K_0K_2}{16}>0.$$
Let
\begin{equation*}\label{eq: C}
K_2= \max\left(k_0^{\beta}(|\bar{\mu}-\underline{\mu}|^2+|\bar{p}-\underline{p}|^2), 8K_3/K_0\right).
\end{equation*}
Then we have $\|\bar{\bx}_k-\bx^*\|^2\leq K_2k^{-\beta}$ for all $1\leq k\leq k_0$, and we can conclude by induction that, for all $k\geq k_0$,
$$\EE[\|\bar{\bx}_k-\bx^*\|^2]\leq K_2k^{-\beta}.$$
%Now, if, additionally, $\eta_k=O(k^{-a}),\delta_k=O(k^{b})$ and $\delta_k=O(k^{-c})$, then by Proposition \ref{prop: bound Bk V_k}, $B_k=O(k^{-2c})$ and $\eta_k\mathcal{V}_k=O(k^{\max(2c-b-a,-a)})$. Choosing $\beta=\max(2c-b-a,-a,-2c)$ leads to the ideal result.
%%%%%%%%%%%%%From x to Regret, now combined in sec.5.2%%%%%%%%%%%%%%%%%%%%%%%%
%			By Lemma \ref{lmm: convexity}, there exists $a\in [0,1]$ such that for $k\geq k_0$,
%			$$|f(\bar{x}_k)-f(x^*)|=|\nabla f(a(\bar{x}^k-x^*)+x^*)^T(\bar{x}_k-x^*)|\leq K_1\|\bar{x}_k-x^*\|^2.$$
%			As a consequence,
%			$$R_{2k}\leq 2T_kK_1K_2k^{-\beta}=O(k^{-\beta}T_k).$$  \hfill $\Box$
%%%%%%%%%%%%%%%%%%%%%%%%%%%%%%%%%%%%%%%%%%%%
\hfill $\Box$

%\subsection{Uniform Moment Bounds}\label{appx: uniform bound}
\subsection{Proofs of Technical Lemmas}\label{appx: lemma proofs}
In addition to the uniform moment bounds for $W_l(t)$ as stated in Lemma \ref{lmm: uniform bound short}, we also need to establish similar bounds for the  so-called observed busy period $X_l(t)$, which will be used in the proof of Lemma \ref{lmm: regret of nonstationarity in warm-up period}. In detail, $X_l(t)$ is the units of time that has elapsed at time point $t$ in cycle $l$ since the last time when the server is idle (probably in a previous cycle). So the value of $X_l(t)$ is uniquely determined by $\{W_l(t)\}$, i.e., $X_l(t)=0$ whenever $W_l(t)=0$ and $dX_l(t) = dt$ whenever $W_l(t)>0$.
\begin{lemma}[Complete Version of Lemma \ref{lmm: uniform bound short}]\label{lmm: uniform bound} Under Assumptions \ref{assmpt: uniform} and \ref{assmpt: light tail}, there exist some constants $\theta_0>0$ and  $M>1$ such that, for any sequence of control parameters $\{(\mu_l,p_l):l\geq 1\}$, 
$$\EE[X_l^m(t)]
\leq M,\quad \EE[W_l(t)^m]\leq M, \quad\EE[W_l(t)^m\exp(2\theta_0 W_l(t))]\leq M,$$
for all $m\in\{0, 1,2\}$, $l\geq 1$ and $0\leq t\leq T_{k}$ with $k=\lceil l/2\rceil$.
\end{lemma}
\begin{proof}{Proof of Lemma \ref{lmm: uniform bound}}
We consider a $M/GI/1$ system under a stationary policy such that $\mu_l\equiv \underline{\mu}$  and $p_l\equiv \underline{p}$ for all $l\geq 1$. We call this system the dominating system and denote its workload process by $W^D_l(t)$. In addition, we set $W_1^D(0)\stackrel{d}{=} W_\infty(\underline{\mu},\underline{p})$ so that $W_l^D(t)\stackrel{d}{=} W_\infty(\underline{\mu},\underline{p})$ for all $l\geq 1$ and $t\in[0,T_k]$. Then, the arrival process in the dominating system is an upper envelop process (UEP) for all possible arrival processes corresponding to any control sequence $(\mu_l,p_l)$ and the service process in the dominating system is a lower envelope process (LEP) for all possible service processes corresponding to any control sequence. In addition, $W_1(0)=0\leq W^D_l(t)$. So we have 
$$W_l(t)\leq_{st} W^D_l(t)\stackrel{d}{=}W_\infty(\underline{\mu},\underline{p}), \text{ for all }l\geq 1\text{ and }t\in[0,T_k].$$
By Theorem 5.2 in the Chapter X of \cite{AsmussenBook}, the stationary workload process $$W_\infty(\underline{\mu},\underline{p})\stackrel{d}{=} Y_1+...+Y_N.$$
Here $N$ is  geometric random variable of mean $1/(1-\bar{\rho})$ and $\bar{\rho} = \lambda(\underline{p})/\underline{\mu}$, and $Y_n$ are I.I.D. random variables independent of $N$. In addition, the density of $Y_n$ is
$$f_Y(t) = \frac{\PP(V_n>t)}{\EE[V_n]},\quad t\in [0,\infty).$$
Under Assumption \ref{assmpt: light tail}, we have
$$\PP(Y_n>t) = \int_t^{\infty} f_Y(s)ds=\int_t^{\infty}\frac{\PP(V_n>s)}{\EE[V_n]}ds\leq \int_t^\infty \frac{\exp(-\eta s)\EE[\exp(\eta V_n)]}{\EE[V_n]}ds=\frac{\EE[\exp(\eta V_n)]}{\eta\EE[V_n]}\cdot\exp(-\eta t).$$
As a consequence, $Y_n$ has finite moment generating function around the origin. As $W_\infty(\underline{\mu},\underline{p})$ is a geometric compound of $Y_n$, it also has finite moment generating function around the origin. So we can conclude that, there exists some constants $\theta_0\in(0,\theta/2)$ and $C\geq1$ such that
$$\EE[W_l(t)^m]\leq \EE[W_\infty(\underline{\mu},\underline{p})^m]\leq C,\quad \EE[W_l(t)^m\exp(2\theta_0 W_l(t))]\leq \EE[W_\infty(\underline{\mu},\underline{p})^m\exp(2\theta_0 W_\infty(\underline{\mu},\underline{p}))]\leq C,$$
for $m = 1, 2$.

To deal with the observed busy period, we need to do a time-change. In detail, for each cycle $l$ and control parameter $(\mu_l,p_l)$, we ``slow down" the clock by $\lambda(p_l)$ times so that the arrival rate is normalized to 1 and mean service time to $\lambda(p_l)/\mu_l$. We denote the time-changed workload and observed busy period by $\tilde{W}_l(t)$ and $\tilde{X}_l(t)$ for $t\in [0,\lambda(p_l)T_k]$. Then, for all $t\in[0,T_k]$,
$$W_l(t)\leq \frac{1}{\lambda(\bar{p})}\tilde{W}_l\left(\lambda(p_l)t\right),\quad X_l(t)\leq \frac{1}{\lambda(\bar{p})}\tilde{X}_l\left(\lambda(p_l)t\right).$$
We denote by $\tilde{X}^D_l(t)$ the time-changed observed busy period corresponding to the dominating system. Then, since $\lambda(p_l)/\mu_l/\leq \lambda(\underline{p})/\underline{\mu}$ for all possible values of $(\mu_l,p_l)$, we can conclude that $\tilde{X}_l(t)\leq_{st}\tilde{X}^D_l(t)$. Following \cite{nakayama2004finite}, $\EE[\tilde{X}^D_l(t)]\leq \EE[X_\infty(1,\underline{\mu}/\lambda(\underline{p}))]<\infty$. Let $M=C\vee \left(\EE[X_\infty(1,\underline{\mu}/\lambda(\underline{p}))]/\lambda(\bar{p})\right)$ and we can conclude that 
$\EE[X_l(t)]\leq M$.
\hfill $\Box{}$
\end{proof}
\vspace{1ex}
%\subsection{Proof of Proposition \ref{prop: delay observation}}\label{appx: delay observation}
%\begin{proof}{Proof of Proposition \ref{prop: delay observation}}
%	Following \eqref{eq: hat W},
%$$\EE_l\left[|\hat{W}_l(t)-W_l(t)|\right]=\EE_l\left[W_l(t)\cdot\ind{t>A_{N^D(l)}^{C^D(l)}}\right].$$
%By definition, the event $\{t\leq A_{N^D(l)}^{C^D(l)}\}$ happens as long as the workload $W_l(t)$ are all served by $T_k$, i.e. $W_l(t)\leq \mu_l(T_k-t)$. Consequently,
%	\begin{align*}	%\EE_l\left[|\hat{W}_l(t)-W_l(t)|\right]&=\EE_l\left[W_l(t)\cdot\ind{t>A_{N^D(l)}^{C^D(l)}}\right]=
%\EE_l\left[W_l(t)\cdot\ind{W_l(t)>\mu_l(T_k-t)}\right]&\leq \EE_l\left[W_l(t)^2\right]^{1/2}\PP_l\left(W_l(t)>\mu_l(T_k-t)\right)^{1/2}\\
%&\leq \EE_l\left[W_l(t)^2\right]^{1/2}\cdot\exp\left(-\frac{1}{2}\theta_0\mu_l(T_k-t)\right)\EE_l\left[\exp(\theta_0W_l(t))\right]^{1/2}.
%\end{align*}
%As a result, by  Lemma \ref{lmm: uniform bound},
%$$\EE\left[|\hat{W}_l(t)-W_l(t)|\right]\leq \exp\left(-\frac{1}{2}\theta_0\underline{\mu}(T_k-t)\right)M.$$
%\hfill $\Box$
%\end{proof}

%\subsection{Proofs of Auxiliary Results in Section \ref{subsec: proof of estimation error}}\label{appx: proof of estimation error}
\begin{proof}{Proof of Lemma \ref{lmm: ergodicity}}
	Let $N(t)$ be the arrival process under control parameter $(\mu,p)$, which is a Poisson process with rate $\lambda(p)$. Define an auxiliary L\'evy process as $R(t)=\sum_{i=1}^{N(t)}V_i-\mu t$. For the workload processes $W(t)$ and $\bar{W}(t)$, define two hitting times $\tau$ and $\bar{\tau}$  as
	$$\tau\equiv \min_{t\geq 0}\{t: W(0)+R(t)=0\},\quad\text{and}\quad \bar{\tau}\equiv \min_{t\geq 0}\{t: \bar{W}(0)+R(t)=0\}.$$
	%Because $W(t)$ and $\bar{W}(t)$ share the same arrivals and individual workloads, one can check that the difference
	%$|W(t)-\bar{W}(t)|$
	%is non-increasing in $t$. Moreover, for any $t>\max(\tau,\bar{\tau})$,
	%$$0\leq |W(t)-\bar{W}(t)|\leq |W(\tau\vee\bar{\tau})-\bar{W}(\tau\vee\bar{\tau})|=0.$$
	Following Lemma 2 of \cite{ChenLiuHong1}, we have 
	\begin{equation}
		\label{eq: pathwise w-bar_w}
		|W(t)-\bar{W}(t)|\leq |W(0)-\bar{W}(0)|\ind{t< \tau\vee\bar{\tau}}.
	\end{equation}
	Next, we give a bound for the probability $\PP(\tau>t)$ by constructing an exponential supermartingale. Define
	$$M(t)=\exp\left(\theta_0(W(0)+R(t))+\gamma t\right),$$
	where $\theta_0$ is defined in Lemma \ref{lmm: uniform bound} and the value of $\gamma$ will be specified in \eqref{eq: gamma choice}.  Let $\{\mathcal{F}_t\}_{t\geq 0}$ be the natural filtration associated to $R(t)$.
	%As $R(t)=\sum_{i=1}^{N(t)}V_i-\mu t$ is a compound Poisson process with negative drift, 
	For any $t,s>0,$
	\begin{align*}
		\EE[M(t+s)|\mathcal{F}_t]&=\EE[M(t)\exp(\theta_0(R(t+s)-R(t))+\gamma s)|\mathcal{F}_t]=M(t)\EE[\exp(\theta_0 R(s)+\gamma s)]\\
		&=M(t)\EE\left[\exp \left(\theta_0\sum_{i=1}^{N(s)}V_i -\theta_0 \mu s +\gamma s\right)\right]=M(t)\EE\left[\EE[\exp(\theta_0 V_i)]^{N(s)}\right]e^{-\theta_0 \mu s+\gamma s}\\
		&=M(t)\exp\left( s\left(\lambda \EE[\exp(\theta_0 V_i)]-\lambda -\mu\theta_0+\gamma\right)\right).
	\end{align*}
	%Next, we specify the value of $\gamma>0$ such that $\lambda \EE[\exp(\theta_0 V_i)]-\lambda -\mu\theta_0+\gamma<0$. Recall that, 
	According to Assumption \ref{assmpt: light tail}, $\phi(\theta)<\log(1+\underline{\mu}\theta/\bar{\lambda})-\gamma_0$ for some $\theta, \gamma_0>0$. Besides, the function	$h(x)\equiv \phi(x)-\log(1+\underline{\mu}x/\bar{\lambda})$
	is convex on $[0,\theta]$. As $0<\theta_0<\theta$, we have
	$$h(\theta_0)\leq (1-\theta_0/\theta) h(0)+ \frac{\theta_0}{\theta}h(\theta)<-\frac{\theta_0}{\theta}\gamma_0.$$
	We choose
	\begin{equation} \label{eq: gamma choice}
		\gamma=\underline{\lambda}\left(1-e^{-\frac{\theta_0\gamma_0}{\theta}}\right)\left(1+\underline{\mu}\theta_0/\bar{\lambda}\right).
	\end{equation}
	Then, it satisfies that
	\begin{align*}
		\lambda \EE[\exp(\theta_0 V_i)]-\lambda -\mu\theta_0+\gamma&=\lambda\left( e^{\phi(\theta_0)}-(1+\frac{\mu\theta_0}{\lambda})+\frac{\gamma}{\lambda}\right)<\lambda\left(e^{-\frac{\theta_0}{\theta}\gamma_0}(1+\underline{\mu}\theta_0/\bar{\lambda})-(1+\mu\theta_0/\lambda)+\frac{\gamma}{\lambda}\right)\\
		&<\lambda\left(-\left(1-e^{\frac{\theta_0\gamma_0}{\theta}}\right)(1+\underline{\mu}\theta_0/\bar{\lambda})+\frac{\gamma}{\underline{\lambda}}\right)=0.
	\end{align*}
	Now, we can conclude that $M(t)$ is an non-negative supermartingale with $\gamma$ as given by \eqref{eq: gamma choice}. By Fatou's lemma, 
	\begin{align*}
		\PP(\tau>t|W(0))&\leq e^{-\gamma t}\EE[\exp(\gamma \tau)|W(0)]=e^{-\gamma t}\EE[\liminf_{n\rightarrow\infty}M(\tau\wedge n)| W(0)]\\
		&\leq e^{-\gamma t}\liminf_{n\rightarrow\infty}\EE[M(\tau\wedge n) | W(0)]\leq e^{-\gamma t}\EE[M(0)|W(0)]=e^{-\gamma t}\exp(\theta_0 W(0)).
	\end{align*}
	Similarly, $\PP(\bar{\tau}>t|\bar{W}(0))\leq e^{-\gamma t} \exp(\theta_0\bar{W}(0))$. Combining these bounds with \eqref{eq: pathwise w-bar_w}, we can conclude that
	\begin{align*}
		\EE\left[|W(t)-\bar{W}(t)|^m|W(0),\bar{W}(0)\right]&\leq |W(0)-\bar{W}(0)|^m \PP(\tau\vee\bar{\tau}>t|W(0),\bar{W}(0))\\
		&\leq |W(0)-\bar{W}(0)|^m \left(\PP(\tau>t|W(0))+\PP(\bar{\tau}>t|\bar{W}(0))\right)\\
		&\leq |W(0)-\bar{W}(0)|^m \left(e^{\theta_0 W(0)+\theta_0 \bar{W}(0)}\right)e^{-\gamma t}.
	\end{align*}
	\hfill $\Box{}$
\end{proof}

\begin{proof}{Proof of Lemma \ref{lmm: autocovariance W}}
	We first analyze the conditional expectation
	$
	\EE_l[({W}_l(t)-w_l)({W}_l(s)-w_l)]
	$
	for each given pair of $(s,t)$ such that $0\leq  s\leq t\leq T_k$. To do this, we synchronously couple with $\{W_l(r): s\leq r\leq T_k\}$ a stationary workload process $\{\bar{W}_l^s(r):s\leq r\leq T_k\}$. In particular, $\bar{W}_l^s(s)$ is independently drawn from the stationary distribution $W_\infty(\mu_l,p_l)$. As a result, $\bar{W}_l^s(r)$ is independent of $W_l(s)$ for all $s\leq r\leq T_k$, and hence
	$$ \EE_l[W_l(s)(\bar{W}^s_l(t)-w_l)]=\EE_l[W_l(s)]\left(\EE_l[\bar{W}^s_l(t)]-w_l\right)=0.$$
	Then, we have
	\begin{align*}
		&\EE_l[(W_l(t)-w_l)(W_l(s)-w_l)] =\EE_l[(W_l(t)-\bar{W}^s_l(t))W_l(s)] - w_l\EE_l[W_l(s)-w_l]. 
	\end{align*}
	By Lemma \ref{lmm: ergodicity}, 
	$$\EE_l[(W_l(t)-\bar{W}^s_l(t))W_l(s)|W_l(s),\bar{W}_l^s(s)]\leq \exp(-\gamma(t-s))(e^{\theta_0 W_l(s)}+e^{\theta_0\bar{W}^s_l(s)})(W_l(s)+\bar{W}^s_l(s))W_l(s).$$
	As $\bar{W}_l^s(s)$ is independent of $W_l(s)$,
	\begin{align*}
		&\EE_l[(W_l(t)-\bar{W}^s_l(t))W_l(s)|W_l(s)]\\
		\leq~&\exp(-\gamma(t-s))\EE_l\left[(e^{\theta_0 W_l(s)}+e^{\theta_0\bar{W}^s_l(s)})(W_l(s)+\bar{W}^s_l(s))W_l(s)|W_l(s)\right]\\
		=~&\exp(-\gamma(t-s))(e^{\theta_0 W_l(s)}W_l(s)^2 +e^{\theta_0 W_l(s)}W_l(s)\EE[\bar{W}^s_l(s)]+W_l(s)^2\EE[e^{\theta_0\bar{W}^s_l(s)}]+W_l(s)\EE[e^{\theta_0\bar{W}^s_l(s)}\bar{W}^s_l(s)] )\\
		\leq ~&\exp(-\gamma(t-s))(e^{\theta_0 W_l(s)}W_l(s)^2 +Me^{\theta_0 W_l(s)}W_l(s)+MW_l(s)^2+MW_l(s)).
	\end{align*}
	One can check that $W_l(s)\leq W_l(0)+\bar{W}_l(s)$, where $\bar{W}_l(s)$ is a stationary workload process synchronously coupled with $W_l(t)$ having an independent drawn initial $\bar{W}_l(0)$. Therefore,  
	\begin{align*}
		\EE_l\left[e^{\theta_0W_l(s)}W_l(s)^2\right]&\leq e^{\theta_0 W_l(0)}\EE_l\left[(W_l(0)+\bar{W}_l(s))^2e^{\theta_0\bar{W}_l(s)}\right]\\
		&=e^{\theta_0 W_l(0)}\left(W_l(0)^2\EE_l[e^{\theta_0\bar{W}_l(s)}]+2W_l(0)\EE_l\left[\bar{W}_l(s)e^{\theta_0\bar{W}_l(s)}\right]+\EE_l\left[W_l(s)^2e^{\theta_0\bar{W}_l(s)}\right]\right)\\
		&\leq 2Me^{\theta_0W_l(0)}(1+W_l(0)^2),\\
		\EE_l\left[e^{\theta_0W_l(s)}W_l(s)\right]&\leq e^{\theta_0W_l(0)}\EE_l\left[W_l(0)e^{\theta_0\bar{W}_l(s)}+\bar{W}_l(s)e^{\theta_0\bar{W}_l(s)}\right]\leq e^{\theta_0 W_l(0)}M(1+W_l(0))\\
		&\leq \frac{3M}{2}e^{\theta_0W_l(0)}(1+W_l(0)^2),
	\end{align*}	
	where the last inequality holds because the constant $M\geq 1$ and $W_l(0)\leq (1+W_l(0)^2)/2$. Note that $W_l(s)^2\leq e^{\theta_0 W_l(s)}W_l(s)^2$ and $W_l(s)\leq W_l(s)e^{\theta_0W_l(s)}$, we have 
	$$\EE_l[(W_l(t)-\bar{W}_l^s(t))W_l(s)]\leq e^{-\gamma (t-s)}e^{\theta_0W_l(0)}(1+W_l(0)^2)(2M+5M^2).$$
	On the other hand, by Lemma \ref{lmm: ergodicity},
	\begin{align*}
		|\EE_l[W_l(s)-w_l]|&\leq \exp(-\gamma s)MW_l(0)(M+W_l(0))\exp(\theta_0 W_l(0))\\
		&\leq e^{-\gamma s}e^{\theta_0W_l(0)}M^2(1+W_l(0))^2\leq 2M^2e^{-\gamma s}e^{\theta_0W_l(0)}(1+W_l(0)^2).
	\end{align*}
	As a consequence,
	\begin{align*}
		\EE_l[(W_l(t)-w_l)(W_l(s)-w_l)]&=\EE_l[(W_l(t)-\bar{W}^s_l(t))W_l(s)] - w_l\EE_l[W_l(s)-w_l]\\
		&\leq (e^{-\gamma (t-s)}+e^{-\gamma s})e^{\theta_0W_l(0)}(1+W_l(0)^2)(2M+5M^2+2M^3).
	\end{align*}
	%	\begin{align*}
		%		\EE[(W_l(t)-w_l)(W_l(s)-w_l)]=~ &\EE[\EE_l[(W_l(t)-w_l)(W_l(s)-w_l)]]\\
		%		\leq ~&\EE[\EE_l[\exp(-\gamma(t-s))(e^{\theta_0 W_l(s)}W_l^2(s) +Me^{\theta_0 W_l(s)}W_l(s)+MW_l(s)^2+MW_l(s))]\\
		%		&+\exp(-\gamma s)MW_l(0)(M+W_l(0))\exp(\theta_0 W_l(0))]\\
		%		=~&\EE[\exp(-\gamma(t-s))(e^{\theta_0 W_l(s)}W_l^2(s) +Me^{\theta_0 W_l(s)}W_l(s)+MW_l(s)^2+MW_l(s))]\\
		%		&+\EE[\exp(-\gamma s)MW_l(0)(M+W_l(0))\exp(\theta_0 W_l(0))]
		%	\end{align*}
	%	By Lemma \ref{lmm: uniform bound}, there exists a universal constant $K_V>0$ such that
	%	$$\EE[e^{\theta_0 W_l(s)}W_l^2(s) +Me^{\theta_0 W_l(s)}W_l(s)+MW_l(s)^2+MW_l(s)], \EE[MW_l(0)(M+W_l(0))\exp(\theta_0 W_l(0))]\leq K_V.$$
	and we can conclude \eqref{eq: covariance bound} with $K_V=2M+5M^2+2M^3$.
	\hfill $\Box{}$
	%	Recall that $\bar{W}_l(t)$  the stationary workload process synchronously coupled with $W_l(t)$ since $t=0$. We have $W_l(s)\leq W_l(0)+\bar{W}_l(s)$. We have
	%	$$\EE_l[(W_l(t)-\bar{W}^s_l(t))W_l(s)|W_l(s)]\leq \red{\exp(-\gamma(t-s))M^2(W_l(0)+M)(W_l(0)+2M)\exp(\eta W_l(0)).}$$
	%	\begin{center}
		%		\red{Assuming $\EE[\bar{W}_l(0)^2\exp(\eta\bar{W}_l(0))]\leq M$.}
		%	\end{center}
	%	
	%	Finally, we have proved in Lemma \ref{lmm: bias} that
	%	
	%	
	%	In Summary, we can conclude that there exists a universal  $K_V>0$ such that 
	%	$$\EE_l[(W_l(t)-w_l)(W_l(s)-w_l)]\leq K_V(\exp(-\gamma(t-s))+\exp(-\gamma s))W_l(0)^2\exp(\eta W_l(0)).$$
\end{proof}
\vspace{2ex}
%\subsection{Proofs of Auxiliary Results in Section \ref{subsec: proof of convergence}}\label{appx: proofs of convergence}

\begin{proof}{Proof of Lemma \ref{lmm: fd approximation}}
	By the mean value theorem, 
	\begin{small}
	\begin{align*}
		f(\mu_1,p) = f\left(\frac{\mu_1+\mu_2}{2},p\right) +\frac{\mu_1-\mu_2}{2}\partial_\mu f\left(\frac{\mu_1+\mu_2}{2},p\right)
		+\frac{(\mu_1-\mu_2)^2}{8}\partial^2_\mu f\left(\frac{\mu_1+\mu_2}{2},p\right)+\frac{(\mu_1-\mu_2)^3}{48}\partial^3_\mu f\left(\xi_1,p\right)\\
		f(\mu_2,p) = f\left(\frac{\mu_1+\mu_2}{2},p\right) +\frac{\mu_2-\mu_1}{2}\partial_\mu f\left(\frac{\mu_1+\mu_2}{2},p\right)
		+\frac{(\mu_1-\mu_2)^2}{8}\partial^2_\mu f\left(\frac{\mu_1+\mu_2}{2},p\right)+\frac{(\mu_2-\mu_1)^3}{48}\partial^3_\mu f\left(\xi_2,p\right),
	\end{align*}
	\end{small}
	where $\xi_1$ and $\xi_2$ take values between $\mu_1$ and $\mu_2$. As a consequence, we have
	$$	\left\vert\frac{f(\mu_1,p)-f(\mu_2,p)}{\mu_1-\mu_2}-\partial_\mu f\left(\frac{\mu_1+\mu_2}{2},p\right)\right\vert\leq c(\mu_1-\mu_2)^2,$$
	with $c=(\max_{(\mu,p)\in \mathcal{B}}|\partial_\mu^3 f(\mu,p)|\vee |\partial_p^3 f(\mu,p)|)/24$. Following the same argument, we have
	$$\left\vert\frac{f(\mu,p_1)-f(\mu,p_2)}{p_1-p_2}-\partial_\mu f\left(\mu, \frac{p_1+p_2}{2}\right)\right\vert\leq c(p_1-p_2)^2.$$
	\hfill$\Box{}$
\end{proof}
\vspace{2ex}

\begin{proof}{Proof of Lemma \ref{lmm: convexity}}
	By Pollaczek-Khinchin formula and PASTA, 
	%$$\EE[W_\infty(\mu,p)]=\frac{\lambda(p)(1+c_V^2)}{2(\mu-\lambda(p))},$$
	%and so
	$$f(\mu,p) = \frac{h_0(1+c_V^2)}{2}\cdot\frac{\lambda(p)}{\mu-\lambda(p)}+ c(\mu) -p\lambda(p).$$ We intend to show that $f(\mu,p)$ is strongly convex in $\mathcal{B}$. For ease of notation, denote $C=\frac{1+c_V^2}{2}$ and
	$$g(\mu,\lambda)=\frac{\lambda}{\mu-\lambda}.$$
	Write $\lambda(p)$,$\lambda'(p)$ and $\lambda''(p)$ as $\lambda$, $\lambda'$ and $\lambda''$ respectively. By direct calculation, we have 
	$$\partial_\lambda g=\frac{\mu}{(\mu-\lambda)^2},\partial_\mu g=\frac{\lambda}{(\mu-\lambda)^2}, \partial^2_{\lambda\lambda}g=\frac{2\mu}{(\mu-\lambda)^3}, \partial^2_{\lambda \mu}g=-\frac{\mu+\lambda}{(\mu-\lambda)^3},\partial^2_{\mu\mu}g=\frac{2\lambda}{(\mu-\lambda)^3}.$$
	The second-order derivatives are
	\begin{align*}
		\partial_{pp} f&=\frac{h_0C\mu}{(\mu-\lambda)^3}\left(2(\lambda')^2+(\mu-\lambda)\lambda''\right)-p\lambda''-2\lambda'\\
		\partial_{p\mu}f&=-\frac{h_0C (\mu+\lambda)}{(\mu-\lambda)^3},\quad \partial_{\mu\mu}f=\frac{2h_0C\lambda}{(\mu-\lambda)^3}+c''(\mu).
	\end{align*}
	By Condition (a) of Assumption \ref{assmpt: uniform}, we have $$-p\lambda''-2\lambda'>0\quad \text{and}\quad 2(\lambda')^2+(\mu-\lambda)\lambda''>0\quad\Rightarrow \quad \partial_{pp} f>0.$$
	It is easy to check that $\partial_{\mu\mu}f>0$ as $c(\mu)$ is convex. So, to verify the convexity of $f$, we only need to show that the determinant of Hessian metric $\bH_f$ is positive in $\mathcal{B}$. By direct calculation, 
	\begin{align*}
		|\bH_f|&=\frac{h_0^2C^2}{(\mu-\lambda)^5}\left(2\mu\lambda\lambda''-(\mu-\lambda)(\lambda')^2\right)+(-p\lambda''-2\lambda')\frac{2h_0C\lambda}{(\mu-\lambda)^3}+c''(\mu)\partial_{pp}f\\
		&\geq \frac{h_0^2C^2}{(\mu-\lambda)^5}\left(2\mu\lambda\lambda''-(\mu-\lambda)(\lambda')^2\right)+(-p\lambda''-2\lambda')\frac{2h_0C\lambda}{(\mu-\lambda)^3}\\
		&=\frac{h_0C}{(\mu-\lambda)^5}\left[h_0C(2\mu\lambda\lambda''-(\mu-\lambda)(\lambda')^2)+2\lambda(\mu-\lambda)^2(-p\lambda''-2\lambda')\right]\\
		&=-\frac{h_0C\lambda'}{(\mu-\lambda)^4}\left[h_0C\lambda'+4\lambda(\mu-\lambda)-2\frac{h_0C\mu-p(\mu-\lambda)^2}{\mu-\lambda}\frac{\lambda''\lambda}{\lambda'}\right].
	\end{align*}
	As $-\lambda'>0$, we need to prove the term in bracket is positive. Note that the term
	$$\frac{h_0C\mu-p(\mu-\lambda)^2}{\mu-\lambda}=h_0C+\frac{h_0C\lambda}{\mu-\lambda}-p(\mu-\lambda)$$
	is monotonically decreasing in $\mu$. By Assumption \ref{assmpt: uniform}, we have, for all $\mu\in[\underline{\mu},\bar{\mu}]$ and $\lambda\in[\underline{\lambda},\bar{\lambda}]$, 
	\begin{align*}
		&h_0C\lambda'+4\lambda(\mu-\lambda)-2\frac{h_0C\mu-p(\mu-\lambda)^2}{\mu-\lambda}\frac{\lambda''\lambda}{\lambda'}\\
		\geq & h_0C\lambda' +4\lambda(\underline{\mu}-\lambda)-2\left(h_0C+\frac{h_0C\lambda}{\mu-\lambda}-p(\mu-\lambda)\right)\frac{\lambda''\lambda}{\lambda'}\\
		\geq &h_0C\lambda'+4\lambda(\underline{\mu}-\lambda)-2h_0C \frac{\lambda''\lambda}{\lambda'}-2\max\left\{\left(\frac{h_0C\lambda}{\underline{\mu}-\lambda}-p(\underline{\mu}-\lambda)\right)\frac{\lambda''\lambda}{\lambda'},\left(\frac{h_0C\lambda}{\bar{\mu}-\lambda}-p(\bar{\mu}-\lambda)\right)\frac{\lambda''\lambda}{\lambda'}\right\}\\
		>&0.
	\end{align*}
	As $\mathcal{B}$ is compact, we can conclude that $f(\mu,p)$ is strongly convex on $\mathcal{B}$. Then by Taylor's expansion, Statement $(a)$ holds for some $1\geq K_0>0$. Statement (b) follows immediately after Assumption \ref{assmpt: uniform}.
	%		We want to show that $\partial^2_\mu f(\mu,p)>K_0$ and  $\partial^2_p f(\mu,p)>K_0$. 
	%		To compute  $\partial^2_\mu f(\mu,p)$, let
	%		$$g(\lambda) = \frac{\lambda C}{\mu-\lambda}, \text{ with }C=h_0(1+c^2_V)/2.$$
	%		We compute 
	%		$$\partial_\lambda g(\lambda) = \frac{C\mu}{(\mu-\lambda)^2}, \quad \partial^2_\lambda g(\lambda) = \frac{2C\mu}{(\mu-\lambda)^3}.$$
	%		Then, by Assumption \ref{assmpt: uniform}, 
	%		\begin{align*}
		%			\partial^2_p f(\mu, p)& =  \left(\partial^2_\lambda g(\lambda(p))\cdot (\lambda'(p))^2+\partial_\lambda f(\lambda(p))\cdot\lambda''(p)\right) - \left(2\lambda'(p) +p\lambda''(p)\right)\\
		%			& = \frac{C\mu}{(\mu-\lambda(p))^3} \left( 2(\lambda'(p))^2-(\mu-\lambda(p))\lambda''(p)\right) - \left(2\lambda'(p) +p\lambda''(p)\right)\\
		%			&>0.
		%		\end{align*}
	%		On the other hand,
	%		\begin{align*}
		%			\partial_\mu^2 f(\mu,p)&= \frac{2C\lambda(p)}{(\mu-\lambda(p))^3} + c''(\mu)>0.
		%		\end{align*}
	%		Since $\mathcal{B}$ is a closed set, there must exist $K_0>0$ such that
	%		$\partial^2_\mu f(\mu,p), \partial^2_p f(\mu,p)\geq K_0$ and hence statement (a) holds.  Statement (b) follows immediately after Assumption \ref{assmpt: uniform}.
	
	\hfill $\Box$
	%or equivalently
	%$$-\lambda'(p)> \max\left(\sqrt{\frac{0\vee\lambda''(p)}{2(\underline{\mu}-\bar{\lambda})}}~,~ \frac{p\lambda''(p)}{2}\right).$$
\end{proof}
\vspace{2ex}
\begin{proof}{Proof of Lemma \ref{coro: regret of nonstationarity after warm-up period}}
	By Lemma \ref{lmm: ergodicity}, conditional on $\mu_l, p_l$ and $W_l(0)$, we have
	\begin{align*}
		\EE_l[|W_l(t)-\bar{W}_l(t)|]&\leq \exp(-\gamma t)\EE_l\left[|W_l(0)-\bar{W}_l(0)|(\exp(\theta_0 W_l(0))+\exp(\theta_0 \bar{W}_l(0)))\right]\\
		&\leq \exp(-\gamma t) \left(W_l(0)\exp(\theta_0 W_l(0)) + MW_l(0)+M\exp(\theta_0 W_l(0)) + M\right)\\
		&\leq \exp(-\gamma t)M(M+W_l(0))\exp(\theta_0 W_l(0)).
	\end{align*}
	As a consequence, for $t\geq t_k$,
	\begin{align*}
		\EE[|W_l(t)-\bar{W}_l(t)|]&\leq\EE[\exp(-\gamma t)M(M+W_l(0))\exp(\theta_0W_l(0))]\\
		&=\exp(-\gamma t)\left(M^2\EE[\exp(\theta_0W_l(0))]+M\EE[W_l(0)\exp(\theta_0W_l(0))]\right)\leq \exp(-\gamma t)\cdot (M^2+M^3)%\leq k^{-1}(M^2+M^3).
	\end{align*}
	Therefore,
	\begin{align*}
		\EE\left[\int_{t_k}^{T_k}(W_l(t)-w_l)dt\right]&=\int_{t_k}^{T_k}\EE[W_l(t)-w_l]dt~\leq \int_{t_k}^{T_k}\EE[|W_l(t)-\bar{W}_l(t)|]dt\\
		&\leq \int_{t_k}^{T_k}\exp(-\gamma t)\cdot (M^2+M^3)dt
		~\leq  \exp(-\gamma t_k)\cdot\frac{M^2+M^3}{\gamma} \\
		&\leq k^{-1}\cdot\frac{M^2+M^3}{\gamma}=O(k^{-1}).
	\end{align*}
	\hfill$\Box{}$
\end{proof}
\vspace{2ex}

\begin{proof}{Proof of Lemma \ref{lmm: regret of nonstationarity in warm-up period}}
	
	Statement (1) is a direct corollary of Pollaczek–Khinchine formula. The proof of Statement (2) involves coupling workload processes with different parameters. Let us first explain the coupling in detail. Suppose $W^1(t)$ and $W^2(t)$ are two workload processes on $[0, T]$ with parameters $(\mu_1,\lambda_1)$ and $(\mu_2, \lambda_2)$ respectively. Let $W^1(0)$ and $W^2(0)$ be the given initial states. We construct two workload processes $\tilde{W}^1(t)$ and $\tilde{W}^2(t)$ on $[0,\infty)$ with parameters $(\mu_1/\lambda_1,1)$ and $(\mu_2/\lambda_2, 1)$ such that $\tilde{W}^i(0)=W^i(0)$ for $i=1,2$. The two processes $\tilde{W}^1(t)$ and $\tilde{W}^2(t)$ are coupled such that they share the same Poisson arrival process $N(t)$ with rate 1 and the same sequence of individual workload $V_n$.
	
	Then, we can couple $W^i(t)$ with $\tilde{W}(t)$ via a change of time, i.e. $W^i(t)=\tilde{W}^i(\lambda_it)$ and obtain
	$$\int_0^T W^i(t)dt = \frac{1}{\lambda_i}\int_0^{\lambda_i T}\tilde{W}^i(t) dt, \text{ for } i=1,2.$$
	Without loss of generality, assuming $\lambda_1\geq \lambda_2$ and we have
	\begin{align}\label{eq: couple different arrival rates}
		&\left|\int_{0}^T W^1(t)dt - \int_{0}^T W^2(t) dt~\right| \notag\\ 
		\leq & \frac{1}{\lambda_1}\left|\int_{0}^{\lambda_2 T}(\tilde{W}^1(t) - \tilde{W}^2(t))dt \right|+ \left|\frac{1}{\lambda_2} -\frac{1}{\lambda_1} \right|\int_{0}^{\lambda_2 T}\tilde{W}^2(t)dt + \frac{1}{\lambda_1}\int_{\lambda_2 T}^{\lambda_1 T}\tilde{W}^1(t)dt. 
	\end{align}
	Following a similar argument as in the proof of Lemma 3 in \cite{ChenLiuHong1}, we have that 
	$$|\tilde{W}^1(t)-\tilde{W}^2(t)|\leq \left|\frac{\mu_1}{\lambda_1}-\frac{\mu_2}{\lambda_2}\right|\max(\tilde{X}^1(t), \tilde{X}^2(t)) + |W^1(0) - W^2(0)|,$$
	where $\tilde{X}^i(t)$ is the observed busy period at time $t$, i.e. $$\tilde{X}^i(t) = t -\sup\{s: 0\leq s\leq t, \tilde{W}^i(s)=0\}.$$
	%	To bound the difference $\EE[W_l(t)]-w_l$,   we write
	%	$$\EE[W_l(t)]-w_l=(\EE[W_l(t)]-w_{l-1}) +(w_{l-1}-w_l).$$
	%	The second term is bounded by the smoothness of $\EE[W_\infty(\mu_l,p_l)]$ and we have
	%	$$\EE[|w_{l-1}-w_l|] \leq  C\EE\left[\|x_l-x_{l-1}\|\right]\leq C\max(\eta_k\sqrt{\mathcal{V}_k},\delta_k) .$$
	To apply \eqref{eq: couple different arrival rates} to bound $\EE[W_l(t)-w_{l-1}]$, we construct a stationary workload process $\bar{W}_{l-1}(t)$ with control parameter $(\mu_{l-1}, p_{l-1})$ synchronously coupled with $W_{l-1}(t)$ since the beginning of cycle $l-1$. In particular, $\bar{W}_{l-1}(0)$ is independently drawn from the stationary distribution of $W_\infty(\mu_{l-1},p_{l-1})$. We extend the sample path $\bar{W}_{l-1}(t)$ to cycle $l$, i.e. for $t\geq T_{k(l-1)}$ with $k(l-1)=\lceil (l-1)/2\rceil$, and couple it with $W_l(t)$ following the procedure described above. Then we have
	$$\EE\left[\int_0^{t_k}(W_l(t)-w_{l-1})dt\right]\leq \EE\left[\left|\int_0^{t_k} W_l(t)dt - \int_0^{t_k}\bar{W}_{l-1}(T_{k(l-1)}+t)dt~\right|\right].$$
	Without loss of generality, assume $\lambda_l\geq \lambda_{l-1}$. Then following \eqref{eq: couple different arrival rates}, we have 
	\begin{align*}
		&\left|\int_0^{t_k} W_l(t)dt - \int_0^{t_k}\bar{W}_{l-1}(T_{k(l-1)}+t)dt~\right|\\
		\leq~ & \frac{1}{\lambda_{l}}\left|\int_0^{\lambda_{l-1}t_k}(\tilde{W}_l(t)-\tilde{W}_{l-1}(T_{k(l-1)}+t))dt\right|+ \left|\frac{1}{\lambda_l} -\frac{1}{\lambda_{l-1}} \right|\int_{0}^{\lambda_{l-1} t_k}\tilde{W}_{l-1}(t)dt + \frac{1}{\lambda_{l}}\int_{\lambda_{l-1}t_k }^{\lambda_lt_k }\tilde{W}_l(t)dt \\
		\leq~ & \frac{1}{\lambda_{l}}\int_0^{\lambda_{l-1}t_k}\left|\tilde{W}_l(t)-\tilde{W}_{l-1}(T_{k(l-1)}+t)\right|dt+ \left|\frac{1}{\lambda_l} -\frac{1}{\lambda_{l-1}} \right|\int_{0}^{\lambda_{l-1} t_k}\tilde{W}_{l-1}(t)dt + \frac{1}{\lambda_{l}}\int_{\lambda_{l-1}t_k }^{\lambda_lt_k }\tilde{W}_l(t)dt,
	\end{align*}
	where $\tilde{W}_l(\cdot)$ and $\tilde{W}_{l-1}(\cdot)$ are the time-change version of $W_l(\cdot)$ and $\bar{W}_{l-1}(\cdot)$, respectively, such that their Poisson arrival processes are both of rate 1.
	For the first term,  we have
	\begin{align*}
		&\EE\left[\left|\tilde{W}_l(t)-\tilde{W}_{l-1}(T_{k(l-1)}+t)\right|\right]\\
		\leq~& \EE\left[\left|\frac{\mu_l}{\lambda_l}-\frac{\mu_{l-1}}{\lambda_{l-1}}\right|\max(\tilde{X}_{l}(t) ,\tilde{X}_{l-1}(T_{k(l-1)}+t))+|W_l(0)-\bar{W}_{l-1}(T_{k(l-1)})|\right]\\
		\stackrel{(a)}{\leq} ~&\EE\left[\left|\frac{\mu_l}{\lambda_l}-\frac{\mu_{l-1}}{\lambda_{l-1}}\right|\tilde{X}^D_l(t)\right]+\EE\left[|W_{l-1}(T_{k(l-1)})-\bar{W}_{l-1}(T_{k(l-1)})|\right]\\
		\stackrel{(b)}{\leq}~& \EE\left[\left|\frac{\mu_l}{\lambda_l}-\frac{\mu_{l-1}}{\lambda_{l-1}}\right|\tilde{X}^D_l(t)\right] +O(k^{-1})\\
		\leq~& \EE\left[\left|\frac{\mu_l}{\lambda_l}-\frac{\mu_{l-1}}{\lambda_{l-1}}\right|^2\right]^{1/2}\EE\left[\tilde{X}^D_l(t)^2\right]^{1/2} +O(k^{-1})\\
		\stackrel{(c)}{=}~& O(\max(\eta_k\sqrt{\mathcal{V}_k},\delta_k))+O(k^{-1})= O(\max(\eta_k\sqrt{\mathcal{V}_k},\delta_k)),
	\end{align*}
	where $\tilde{X}^D_l(\cdot)$ is the dominant observed busy period defined in the proof of Lemma \ref{lmm: uniform bound}.
	Here inequality $(a)$ follows from the definition of $\tilde{X}^D_l(\cdot)$, inequality $(b)$ from Lemma \ref{coro: regret of nonstationarity after warm-up period}  and equality $(c)$ from Lemma \ref{lmm: uniform bound} and the fact that 
	$$\|\bx_l-\bx_{l-1}\|=\begin{cases}
		\delta_k&\text{ for }l = 2k\\
		\eta_k\|\bH_{k-1}\|&\text{ for }l = 2k-1.
	\end{cases}$$
	For the second term,
	\begin{align*}
		&\EE\left[\left|\frac{1}{\lambda_l} -\frac{1}{\lambda_{l-1}} \right|\int_{0}^{\lambda_{l-1} t_k}\tilde{W}_{l-1}(t)dt \right]=\EE\left[\left|1 -\frac{\lambda_{l-1}}{\lambda_{l}} \right|\int_{0}^{t_k}W_{l-1}(t)dt \right]\\
		\leq~& \frac{1}{\underline{\lambda}}\EE\left[(\lambda_l-\lambda_{l-1})^2\right]^{1/2}\EE\left[\left(\int_0^{t_k}W_{l-1}(t)dt\right)^2\right]^{1/2}= O(\max(\eta_k\sqrt{\mathcal{V}_k} ,\delta_k)t_k).
	\end{align*}
	Following a similar argument, we have that 
	$$\EE\left[\frac{1}{\lambda_{l}}\int_{\lambda_{l-1}t_k }^{\lambda_lt_k }\tilde{W}_l(t)dt\right]=\EE\left[\int^{t_k}_{\frac{\lambda_{l-1}}{\lambda_l}t_k}W_l(t)dt\right]=O (\max(\eta_k\sqrt{\mathcal{V}_k} ,\delta_k)t_k).$$
	In summary, we can conclude that there exists a constant $C_0>0$ such that
	\small{$$\EE\left[\int_0^{t_k}(W_l(t)-w_l)dt\right]\leq t_k\EE\left[|w_l-w_{l-1}|\right]+\EE\left[\left|\int_0^{t_k}(W_l(t)-\bar{W}_{l-1}(T_{k(l-1)}+t))dt\right|\right]\leq C_0\max(\eta_k\sqrt{\mathcal{V}_k} ,\delta_k)t_k.$$}
	As a consequence,
	$$\EE\left[\int_{0}^{t_k} (W_l(t)-w_l) dt\right]\leq C_0\max(\eta_k\sqrt{\mathcal{V}_k},\delta_k)t_k=O\left(\max(\eta_k\sqrt{\mathcal{V}_k},\delta_k)\log(k)\right).$$
	\hfill $\Box{}$
\end{proof}
\vspace{2ex}
\begin{proof}{Proof of Lemma \ref{lmm: exploration cost}}
	By Taylor's expansion and the mean value theorem, 
	\begin{align*}%\label{eq: R_3k}
		R_{3k}=\EE[T_k\left(f(\bx_{2k-1})+f(\bx_{2k})-2f(\bar{\bx}_k)\right)]= \EE[T_k(f''(\bx')+f''(\bx''))\delta_k^2]\leq K_4T_k\delta_k^2,
	\end{align*}
	where $\bx', \bx''\in\mathcal{B}$ and the last inequality follows from Lemma \ref{lmm: convexity}. \hfill$\Box$
\end{proof}

\subsection{Proof of Theorem \ref{thm: heavy traffic}}\label{sec: proofs in heavy traffic analysis}
The proof of Theorem \ref{thm: heavy traffic} follows a structure similar to that of the proof of Theorem \ref{thm: upper bound}. We first need to build bounds on (i) moments; (ii) transient bias of the queueing data; (iii) variance of the queueing data; (iv) and FD approximation error of the gradient in terms of the parameter $h$ which corresponds to Lemmas \ref{lmm: moment bounds heavy traffic} to \ref{lmm: FD error heavy traffic}. Based on the results, we could bound the bias and variance of our gradient estimator in Lemma \ref{lmm: bias and variance heavy traffic} and the order of strong-convexity coefficient in Lemma \ref{lmm: convexity heavy traffic}. Then, following the regret decomposition in the main paper, we bound the regret of suboptimality, nonstationary and finite difference in Lemmas \ref{lmm: suboptimal regret heavy traffic} to \ref{lmm: FD regret heavy traffic}, which complete the proof of Theorem \ref{thm: heavy traffic}.

For M/M/1 queue with unit service rate, the mean stationary workload is equal to mean stationary queueing length (including the customer in service). So, one could estimate the objective function using the observed queue length data, and hence, entirely eliminate the bias of delayed observation. In the following analysis, we use $Q^h_l(t)$ and $p_l^h$  to denote the observed queueing length and control price, respectively, in cycle $l$ when applying LiQUAR to the $h$-th system. 

In addition, when applying LiQUAR to the $h$-th system, we denote the gradient estimator in iteration $k$ as 
$$H_k^h=\frac{1}{2\delta_k^h}\left[-p^h_{2k-1}\frac{N^h_{2k-1}}{T_k}+p^h_{2k}\frac{N^h_{2k}}{T^h_k}+h\int_{\alpha T^h_k}^{T^h_k}Q^h_{2k-1}(t)-Q^h_{2k}(t)dt\right]$$ and the corresponding finite difference 
$$\frac{f_h(p^h_{2k-1})-f_h(p^h_{2k+1})}{2\delta_k^h}\equiv Df_h(\bar{p}^h_k),$$
where
$$p^h_{2k-1}=\bar{p}^h_k+\delta_k^h, \quad p^h_{2k}=\bar{p}^h_k-\delta_k^h.$$
Following the main paper, we define the bias and variance of the gradient estimator as 
$$B_k^h\equiv \EE[(\EE[H_k^h-f'(\bar{p}_k^h)|\mathcal{F}_k])^2],\quad \mathcal{V}_k^h \equiv \EE[(H_k^h)^2].$$
For the simplicity of notation, we will denote all positive constants that are independent of $h$ and $T_0$ by $C$ in the following analysis.
\begin{lemma}[Moment Bounds]\label{lmm: moment bounds heavy traffic} Under any control sequence $p_l^h$, 
	$$\EE\left[(Q_l^h(t))^m\right]\leq Ch^{-m/2}, \text{ for all }l\geq 1\text{ and }t\in[0,T_k].$$
\end{lemma}
\begin{proof}{Proof of Lemma \ref{lmm: moment bounds heavy traffic}}
	Let $\tilde{Q}_h(\cdot)$ be the stationary queue length process of an $M/M/1$ queue with service rate $1$ and arrival rate $\lambda(p^*+c_1\sqrt{h})$. Then, for arbitrary control sequence $p_l^h$, we have $$Q_l^h(t)\leq_{st}\tilde{Q}_h(t),$$
	for all $t\geq 0$. Therefore, it is sufficient to show that $$\EE[\tilde{Q}_h(t)^m]\leq Ch^{-m/2}$$
	for some $C>0$ and $1\leq m \leq 4$. By Taylor expansion, $\lambda(p^*+c_1 \sqrt{h})=1+\lambda'(p^*+\theta c_1\sqrt{h})c_1\sqrt{h}$ with some $\theta\in(0,1)$, so the corresponding traffic intensity satisfies 
	$$1-\rho=-\lambda'(p^*+\theta c_1\sqrt{h})c_1\sqrt{h}\leq c_1\cdot C_0 \sqrt{h},$$
	with $C_0=-\arg\min_{p\in \mathcal{B}_1} \lambda'(p)$. Then, by the stationary distribution of $M/M/1$ queue, the moment bounds are valid. $\hfill\Box{}$
\end{proof}

\begin{lemma}[Transient Bias Bound]\label{lmm: transient bias heavy traffic}Suppose $\bar{Q}^h_l(\cdot)$ is a stationary queue length process synchronously coupled with $Q^h_l(\cdot)$. Then, conditional on their initial values, 
	$$\EE[|Q_l^h(t)-\bar{Q}_l^h(t)||Q_l^h(0),\bar{Q}_l^h(0)]\leq |Q_l^h(0)^2-\bar{Q}_l^h(0)^2|\cdot \frac{2{C}t^{-3/2}}{h}\exp(-ht/2{C}).$$
	%As a consequence,
	%$$\EE[(\EE[H_k^h-Df_h(p_k))^2]^{1/2}\leq Ch^{-1/2}k^{-2/3}.$$
\end{lemma}

\begin{proof}{Proof of Lemma \ref{lmm: transient bias heavy traffic}}
	Consider an M/M/1 queue with traffic intensity $\rho$ and $i$ customers in the system at time 0. Let $\tau$ be the first hitting time when the system gets empty. Following theorem 3.1 in \cite{abate1988transient},
	\begin{align*}
		\PP((1-\rho)^{2}\tau>t)=\int_t^\infty f(s;i,0)ds,
	\end{align*}
	with
	\begin{align*}
		f(t;i,0)=(i/t)\rho^{1/2}\exp(-2t/(1+\sqrt{\rho})^2)\exp(-4\rho^{1/2}t/(1-\rho)^2)I_i(4\rho^{1/2}t/(1-\rho)^2).
	\end{align*}
	Here $I_i(x)$ is the modified Bessel function of the first kind such that
	$I_i(x)\leq I_0(x)$ for any integer $i\geq 0$. By \cite{olivares2018simple}, for all $x>0$,
	$$I_0(x)\leq   1.006\cdot \frac{e^x+e^{-x}}{2(1+x^2/4)^{1/4}}\frac{1+0.24273x^2}{1+0.43023x^2}\leq 1.006\cdot\frac{e^x}{(1+x^2/4)^{1/4}}\leq 1.006\cdot e^x\cdot (1\wedge \sqrt{2/x}).$$
	We bound $f(t;i,0)$ by
	$$f(t;i,0)\leq 1.006\cdot (i/t)\exp(-t/2)\cdot\left(1\wedge (1-\rho)\sqrt{1/t}\right),$$
	if $\rho>1/4$. Therefore, for $t\geq 1$,
	\begin{align*}
		\PP(\tau>t)&=\PP\left((1-\rho)^2\tau>(1-\rho)^2t\right)=\int_{(1-\rho)^{2}t}^\infty f(s;i,0)ds\\
		&\leq \int_{(1-\rho)^{2}t}^\infty 1.006\cdot\frac{i}{s}\exp(-s/2)(1-\rho)\sqrt{1/s}~ds\\
		&\leq 2.012 (1-\rho)is^{-3/2}\exp(-s/2)\huge{|}_{s=(1-\rho)^2t}\\
		&=\frac{2.012i}{(1-\rho)^2}t^{-3/2}\exp(-(1-\rho)^2t/2)
	\end{align*}
	{The last inequality comes from integral by part.} Suppose we synchronously couple an M/M/1 queue length process $Q(t)$ with a stationary one $\bar{Q}(t)$ and denote by $\bar{\tau}$ the first hitting time to 0 of $\bar{Q}(t)$. Then, we have
	\begin{align*}
		\EE[|Q(t)-\bar{Q}(t)| | Q(0)=i] &\leq \EE[|i-\bar{Q}(0)|1(\tau\vee\bar{\tau}>t)]\\
		&\leq \EE\left[\frac{2.012|i-\bar{Q}(0)|(i+\bar{Q}(0))}{(1-\rho)^2}t^{-3/2}\exp(-(1-\rho)^2t/2)\right]\\
		&\leq \EE[|i^2 -\bar{Q}(0)^2|]\cdot \frac{2.012}{(1-\rho)^2}t^{-3/2}\exp(-(1-\rho)^2t/2).
	\end{align*} 
	Note that for $p\in\mathcal{B}_h$, $1-\rho=O(\sqrt{h})$ . Then, setting $Q(t),\bar{Q}(t)$ being the $Q_l^h(t),\bar{Q}_l^h(t)$ closes the proof. $\hfill\Box{}$
\end{proof}

\begin{lemma}[Variance Bound]\label{lmm: variance bound heavy traffic} For all $h$ and $l$, the stationary queue satisfies
	$$Var\left[\int_0^T \bar{Q}_l^h(t)ds\right]\leq \frac{CT}{h^2}.$$
\end{lemma}
\begin{proof}{Proof of Lemma \ref{lmm: variance bound heavy traffic}}
	Let $c_q(t)= corr(\bar{Q}_l^h(0),~\bar{Q}_l^h(2t/(1-\rho)^2)~)$, with $\rho = 1-\lambda(p^h_l)$ and thus $1-\rho\geq C\sqrt{h}$. According to corollary 5 of \cite{abate1988correlation}, 
	$$\int_0^\infty c_q(t)dt = \frac{1+\rho}{2}\leq 1.$$
	Consequently, we have
	$$\int_{0}^{\infty} Cov(\bar{Q}_l^h(0),\bar{Q}_l^h(2t/(1-\rho)^2))dt\leq \EE[\bar{Q}_l^h(0)]^2=\frac{\rho(1+\rho)}{(1-\rho)^2}\leq \frac{C}{h}.$$
	By changing of variables, we would see
	$$\int_{0}^{\infty} Cov(\bar{Q}_l^h(0),\bar{Q}_l^h(t))dt\leq  \frac{C}{h^2}.$$
	Now, we have 
	\begin{align*}
		Var\left[\int_0^T \bar{Q}_l^h(t)ds\right]&=\int_{0}^{T}\int_{0}^{T}Cov(\bar{Q}_l^h(t),\bar{Q}_l^h(s))dtds\\
		&\leq 2\int_{0}^{T}\int_{0}^{\infty}Cov(\bar{Q}_l^h(t),\bar{Q}_l^h(t+s))dsdt\leq \frac{CT}{h^2}.
	\end{align*} $\hfill\Box{}$
\end{proof}
\begin{lemma}[FD Approximation Error Bound]\label{lmm: FD error heavy traffic}
	$$|Df_h(p^h_k)-f'_h(p^h_k)|\leq Ck^{-2/3}.$$
	%$$Df_h(p)(p-p_h)\geq h^{-1/2}C(p-p_h)^2.$$
\end{lemma}
\begin{proof}{Proof of Lemma \ref{lmm: FD error heavy traffic} }
	For fixed $h$, $p\in \mathcal{B}_h$ and $\delta>0$,
	%\begin{align*}
	%	&\frac{f_h(p+\delta)-f_h(p-\delta)}{2\delta}=\frac{(p+\delta)\lambda(p+\delta)-(p-\delta)\lambda(p-\delta) -\left(\frac{h}{1-\lambda(p+\delta)}-\frac{h}{1-\lambda(p-\delta)}\right)}{2\delta}\\
	%	=~&\frac{p(\lambda(p+\delta)-\lambda(p-\delta)}{2\delta} +\frac{\lambda(p+\delta)+\lambda(p-\delta)}{2}+\frac{h(\lambda(p+\delta)-\lambda(p-\delta))}{2\delta(1-\lambda(p+\delta))(1-\lambda(p-\delta))}\\
	%	=~& p\dot{\lambda}(p_1) +\lambda(p_2) + \frac{h\dot{\lambda}(p_1)}{(1-\lambda(p+\delta))(1-\lambda(p-\delta))},
	%\end{align*}
	%where $p_1,p_2\in [p-\delta,p+\delta]$.
	\begin{align*}
		f_h(p+\delta)=f_h(p) +\delta f_h'(p)+\frac{\delta^2}{2}f_h''(p)+\frac{\delta^3}{6}f_h'''(p_1)\\
		f_h(p-\delta)=f_h(p) -\delta f_h'(p)+\frac{\delta^2}{2}f_h''(p)-\frac{\delta^3}{6}f_h'''(p_2)
	\end{align*}
	Therefore,
	\begin{align*}
		\frac{f_h(p+\delta)-f_h(p-\delta)}{2\delta}=f_h'(p) +\frac{\delta^2f_h'''(p_3)}{6}.
	\end{align*}
	Note that
	$$f_h'''(p)=3\lambda''(p)+p\lambda'''(p)-\frac{6h\lambda'(p)^3}{(1-\lambda(p))^4}-\frac{6h\lambda''(p)\lambda(p)}{(1-\lambda(p))^3}-\frac{h\lambda'''(p)}{(1-\lambda(p))^2}.$$
	As $1-\lambda(p)=O(\sqrt{h})$, we can conclude that
	$$f_h'''(p)=O(h^{-1}).$$
	As $\delta=O(\sqrt{h}k^{-1/3})$, we conclude that the FD approximation error is of order $O(k^{-2/3})$. $\hfill\Box{}$
\end{proof}
\begin{lemma}[Bounds on Gradient Estimator Bias and Variance]\label{lmm: bias and variance heavy traffic}
	For all $h$ and $k$,
	$$B_k^h\leq C\cdot k^{-2/3},\quad \mathcal{V}_k^h\leq C.$$
\end{lemma}
\begin{proof}{Proof of Lemma \ref{lmm: bias and variance heavy traffic}}
	We first prove the bias term and then we prove the variance term. 
	\paragraph{Bias term}
	By definition, the bias is defined by 
	$$(B_k^h)^2= \EE[(\EE[H_k^h-f'(\bar{p}_k^h)|\mathcal{F}_k])^2]\leq2\EE[\EE[f_h'(p_k^h)-Df_h(p_k^h)|\mathcal{F}_k]^2]+ 2\EE[\EE[H_k^h-Df_h(p_k^h)|\mathcal{F}_k]^2].$$
	By Lemma \ref{lmm: FD error heavy traffic}, we have following bound for the first term. 
	$$\EE[\EE[f_h'(p_k^h)-Df_h(p_k^h)|\mathcal{F}_k]^2]\leq C k^{-4/3}.$$
	We next bound the second term. By Lemma \ref{lmm: transient bias heavy traffic}, we have 
	\begin{align*}
		\EE[|Q_l^h(t)-\bar{Q}_l^h(t)||Q_l^h(0),\bar{Q}_l^h(0)]\leq |Q_l^h(0)^2-\bar{Q}_l^h(0)^2|\cdot \frac{2{C}t^{-3/2}}{h}\exp(-ht/2{C}).
	\end{align*}
	Consequently, we have 
	\begin{align*}
		\EE[\hat{f}_h(p_l^h)-f_h(p_l^h)\big|\mathcal{G}_l]&=\frac{h}{(1-\alpha)T_k^h}\EE\left[\int_{\alpha T_k^h}^{T^h_k}Q_l^h(t)-\bar{Q}_l^h(t)dt\Big|\mathcal{G}_l\right]\\
		&\leq  \frac{|Q_l^h(0)^2-\bar{Q}_l^h(0)^2|}{(1-\alpha )T^h_k}\cdot \int_{\alpha T^h_k}^{T^h_k}2 {C}t^{-3/2}\exp(-ht/2{C})dt\\
		&\leq  \frac{{C}|Q_l^h(0)^2-\bar{Q}_l^h(0)^2|}{\alpha^{3/2} (T_k^h)^{3/2}}\exp(-\alpha h T^h_k/2{C}),
	\end{align*}
	where the last inequality holds due to the monotonicity of $t^{-3/2}\exp(-ht/2{C})$.
	Therefore, by our choice of $T_k^h,\delta_k^h$, we have 
	\begin{align*}
		&\EE[H_k^h-Df_h(p_k^h)|\mathcal{F}_k]=\frac{{C}\EE[Q_l^h(0)^2-\bar{Q}_l^h(0)^2|\mathcal{F}_k]}{\delta_k (T_k^h)^{3/2}}\exp (-\alpha h T^h_k/2{C})\\
		\leq& C\frac{h^{-1}}{\sqrt{h} k^{-1/3} h^{-3/2}k^{1/2}} \exp(-\alpha k^{1/3}/2{C})\leq C\cdot k^{-2/3},
	\end{align*}
	for sufficient large $k$. This closes the proof of Bias.
	\paragraph{Variance Term} For the variance term, we have 
	$$\EE[H_k^2]\leq 3(\delta_k^h)^{-2}\sum_{l=2k-1}^{2k}\EE[\hat{f}_h(p_l^h)-f_h(p_l^h)^2]+3(\delta_k^h)^{-2} \EE[f_h(p_{2k}^h)-f_h(p_{2k-1}^h)^2].$$
	For the second term, we calculate that for $p\in\mathcal{B}_h$,
	$$f'_h(p)=-p\lambda'(p)-\lambda(p)+h\frac{\lambda'(p)}{(1-\rho(p))^2}=O(1).$$
	Consequently, we have
	$$(\delta_k^h)^{-2} \EE[f_h(p_{2k}^h)-f_h(p_{2k-1}^h)^2]\leq \max_{p\in\mathcal{B}_h}\|f'_h(p)\|=O(1).$$
	For the first term, we have 
	$$\EE[(\hat{f}_h(p_l^h)-f_h(p_l^h))^2]\leq 2\EE\left[\left(p_l^h\frac{N_l}{T_k}-p_l^h\lambda(p_l^h)\right)^2\right]+2\frac{h^2}{((1-\alpha)T_k^h)^2}\EE\left[\left(\int_{\alpha T_k^h}^{T_k^h}Q_l^h(t)-\EE[\bar{Q}_l^h(t)]\right)^2\right].$$
	Let's denote $\bar{Q}_l^h(t)$ as a stationary version of queueing process synchronously coupled with $Q^h_l(t)$, and define $\tau,\bar{\tau}$ as the first hitting time of them to the empty states. Note that 
	\begin{align*}
		&\EE\left[\left(\int_{\alpha T_k^h}^{T_k^h}Q_l^h(t)-\EE[\bar{Q}_l^h(t)]dt\right)^2\right]\\
		\leq& \EE\left[\left(\int_{\alpha T_k^h}^{T_k^h}\bar{Q}_l^h(t)-\EE[\bar{Q}_l^h(t)]dt\right)^2\right]+\EE\left[\left(\int_{\alpha T_k^h}^{T_k^h}Q_l^h(t)-\EE[\bar{Q}_l^h(t)]dt\right)^2\textbf{1}(\tau\vee \bar{\tau}>\alpha T_k^h)\right]\\
		\stackrel{(a)}{\leq}&\frac{C(1-\alpha)T_k^h}{h^2}+(1-\alpha )T_k^h \EE\left[\int_{\alpha T_k^h}^{T_k^h}(Q_l^h(t)-\EE[\bar{Q}_l^h](t))^2dt\textbf{1}(\tau\vee \bar{\tau}>\alpha T_k^h)\right]\\
		\leq &\frac{C(1-\alpha)T_k^h}{h^2}+CT_k^h\cdot \frac{T_k^h}{h}\cdot \PP(\tau\vee \bar{\tau}>\alpha T_k^h)^{1/2}\\
		\leq &\frac{C(1-\alpha)T_k^h}{h^2}+\frac{C}{h^2}T_k^h\cdot hT_k^h\cdot \frac{\sqrt{h}\EE[Q_l^h(0)+\bar{Q}_l^h(0)]}{(hT_k^{h})^{3/2}}e^{-hT_k^h/2{C}}\\
		\leq &\frac{C}{h^2}T_k^h.
	\end{align*}
	Here, the inequality (a) comes from Lemma \ref{lmm: variance bound heavy traffic} and the Cauchy-Schwartz inequality, and the last inequality comes from the fact that $hT_k^h\rightarrow\infty$ and $\sqrt{h}\EE[Q_l^h(0)+\bar{Q}_l^h(0)]=O(1)$.
	Consequently, we have 
	$$\EE[(\hat{f}_h(p_l^h)-f_h(p_l^h))^2]\leq\frac{C}{T_k},$$
	for some $C$ large enough.
	Therefore, we have 
	$$\EE[H_k^2]\leq  \max\left(\frac{C}{T_k^h\delta_k^2},C\right)=C.$$ $\hfill\Box{}$
\end{proof}
\begin{lemma}[Convexity]\label{lmm: convexity heavy traffic}
	There exists a constant $K_0>0$ independent of $h$ such that, for all $p\in\mathcal{B}_h$, 
	$$f''_h(p)>h^{-1/2}K_0.$$
\end{lemma}
\begin{proof}{Proof of Lemma \ref{lmm: convexity heavy traffic}}
	Note that for all $p\in\mathcal{B}_h$, the traffic intensity $1-\rho(p)=O(1/\sqrt{h})$. Then, by direct calculation and Polleczk-Khinchine formula, we have 
	$$f''_h(p)=(-p\lambda(p))''+\frac{h}{(1-\rho(p))^3}\left(2(\lambda'(p)^2+(1-\rho(p))\lambda''(p))\right)>h^{-1/2}K_0,$$
	with $K_0=2\min_{p\in \mathcal{B}_1} 2|\lambda'(p)|^2$. $\hfill\Box{}$
\end{proof}

Given Lemmas \ref{lmm: bias and variance heavy traffic} and \ref{lmm: convexity heavy traffic}, we are ready to provide an upper bound on the $L_2$ distance $\EE[(\bar{p}_k^h-p^*_h)^2]$ following the analysis of main paper.
\begin{lemma}[Suboptimal Regret]\label{lmm: suboptimal regret heavy traffic}
	The suboptimal regret could be bounded by 
	$$R_1^h(L)\leq C\cdot \frac{L^{2/3}}{\sqrt{h}}.$$
\end{lemma}
\begin{proof}{Proof of Lemma \ref{lmm: suboptimal regret heavy traffic}}
	
	For all $h>0$ and $k\geq 1$, we denote 
	$$b^h_k\equiv h^{-1}(\bar{p}_k^h-p_h)^2.$$
	For a given $h$ small enough, we omit the superscript $h$ for the simplicity of notation and obtain
	\begin{align*}
		hb_{k+1}& = \EE[(\bar{p}_{k+1}-p^*)^2]\leq \EE[(\bar{p}_k-p^*-\eta_kH_k)]\\
		&= \EE[(\bar{p}_k-p^*)^2-2\eta_kf'(\bar{p}_k)\cdot(\bar{p}_k-p^*)]-2\eta_k\EE[(H_k-f'(\bar{p}_k))\cdot(\bar{p}_k-p^*)]+2\eta_k^2\EE[H_k^2]\\
		&\leq (1-2\eta_kh^{-1/2}K_0)\EE[(\bar{p}_k-p^*)^2]+\sqrt{h}\eta_kB_k (1+b_k)+2\eta_k^2V_k\\
		&= (1-2c_\eta K_0k^{-1})hb_k + hc_\eta k^{-1} B_k + hc_\eta k^{-1} B_kb_k + 2hc_\eta^2V_k\\
		&\leq h\cdot \left[(1-2c_\eta K_0k^{-1})b_k + Ck^{-5/3} + Ck^{-5/3}b_k+Ck^{-2} \right].
	\end{align*}
	Following the proof of theorem 2 in the main paper, we can prove by induction that, there exists a constant $C>0$ independent of $h$ such that
	$b_{k}\leq Ck^{-2/3}$, and therefore, we can conclude
	$$\EE[(\bar{p}_k^h-p^*_h)^2]\leq C\cdot hk^{-2/3}.$$
	As a result, we have 
	\begin{align*}
		R_1^h(L)=&\sum_{k=1}^{L}\EE\left[(f(\bar{p}_k^h)-f(p^*_h))T_k^h\right]\\
		\leq&\sum_{k=1}^{L}\EE\left[\nabla^2 f(p^*+c_1 \sqrt{h})(\bar{p}_k^h-p_h^*)^2T_k^h\right]\\
		\leq &\sum_{k=1}^{L} C\sqrt{h}k^{-2/3} T_k^h\leq C\cdot  \frac{L^{2/3}}{\sqrt{h}}
	\end{align*} $\hfill\Box{}$
\end{proof}
\begin{lemma}[Non-stationary Regret]\label{lmm: non-stationary regret heavy traffic}
	The non-stationary regret could be bounded by 
	$$R_2^h(L)\leq C\cdot \frac{L^{2/3}\log(L)}{\sqrt{h}}.$$
\end{lemma}
\begin{proof}{Proof of Lemma \ref{lmm: non-stationary regret heavy traffic}}
	Following the decomposition of non-stationary regret in the main paper, we have 
	\begin{align*}
		R_{2k}^h&=\sum_{l=2k-1}^{2k} h\EE\left[\int_{0}^{T_k^h}Q_l^h(t)-\bar{Q}^h_l(t)dt\right]\\
		&=\sum_{l=2k-1}^{2k} h\EE\left[\int_{0}^{t_k^h}Q_l^h(t)-\bar{Q}^h_l(t)dt\right]+h\EE\left[\int_{t_k^h}^{T_k^h}Q_l^h(t)-\bar{Q}^h_l(t)dt\right],
	\end{align*}
	with $t_k^h=\frac{2\log k}{h}$. In this way, we following the similar analysis in our main paper. For the second term, by Lemma \ref{lmm: transient bias heavy traffic}, we have 
	\begin{equation}\label{eq: heavy-R2k-late}
		h\int_{t_k^h}^{T_k^h}\EE[|Q_l^h(t)-\bar{Q}_l^h(t)|]dt\leq \int_{t_k^h}^{\infty}\frac{{C}\EE[|Q^h_l(0)^2-\bar{Q}_l^h(0)|^2]}{ht_k^{3/2}}\exp(-ht/2{C})dht\stackrel{(b)}{\leq}\frac{C}{h^2t_k^{3/2}}\cdot k^{-1} \leq \frac{C}{\sqrt{h}k}.
	\end{equation}
	The inequality (b) comes from the fact that $\EE[Q_l^h(0)^2],\EE[\bar{Q}_l^h(0)^2]=O(h^{-1})$.
	For the first term, we decompose $\EE[Q_l^h(t)-Q_l^h(t)]$ into $\EE[\bar{Q}^h_{l-1}(t)-\bar{Q}_l^h(t)]$ and $\EE[Q_l^h(t)-\bar{Q}_{l-1}(t)]$ as we did in the main paper. By Polleczk-Khinchine formula, we have 
	$$\EE[\bar{Q}^h_{l-1}(t)-\bar{Q}_l^h(t)]= \EE\left[\frac{\lambda(p_l^h)-\lambda(p_{l-1}^h)}{(1-\lambda(p_l^h))(1-\lambda'(p^h_{l-1}))}\right]\leq \frac{C}{h}\EE[|p_l^h-p^h_{l-1}|\leq \frac{C}{h}\EE[|p_l^h-p^h_{l-1}|^2]^{1/2}.$$
	Next, following the same argument in Lemma 7 in the main paper, we define $\tilde{Q}_l^h(\cdot)$ and $\tilde{X}_l^h(\cdot)$ as the queue length and busy period process with arrival rate $1$ and service rate $1/\lambda(p_l^h)$. Then, by the same analysis in Lemma 7 in the main paper,
	\begin{align*}
		&\int_{0}^{t_k^h}\EE[Q^h_l(t)-\bar{Q}_{l-1}^h(t)]dt\\
		\leq & \frac{1}{\lambda_l}\int_{0}^{\lambda_{l-1}t_k^h}\EE[|\tilde{Q}_l^h(t)-\tilde{Q}_{l-1}(T_{k(l-1)}^h+t)dt|]+\EE\left[\left|\frac{1}{\lambda_l}-\frac{1}{\lambda_{l-1}}\right|\int_{0}^{\lambda_{l-1}t_k^h}\tilde{Q}_{l-1}(t)dt\right]+\EE\left[\frac{1}{\lambda_l}\int_{\lambda_{l-1}t_k^h}^{\lambda_l t_k^h}\tilde{Q}_l(t)dt\right]\\
		\leq & \EE\left[\left|\frac{1}{\lambda(p_l^h)}-\frac{1}{\lambda(p_{l-1}^h)}\right|^2\right]^{1/2}\EE[\tilde{X}_l(t)^2]^{1/2} t_k^h+C \EE\left[\left|\frac{1}{\lambda(p_l^h)}-\frac{1}{\lambda(p_{l-1}^h)}\right|^2\right]^{1/2} \EE[\tilde{Q}_{l-1}^h(t)^2]^{1/2}t_k^h\\
		\leq & C \EE[|p_l^h-p_{l-1}^h|^2]^{1/2} \frac{t_k^h}{h}
	\end{align*}
	Therefore, by Lemma \ref{lmm: variance bound heavy traffic}, we have 
	\begin{equation}\label{eq: heavy-R2k-early}
		h\EE\left[\int_{0}^{t_k^h}Q_l^h(t)-\bar{Q}_l^h(t)dt\right]\leq C\cdot\EE[|p_l^h-p^h_{l-1}|^2]^{1/2}t_k^h\leq  Ct_k^h \cdot \max(\eta_k\sqrt{\mathcal{V}_k},\delta_k)=C\frac{\log k}{\sqrt{h}k^{1/3}}.
	\end{equation}
	Combining equations \eqref{eq: heavy-R2k-early} and \eqref{eq: heavy-R2k-late}, we have $R_{2k}^h\leq C\frac{\log k}{\sqrt{h}k^{1/3}}$. Therefore, we have
	$$R_{2}^h(L)\leq C \frac{L^{2/3}\log L}{\sqrt{h}}.$$ $\hfill\Box{}$
\end{proof}
\begin{lemma}[Finite-Difference Regret]\label{lmm: FD regret heavy traffic}
	The finite-difference regret could be bounded by 
	$$R_3^h(L)\leq C\cdot \frac{L^{2/3}}{\sqrt{h}}.$$
\end{lemma}
\begin{proof}{Proof of Lemma \ref{lmm: FD regret heavy traffic}}
	By calculation, we have
	\begin{align*}
		R_3^h(L)=&\sum_{k=1}^{L}\EE\left[(f(\bar{p}_{2k-1})+f(\bar{p}_{2k})-2f(p^*_h))T_k^h\right]
		\leq C \sum_{k=1}^{L}\frac{1}{\sqrt{h}}(\delta_k^h)^2 T_k^h
		\leq C\cdot \frac{L^{2/3}}{\sqrt{h}}.
	\end{align*} $\hfill\Box{}$
\end{proof}

Summing up all three regrets, the total regret in the first $L$ cycle is 
$$R^h(L)=R^h_1(L)+R^h_2(L)+R^h_3(L)\leq C\frac{L^{2/3}\log L}{\sqrt{h}}.$$
Note that the total time that used is $T^h=\frac{T_0}{h}=\frac{L^{4/3}}{h}$, and therefore,
$$R^h(T_0/h)\leq C \sqrt{h}^{-1}\sqrt{T_0\log T_0}.$$ $\hfill\Box{}$
\subsection{Proof of Proposition \ref{prop: PTO heavy-traffic}}
We neglect the superscribe $h$ in the following analysis to ease the burden of notation. The proof of Proposition \ref{prop: PTO heavy-traffic} basically follows the proof of proposition 2 in \cite{besbes2009dynamic}. Let $\Delta_0=\max_{p\in\mathcal{B}_h} f^h(p) - f^h(p^*)=O(\sqrt{h})$, and denote $p_G^*$ as the optimal points in the testing pricing grid. The regret can be decomposed according to three sources of cost: exploration cost, stochastic error, and discrete grid cost.
% Specifically, we have 
%	\begin{enumerate}
	%		\item Price experimentation cost $O(\Delta_0^h t_0^h)=O(\sqrt{h} t_0^h)$
	%		\item Stochastic Fluctuation Error $|f(p^*)-f(p_G^*)|=O(1/\sqrt{t_0^h/\kappa ^h})=O(\sqrt{\kappa^h/t_0^h})$
	%		\item Deterministic Grid Error $|p^*-p_G^*|=O(|\mathcal{B}_h|/\kappa^h)$,
	%	\end{enumerate}

\begin{align*}
	R^h(T_0)&\leq \Delta_0 t_0+ (T-t_0) \EE[f(\hat{p}^*)-f(p^*)]\\
	&\leq \underbrace{\Delta_0 t_0}_{\text{Exploration Cost}}+\EE[T\cdot [\underbrace{f(\hat{p}^*)-f(p^*_G)}_{\text{Stochastic Error}}+\underbrace{f(p^*_G)-f(p^*)}_{\text{Discrete Grid Cost}}]]
\end{align*}
We treat the first two terms following in the same way as in \cite{besbes2009dynamic}. For the third term, we apply second order Taylor expansion (rather than the first order in \cite{besbes2009dynamic}) as $\nabla f(p^*)=0$ in our problem. Therefore, using the fact that the grid length is at most $|p_G^*-p^*|\leq |\mathcal{B}_h|/\kappa=O(\sqrt{h}/\kappa)$, we have 
%	Note that for stochastic error, the best $\hat{p}^*$ are chosen from the price grids after time $t_0/\kappa$, and following the analysis in \cite{besbes2009dynamic}, it has the regret of the order $O(\sqrt{\log T/ (t_0/\kappa)})$. For the discrete grid cost, note that the grid length is at most $|p_G^*-p^*|\leq |\mathcal{B}_h|/\kappa=O(\sqrt{h}/\kappa)$.
\begin{align*}
	R^h(T_0)&\leq \Delta_0 t_0 + CT\cdot \sqrt{\frac{\kappa \log T}{t_0}} + CT\cdot \frac{\nabla^2 f(\xi) }{2} |p_G^*-p^*|^2\\
	&\leq C\sqrt{h} t_0 + CT\cdot \sqrt{\frac{\kappa \log T}{t_0}}  +C\frac{T}{\sqrt{h}}\cdot \left( \frac{\sqrt{h}}{\kappa}\right)^2.
\end{align*}
By optimizing the regret order, we choose $$t_0=O\left(\frac{T_0^{5/7}\log (T)^{2/7}}{h}\right),\quad \kappa=O\left(\frac{T_0^{1/7}}{\log T}\right),$$ such that 
$$R^h(T_0)=O\left(\frac{T_0^{5/7}\log (T_0/h)^{2/7}}{\sqrt{h}}\right).$$ $\hfill\Box{}$

\section{Regret Lower Bound }\label{sec: lower bound}

In this section, we aim to demonstrate that when the demand function is unknown, the worst-case suboptimal regret for any pricing and capacity-sizing policy is at least of order $\Omega(\sqrt{T})$, where $T$ denotes the total time elapsed. However, deriving a tight lower bound for the regret due to nonstationarity presents significant technical challenges.
First, this type of regret is not necessarily positive in general, making it inherently difficult to analyze. Second, transient errors in the queueing system can become substantial if the control parameters are frequently adjusted, further complicating the task of bounding the nonstationary regret from below for arbitrary learning policies.
Despite these challenges, we provide a partial result in the form of a theoretical lower bound for the suboptimality regret component $R_1$, showing that it scales at the order of $\Omega(\sqrt{T})$.

In particular, we construct a specific demand function with an unknown parameter. The proof is then based on the analysis of KL divergence that measures the uncertainty on this unknown parameter. Intuitively, the proof basically says that, on the one hand, if the uncertainty on the parameter is high, the regret is also high because of the uncertainty (Lemma \ref{lmm: lb uncertainty cost}). On the other hand, it is shown that to reduce uncertainty, a learning cost must be paid (Lemma \ref{lmm: lb learning cost}). As a consequence, there is a lower bound for the regret caused by the uncertainty of the parameter. %This argument is commonly used to prove regret lower bound (see, for example, \cite{Tsybakov2011Introduction} and \cite{Broder2012Dynamic}). 

To make the analysis more intuitive, we consider $T$ as an integer and decompose the total time $T$ into $T$ periods with unit period length. We restrict the policy class so that any admissible policy can only change the price and service capacity at the beginning of each period. This simplification is reasonable because changing policy is usually costly for service providers in reality, and this restriction does not lose generality for our intuition in practice. Note that LiQUAR also belongs to this class. We can formally describe the admissible policies as follows. 
Denote $\omega_0$ as the initial decision $(\mu_0,p_0)$ and $\omega_{t},~t\geq1,$ as the arrivals and corresponding job sizes in $t$-th period and let $\bom_t=(\omega_0,\omega_1,\cdots,\omega_t)$. We denote the corresponding filters as $\{(\Omega_t,\mathcal{F}_t)\}_{t=0}^T$.  An admissible policy is defined by a sequence of decision functions $\pi=\{\pi_1,\cdots,\pi_T\},~\pi_t:\Omega_{t-1}\rightarrow \mathbb{R}_+^2$. We denote these non-anticipating policy class as $\Psi$. 
{
	%We first describe the admissible policies. Let $T$ be total length and we consider the following piecewise constant policy class. Let $U\in\mathbb{U}$ be the initialization random variable defined on a probability space. Let $\pi_1:\mathbb{U}\rightarrow \mathbb{R}_+^3$, then the first round decision is $(p_1,\mu_1,T_1)=\pi_1(U)$. For simplicity, let random variable $\omega_{k-1}$ be the arrival times and corresponding job sizes during cycle $k-1$. Then denote $(\Omega_{k-1},\mathcal{H}_{k-1})$ as the underlying measurable space such that $(U,\omega_1,\cdots,\omega_{k-1})$ is measurable with respect to  $(\Omega_{k-1},\mathcal{H}_{k-1})$. Then let $\pi_k: \mathbb{U}\times\Omega_{k-1}\rightarrow\mathbb{R}_+^3,~k\geq 2$ be the decision functions for $k\geq 2$. The decision functions together $\pi=\{\pi_1,\cdots,\pi_n,\cdots\}$ describe the non-anticipating piecewise constant policy $\pi$ and we denote this policy class by $\Psi$.  
	\begin{thm}\label{thm: lower bound}\textbf{(Theoretic Lower Bound of Regret)}
		There exists a demand function $\lambda(p)$ satisfying Assumption 1 in our main paper and a positive constant $C_2$ such that for any admissible policy $\pi\in\Psi$ and $T\geq2$,
		$$R_1(T)\geq C_2\sqrt{T}.$$
\end{thm}}
%	\begin{proof}{Proof of Thoerem \ref{thm: lower bound}}
	%	TBD
	%\end{proof}}

	%	\red{Put this paragraph after Lemma EC.11} For this demand family, there is an uninformative point $(p_0^*,4)$, at which all demands cross. It's also the optimal point for $z=z_0$. Therefore, to learn the demand, the algorithm needs to step away from the uninformative point, which incurs amounts of cost when $z=z_0$. However, if one algorithm performs well at $z_0$, it seldom learns. So it cannot perform well under other $z$. \red{and we denote $p_0^*=5.5,z_0=1$.}

	Next, we first introduce the demand class and some key properties of problem class $\mathcal{C}$ in Section \ref{subsec: demand class}. Based on these properties, we prove two critical lemmas craving the trade-off between learning cost and uncertainty cost in Section \ref{subsec: lower bound proof sec}. The lower bound is the direct consequence of these two lemmas.
	\subsection{Demand Class and Its Properties} \label{subsec: demand class}
	We consider a parametric problem class $\mathcal{C}$ where the demands are linear functions with slope $z$ as parameter
	\begin{equation}
		\lambda(p;z)=4-z(p-5.5),
		\label{eq: lb dmd}
	\end{equation}
	We set $z\in\mathcal{Z}=[0.95,1.05]$ and $\mathcal{B}=[5,6.5]\times[5.4,6]$ and the queueing system is $M/M/1$. Moreover, we set $h_0=1$ and $c(\mu)=\mu$. In this case, the objective function is 
	$$f(\mu,p;z)=-p\lambda(p;z)+\frac{\lambda(p;z)}{\mu-\lambda(p;z)}+\mu.$$
	Denote optimal decision under demand $\lambda(p;z)$ as
	$$(\mu^*(z),p^*(z))=\arg\min_{(\mu,p)\in\mathcal{B}} f(\mu,p;z).$$
	The corresponding suboptimal regret is
	$$R_1(z,\pi,t)=\EE^{z,\pi}\left[\sum_{k=1}^{t}\left(f(\mu_k,p_k;z)-f(\mu^*(z),p^*(z);z)\right)T_k\right].$$
	In the next lemma, we summarize the key properties of this demand class, which we will use in lower bound analysis.
	\begin{lemma}\label{lmm: lb case properties}
		The problem instance class $\mathcal{C}$ has the following properties:
		\begin{enumerate}
			\item \textbf{Uninformative point.} All demand curves cross an uninformative point, i.e., $\lambda(5.5;z)=4$ for all $z\in\mathcal{Z}$. Moreover, $p^*(1)=5.5.$
			\item \textbf{Strongly convex.} For any $z\in\mathcal{Z}$, the objective function $f(\mu,p;z)$ is strongly convex. As a result, there exists a constant $K_5>0$, such that 
			$$\vert f(\mu^*(z),p^*(z))-f(p,\mu;z)\vert\geq K_5\left((p-p^*(z)^2+(\mu-\mu^*(z))^2\right)$$
			\item \textbf{Uniform stability.} The system is uniformly stable for all problem instances, i.e., 			$$\sup_{p,z}\lambda(p;z)=\lambda(5.4;0.95)<\underline{\mu}.$$
			\item \textbf{Continuity of demand function.} The difference between two demand curves can be represented by difference of $z$ and $z_0$
			$$\vert\lambda(p;z)-\lambda(p;z_0)\vert=\vert(p-5.5)(z-1)\vert.$$
			\item \textbf{{Separability} between optimal solutions.} There exists a constant $K_1$ such that $\vert p^*(z)'\vert\geq K_1$ for all $z\in\mathcal{Z}$. Therefore,
			$$\vert p^*(z)-5.5\vert\geq K_6\vert z-1\vert.$$   %separability is better - Guiyu
		\end{enumerate}
	\end{lemma}
	
	\begin{proof}{Proof of Lemma \ref{lmm: lb case properties}}
		Properties 1, 3 and 4 are obvious by direct calculation. For property 2, notice that the demands $\lambda(p;z)$ are linear functions and by direct calculation, we have the strongly convexity result. For property 5, by the first-order condition {$\nabla f(\mu^*(z),p^*(z))=0$}, the optimal solution is given by
		%		$$
		%		\begin{cases}
			%			\mu^*(z)&=\lambda(p^*(z);z)+\sqrt{\lambda(p^*(z);z)}\\
			%			0&=-z\frac{\mu^*(z)}{(\mu-\lambda(p^*(z);z))^2}-\lambda(p^*(z);z)+p^*(z)z.
			%		\end{cases}	
		%		$$
		$$
		\begin{cases}
			\mu^*(z)&=\lambda(p^*(z);z)+\sqrt{\lambda(p^*(z);z)}\\
			1&=(2p^*(z)-6.5-4z^{-1})^2(4+5.5z-p^*(z)z).
		\end{cases}	
		$$
		To show property 5, we define an auxiliary function $g(p,z)=(2p-6.5+4z^{-1})^2(4+5.5z-pz)$. By direct calculation, there is an $p^*(z)\in[5.4,6]$ satisfying $g(p^*(z),z)=1$. In addition, by direct calculation, in our problem instance,
		\begin{align*}
			\frac{\partial}{\partial p}g(p,z)&=[16+22z-6p+6.5+4z^{-1}](2p-6.5-4z^{-1})>0,\\
			\frac{\partial}{\partial z}g(p,z)&=5.5(2p-6.5-4z^{-1})^2+\frac{8}{z^2}(4+5.5z-p)(2p-6.5-4z^{-1})>0.
		\end{align*}
		Note that 
		$$\frac{d}{d z}g(p^*(z),z)=	\frac{\partial}{\partial p}g(p^*(z),z)p^*(z)'+	\frac{\partial}{\partial z}g(p^*(z),z)=0,$$
		which implies that $p^*(z)'<0$ for all $z\in\mathcal{Z}$. Since $\mathcal{Z}$ is compact, there is a constant $K_6>0$ satisfying the statement in this property. This closes the proof. $\hfill\Box{}$
		%Since $(6.5+4z^{-1})/2<5.4$, there is only one solution in $[5.4,6.5]$. By Cardano's formula $\red{citation?}$ of cubic equation, we have the closed form of this solution. Taking derivative to this solution, we find that $p'(z)<0$ for all $z\in\mathcal{Z}$.  Moreover, $\mathcal{Z}$ is compact, so there is a constant $K_1>0$ satisfying the statement in this property. This closes the proof.
	\end{proof}
	According to Lemma \ref{lmm: lb case properties}, this problem class has an uninformative point at $p=5.5$, where all demands cross. It's also the optimal price for $z=1$. As a consequence,  when $z=1$, the algorithm needs to step away from the uninformative point to learn the demand, which will incur suboptimal cost. On the other hand, if one algorithm performs very well when $z=1$, it seldom learns any information and thus cannot perform well under other $z$. The above observations lead to our proof of the lower bound. 
	\subsection{Proof for the Regret Lower Bound} \label{subsec: lower bound proof sec}
	We denote $p_0^*=5.5$ and $z_0=1$. We shall introduce two lemmas to describe the trade-off between learning cost and the cost of uncertainty. We use Kullback-Leibler divergence to measure the information gain. Let $\PP^{\pi,z}_t$ denote the probability measure of $\bom_t$ under demand $\lambda(p;z)$ with policy $\pi$. We measure the knowledge of demand by 
	$$\mathcal{K}(\PP^{\pi,z_0}_T\|\PP^{\pi,z}_T).$$
	The following lemma craves the learning cost. Denote $\underline{\lambda}\equiv \inf_{z,p}\lambda(p;z)=\lambda(6.5;1.05)=2.95$.
	\begin{lemma}\label{lmm: lb learning cost}
		For any $z\in\mathcal{Z}$, $T>0$ and any piecewise constant policy $\pi\in\Psi$,
		$$\mathcal{K}(\PP^{\pi,z_0}_T\|\PP^{\pi,z}_T)\leq \frac{(z-z_0)^2}{2\underline{\lambda}K_5}R_1(z_0,\pi,T)$$ 
	\end{lemma}
	\begin{proof}{Proof of Lemma \ref{lmm: lb learning cost} }
		We decompose the KL-divergence in $T$ into conditional KL-divergence in each periods. 
		%	Denote $\bom_T=(U,\omega_1,\cdots,\omega_{N_\pi})$ as the path during $[0,T]$. By memoryless property of Poisson process,
		%	\begin{align*}
			%		\mathcal{K}(\PP^{\pi,z_0}_T\|\PP^{\pi,z}_T)&=\int_{\bom_T} \log\left(\frac{d\PP^{\pi,z_0}_T}{d\PP^{\pi,z}_T}(\bom_T)\right)d\PP^{\pi,z_0}_T\\
			%		&=\int_{\bom_T}\log\left(\frac{d\PP^{\pi,z_0}_{T}(\bom_{T}\vert \mathcal{H}_{N_\pi-1})d\PP^{\pi,z_0}_{N_{\pi}-1}(\boldsymbol{\omega}_{N_\pi-1})}{d\PP^{\pi,z}_{T}(\bom_{T}\vert \mathcal{H}_{N_\pi-1})d\PP^{\pi,z}_{N_{\pi}-1}(\boldsymbol{\omega}_{N_\pi-1})} \right)d\PP^{\pi,z_0}_{T}(\bom_{T}\vert \mathcal{H}_{N_\pi-1})d\PP^{\pi,z_0}_{N_{\pi}-1}(\boldsymbol{\omega}_{N_\pi-1})\\
			%		&=\int_{\bom_{N_\pi-1}}\int_{\omega_{N_\pi}}\log \left(\frac{d\PP_{N_\pi}^{\pi,z_0}(\omega_{N_\pi}| \mathcal{H}_{N_\pi-1})}{d\PP_{N_\pi}^{\pi,z}(\omega_{N_\pi}| \mathcal{H}_{N_\pi-1})}\right)d\PP^{\pi,z_0}_{N_\pi}(\omega_{N_\pi}\vert \mathcal{H}_{N_\pi-1}) d\PP^{\pi,z_0}_{N_{\pi}-1}(\bom_{N_\pi-1})\\
			%		&\quad\quad\quad + \int_{\bom_{N_\pi-1}}\log\left(\frac{d\PP_{N_\pi-1}^{\pi,z_0}(\bom_{N_\pi-1})}{d\PP_{N_\pi-1}^{\pi,z}(\bom_{N_\pi-1})}\right)d\PP_{N_\pi-1}^{\pi,z_0}(\bom_{N_\pi-1})
			%	\end{align*}
		By chain rule of KL divergence,
		\begin{align*}
			\KK(\PP_T^{\pi,z_0}\|\PP_T^{\pi,z})&=\sum_{t=1}^{T}\KK(\PP_T^{\pi,z_0}\|\PP_T^{\pi,z}\vert \bom_{t-1})\\
			\KK(\PP_T^{\pi,z_0}\|\PP^{\pi,z}_T|\bom_{t-1})&=\int_{\bom_t}\log\left(\frac{d\PP^{\pi,z_0}_t(\omega_t\vert\bom_{t-1})}{d\PP^{\pi,z}_t(\omega_t\vert\bom_{t-1})}\right)d\PP_t^{\pi,z_0}(\bom_t)
		\end{align*}
		Conditional on $\bom_{t-1}$, the arrivals in cycle $t$ follows Poisson process with rate $\lambda_{t}^z\equiv\lambda(p_{t};z) $ and we denote the density function of individual work load $V$ by $g(\cdot)$. Then, using the conditional density of Poisson arrivals, we have
		%\begin{align*}
		%	&\int_{\omega_{N_\pi}}\log \left(\frac{d\PP_{N_\pi}^{\pi,z_0}(\omega_{N_\pi}| \mathcal{H}_{N_\pi-1})}{d\PP_{N_\pi}^{\pi,z}(\omega_{N_\pi}| \mathcal{H}_{N_\pi-1})}\right)d\PP^{\pi,z_0}_{N_\pi}(\omega_{N_\pi}\vert \mathcal{H}_{N_\pi-1})\\
		%	\stackrel{(a)}{=}&\sum_{i=0}^{\infty}\log\left(\frac{(\lambda_{N_\pi}^{z_0}T_{N_\pi})^ie^{-\lambda_{N_\pi }^{z_0}T_{\pi}}(i!)^{-1} (\frac{1}{T_{N_\pi}})^{i}\prod_{j=0}^{i}g(V_i)}{(\lambda_{N_\pi}^{z}T_{N_\pi})^ie^{-\lambda_{N_\pi }^{z}T_{N_\pi}}(i!)^{-1} (\frac{1}{T_{N_\pi}})^{i}\prod_{j=0}^{i}g(V_i)}\right)\frac{(\lambda_{N_\pi}^{z_0}T_{\pi})^ie^{-\lambda_{N_\pi}^{z_0}T_{N_\pi}}}{i!} \\
		%	=&T_{N_\pi}\left[(\lambda_{N_\pi}^{z}-\lambda_{N_\pi}^{z_0})+\lambda_{N_\pi}^{z_0}\log\left(\frac{\lambda_{N_\pi}^{z_0}}{\lambda_{N_\pi}^{z}}\right)\right]\\
		%	=&T_{N_\pi}\left[(\lambda_{N_\pi}^{z}-\lambda_{N_\pi}^{z_0})-\lambda_{N_\pi}^{z_0}\log\left(1+\frac{\lambda_{N_\pi}^{z}-\lambda_{N_\pi}^{z_0}}{\lambda_{N_\pi}^{z_0}}\right)\right]\\
		%	\stackrel{(b)}{\leq} & T_{N_\pi } \frac{(\lambda_{N_\pi}^z-\lambda_{N_\pi}^{z_0})^2}{2\underline{\lambda}}=\frac{T_{N_\pi}}{2\underline{\lambda}} (p_{N_\pi}-p_0^*)^2 (z-z_0)^2\\
		%	\stackrel{(c)}{\leq }& \frac{T_{N_\pi}}{2\underline{\lambda}K}(z-z_0)^2 (f(\mu_{N_\pi},p_{N_\pi};z_0)-f(u^*(z_0),p^*(z_0);z_0)) 
		%\end{align*}
		\begin{align*}
			&\KK(\PP_t^{\pi,z_0}\|\PP_t^{\pi,z}|\bom_{t-1})\\
			=&\int_{\bom_{t-1}}\int_{\omega_t}\log\left(\frac{d\PP^{\pi,z_0}_t(\omega_t\vert\bom_{t-1})}{d\PP^{\pi,z}_t(\omega_t\vert\bom_{t-1})}\right)d\PP^{\pi,z_0}_t(\omega_t|\bom_{t-1})d\PP^{\pi,z_0}_{t-1}(\bom_{t-1})\\
			=&\int_{\bom_{t-1}}\sum_{k=0}^{\infty}\int_{v_1,\cdots,v_k}\frac{(\lambda_t^{z_0})^k e^{-\lambda_t^{z_0}}}{k!}\log\left( \frac{(\lambda_t^{z_0})^k\exp(-\lambda_t^{z_0})(k!)^{-1} 1^{-k} \prod_{i=1}^{k}g(v_i)}{(\lambda_t^{z})^k\exp(-\lambda_t^{z})(k!)^{-1} 1^{-k} \prod_{i=1}^{k}g(v_i)}  \right)dv_1\cdots dv_kd\PP_{t-1}^{\pi,z_0}(\bom_{t-1})\\
			=&\int_{\bom_{t-1}}(\lambda_t^z-\lambda_t^{z_0})+\lambda_t^{z_0}\log\left(\frac{\lambda_t^{z_0}}{\lambda_t^{z}}\right)d\PP_{t-1}^{\pi,z_0}(\bom_{t-1})\\
			=&\int_{\bom_{t-1}}(\lambda_t^z-\lambda_t^{z_0})-\lambda_t^{z_0}\log\left(1+\frac{\lambda_t^{z}-\lambda_t^{z_0}}{\lambda_t^{z_0}}\right)d\PP_{t-1}^{\pi,z_0}(\bom_{t-1})\\
			\stackrel{(a)}{\leq}&\int_{\bom_{t-1}}\frac{(\lambda_t^z-\lambda_t^{z_0})^2}{2\underline{\lambda}}d\PP^{\pi,z_0}_{t-1}(\bom_{t-1})=\frac{(z-z_0)^2}{2\underline{\lambda}}\int_{\bom_{t-1}}(p_t-p^*_0)^2d\PP_{t-1}^{\pi,z_0}\\
			\stackrel{(b)}{\leq}&\frac{(z-z_0)^2}{2K_5\underline{\lambda}}\EE^{\pi,z_0}\left[f(\mu_{t},p_{t};z_0)-f(u^*(z_0),p^*(z_0);z_0)\right]
		\end{align*}
		
		Here (a) uses the fact that $-\log(1+x)\leq -x+\frac{x^2}{2}$, and (b) uses the strongly convex property (Lemma \ref{lmm: lb case properties}) of our problem case. Therefore, summing up all $t$, we have the result.$\hfill\Box{}$
	\end{proof}
	
	The next lemma describes the cost of uncertainty.
	
	\begin{lemma}\label{lmm: lb uncertainty cost}
		For any integer $T\geq1$, set $z_1=z_0+K_7 T^{-1/4}$ with $K_7$ specified later. Then, for any policy $\pi\in\Psi$, we have 
		$$R_1(z_0,\pi,T)+R_1(z_1,\pi,T)\geq \frac{K_5K_6^2K_7^2}{18}T^{1/2}e^{-\mathcal{K}\left(\PP^{\pi,z_0}_T\|\PP^{\pi,z_1}_T\right)}.$$
	\end{lemma}
	Lemma \ref{lmm: lb uncertainty cost} directly follows lemma 3.4 in \cite{Broder2012Dynamic}, so we omit the proof here.
	With these two lemmas, we now complete the proof of Theorem \ref{thm: lower bound}.
	
	\begin{proof}{Proof of Theorem \ref{thm: lower bound}}
		Let  $z_1={z_0+}K_7T^{-1/4}$ and by Lemma \ref{lmm: lb learning cost}, we have 
		$$R_1(z_0,T,\pi)+R_1(z_1,T,\pi)\geq \frac{2\underline{\lambda}K_5}{K_7^2}\sqrt{T}\mathcal{K}(\PP^{\pi,z_0}_T,\|\PP^{\pi,z_1}_T).$$
		By Lemma \ref{lmm: lb uncertainty cost}, we also have 
		$$R_1(z_0,T,\pi)+R_1(z_1,T,\pi)\geq \frac{K_5K_6^2K_7^2}{18}\sqrt{T}e^{-\mathcal{K}(\PP^{\pi,z_0}_T\|\PP^{\pi,z_1}_T)}.$$
		Therefore, set $C_2=\frac{1}{4}\min\left\{\frac{2\underline{\lambda}K_5}{K_7^2},\frac{K_5K_6^2K_7^2}{18}\right\}$ and we have
		\begin{align*}
			\max_{z\in\{z_0,z_1\}} R_1(1,\pi,T)\geq& \frac{R_1(z_0,\pi,T)+R_1(z_1,\pi,T)}{2}\\
			\geq& \frac{\sqrt{T}}{4}\left(\frac{2\underline{\lambda}K_5}{K_7^2}\mathcal{K}(\PP_T^{\pi,z_0}\|\PP_T^{\pi,z_1})+\frac{K_5K_6^2K_7^2}{18}e^{-\mathcal{K}(\PP_T^{\pi,z_0}\|\PP_T^{\pi,z_1})}\right)\\
			\geq&C_2 \sqrt{T}\left(\mathcal{K}(\PP_T^{\pi,z_0}\|\PP_T^{\pi,z_1})+e^{-\mathcal{K}(\PP_T^{\pi,z_0}\|\PP_T^{\pi,z_1})}\right)\\
			\geq &C_2\sqrt{T}
		\end{align*}
		The last inequality is because $x+e^{-x}\geq 1$ for all $x$. This finishes the proof of the lower bound. $\hfill\Box{}$
	\end{proof}
	
	\section{Examples of the Demand Function}\label{appx: demand function}
	In this part, we verify that the following two inequalities in Condition (a) of Assumption \ref{assmpt: uniform} hold for a variety of commonly-used demand functions.
	\begin{align}\label{eq: lambda prime}
		-\lambda'(p)&> \max\left(\sqrt{\frac{0\vee \left(-\lambda''(p)(\bar{\mu}-\lambda(p)\right))}{2}}~,~ \frac{p\lambda''(p)}{2}\right),\\
		\lambda'(p)&> \max_{\mu\in[\underline{\mu},\bar{\mu}]} \left(2g(\mu)\frac{\lambda''(p)\lambda(p)}{\lambda'(p)}-\frac{4\lambda(p)(\mu-\lambda(p))}{h_0C}\right). \label{eq: Hessian}
	\end{align}
	\vspace{1ex}
	\begin{example}[Linear Demand]
		Consider a linear demand function 
		$$\lambda(p)=a-bp, \quad\text{ with } 0<b<\frac{4\underline{\lambda}(\underline{\mu}-\bar{\lambda})}{h_0C}.$$ 
		Then, inequality \eqref{eq: lambda prime} holds immediately as $\lambda''(p)\equiv 0$. Inequality \eqref{eq: Hessian} is equivalent to
		$$-b>-\frac{4\lambda(p)(\underline{\mu}-\lambda(p))}{h_0C},$$ 
		which also holds as $\lambda(p)(\underline{\mu}-\lambda(p))\geq \underline{\lambda}(\underline{\mu}-\bar{\lambda})$.
	\end{example}
	\vspace{1ex}
	\begin{example}[Quadratic Demand]
		Consider a quadratic demand function 
		$$\lambda(p)=c-ap^2, \quad \text{ with }a,c>0 \text{ and } 0<\frac{\bar{\mu}-c}{3\underline{p}^2}<a<\left(\frac{3(\underline{\mu}-\bar{\lambda})\underline{p}}{h_0C}-\frac{\underline{\mu}}{\underline{\mu}-\bar{\lambda}}\right)\frac{\underline{\lambda}}{\bar{p}^2}.$$
		Inequality \eqref{eq: lambda prime} is equivalent to $3a^2p^2>a(\bar{\mu}-c)$, which holds as $a>\frac{\bar{\mu}-c}{3\underline{p}^2}$. For inequality \eqref{eq: Hessian}, note that $\lambda''=-2a$ and $\lambda'=-2ap$. So, for any $\mu\in[\underline{\mu},\bar{\mu}]$, we have
		\begin{align*}
			&\lambda'(p)-2g(\mu)\frac{\lambda''(p)\lambda(p)}{\lambda'(p)}+\frac{4\lambda(p)(\mu-\lambda(p))}{h_0C}\\
			=&-2ap-2\left(\frac{\mu}{\mu-\lambda}-\frac{(\mu-\lambda)p}{h_0C}\right)\frac{\lambda}{p}+\frac{4\lambda(\mu-\lambda)}{h_0C}\\
			=&2p\left(\frac{\lambda}{p^2}\left(\frac{3(\mu-\lambda)p}{h_0C}-\frac{\mu}{\mu-\lambda}\right)-a\right).
		\end{align*}
		Note that $\frac{3(\mu-\lambda)p}{h_0C}-\frac{\mu}{\mu-\lambda}>\frac{3(\underline{\mu}-\bar{\lambda})\underline{p}}{h_0C}-\frac{\underline{\mu}}{\underline{\mu}-\bar{\lambda}}>0$ by our assumption, and consequently,
		$$\frac{\lambda}{p^2}\left(\frac{3(\mu-\lambda)p}{h_0C}-\frac{\mu}{\mu-\lambda}\right)-a>\left(\frac{3(\underline{\mu}-\bar{\lambda})\underline{p}}{h_0C}-\frac{\underline{\mu}}{\underline{\mu}-\bar{\lambda}}\right)\frac{\underline{\lambda}}{\bar{p}^2}-a>0,$$
		which shows that \eqref{eq: Hessian} holds.
	\end{example} 
	\vspace{1ex}
	\begin{example}[Exponential Demand] Consider an exponential demand function $$\lambda(p)=\exp(a-bp),\quad \text{ with }b>0 \text{ and }
		b\bar{p}<2.$$
		\noindent Then $\lambda'(p)=-b\lambda(p)$ and $\lambda''(p)=b^2\lambda(p)>0$. Therefore, inequality \eqref{eq: lambda prime} is automatically satisfied as $b<2/\bar{p}$.  For inequality \eqref{eq: Hessian}, given that $p\leq\bar{p}<2/b$,  we have, for any $\mu\in[\underline{\mu},\bar{\mu}]$, 
		\begin{align*}
			&\lambda'(p)-2g(\mu)\frac{\lambda''(p)\lambda(p)}{\lambda'(p)}+\frac{4\lambda(p)(\mu-\lambda(p))}{h_0C}\\
			=~&-b\lambda(p)-2\frac{\mu}{\mu-\lambda}\cdot\frac{b^2\lambda^2(p)}{-b\lambda(p)}+\frac{4\lambda(\mu-\lambda)-2bp\lambda(\mu-\lambda)}{h_0C}\\
			>~&-b\lambda(p)+2\frac{\mu}{\mu-\lambda}b\lambda(p)>b\lambda(p)>0.
		\end{align*}
		Therefore, \eqref{eq: Hessian} holds as well.
	\end{example}
	\vspace{1ex}
	\begin{example}[Logit Demand] Consider a logit demand function $$\lambda(p)=c\cdot \exp(a-bp)/(1+\exp(a-bp)),\quad \text{ with }
		a-b\bar{p}<\log(1/2)\text{ and } 0<b<2/\bar{p}.$$
		We have
		$$\lambda'(p) = -\frac{b}{1+e}\lambda(p), ~ \lambda''(p)=\frac{b^2(1-e)}{(1+e)^2}\lambda(p), \text{ with }e \equiv \exp(a-bp).$$
		As a result, inequality \eqref{eq: lambda prime} becomes
		$2>bp(1-e)/(1+e)$ if $e<1$. Since $a-bp<\log(1/2)$, $e<1/2$ and \eqref{eq: lambda prime} holds accordingly. We next show that \eqref{eq: Hessian} holds as well. For any $\mu\in[\underline{\mu},\bar{\mu}]$, 
		\begin{align*}
			&\lambda'(p)-2g(\mu)\frac{\lambda''(p)\lambda(p)}{\lambda'(p)}+\frac{4\lambda(p)(\mu-\lambda(p))}{h_0C}\\
			=~& \left(-\frac{b}{1+e} + \frac{2\mu(1-e)b}{(\mu-\lambda)(1+e)} - \frac{2p(\mu-\lambda)}{h_0C}\cdot\frac{b(1-e)}{1+e} +\frac{4(\mu-\lambda)}{h_0C}\right)\cdot\lambda\\
			>~& \left(-\frac{b}{1+e} + \frac{\mu b}{(\mu-\lambda)(1+e)} - \frac{2bp(1-e)}{1+e}\frac{(\mu-\lambda)}{h_0C} +\frac{4(\mu-\lambda)}{h_0C}\right)\cdot\lambda>0,
		\end{align*}
		where the first inequality holds as $0<e<1/2$ and the second inequality holds as long as $b<2/p$. So \eqref{eq: Hessian} holds as well. In summary, we can conclude that \eqref{eq: lambda prime} and \eqref{eq: Hessian} hold if $0<b<2/\bar{p}$ and $a-b\bar{p}<\log(1/2)$.
	\end{example}
	%\end{APPENDICES}
	\section{Additional Numerical Experiments}\label{appx: additional numeric}
	%In this section, we conduct additional numerical experiments to confirm the practical effectiveness of our algorithm. In what follows, we first test the performance of LiQUAR under exponential demand. Then, we test the robustness of LiQUAR under different model settings.
	\subsection{Robustness of LiQUAR}
	In this section, we give more discussion on the robustness of LiQUAR via numerical examples. Specifically, we test the performance of LiQUAR in a set of model settings with different values of optimal traffic intensity $\rho^*$ and service time distributions.

	We consider an $M/GI/1$ model with phase-type service-time distribution and the logistic demand function in \eqref{eq: logit demand} with $M_0=10, a=4.1$ and $b=1$. We fix staffing cost coefficient $c_0=1$ in \eqref{eq: cost function} in this experiment. By PK formula and PASTA, the service provider's problem reduces to 
	$$\min_{\mu,p}\left\{f(\mu,p)=-p\lambda(p)+\frac{h_0(1+c_s^2)}{2}\cdot\frac{\lambda(p)/\mu}{1-\lambda(p)/\mu}+\mu\right\},$$
	where $c_s^2$ is SCV of the service time. We investigate the impact on performance of LiQUAR of the following two factors: (i) the optimal traffic intensity $\rho^*$ (which measures the level of heavy traffic), and the service-time SCV $c_s^2$ (which quantifies the stochastic variability in service and in the overall system).
	\begin{figure}[h]
		\vspace{-0.2in}
		\includegraphics[width=\linewidth]{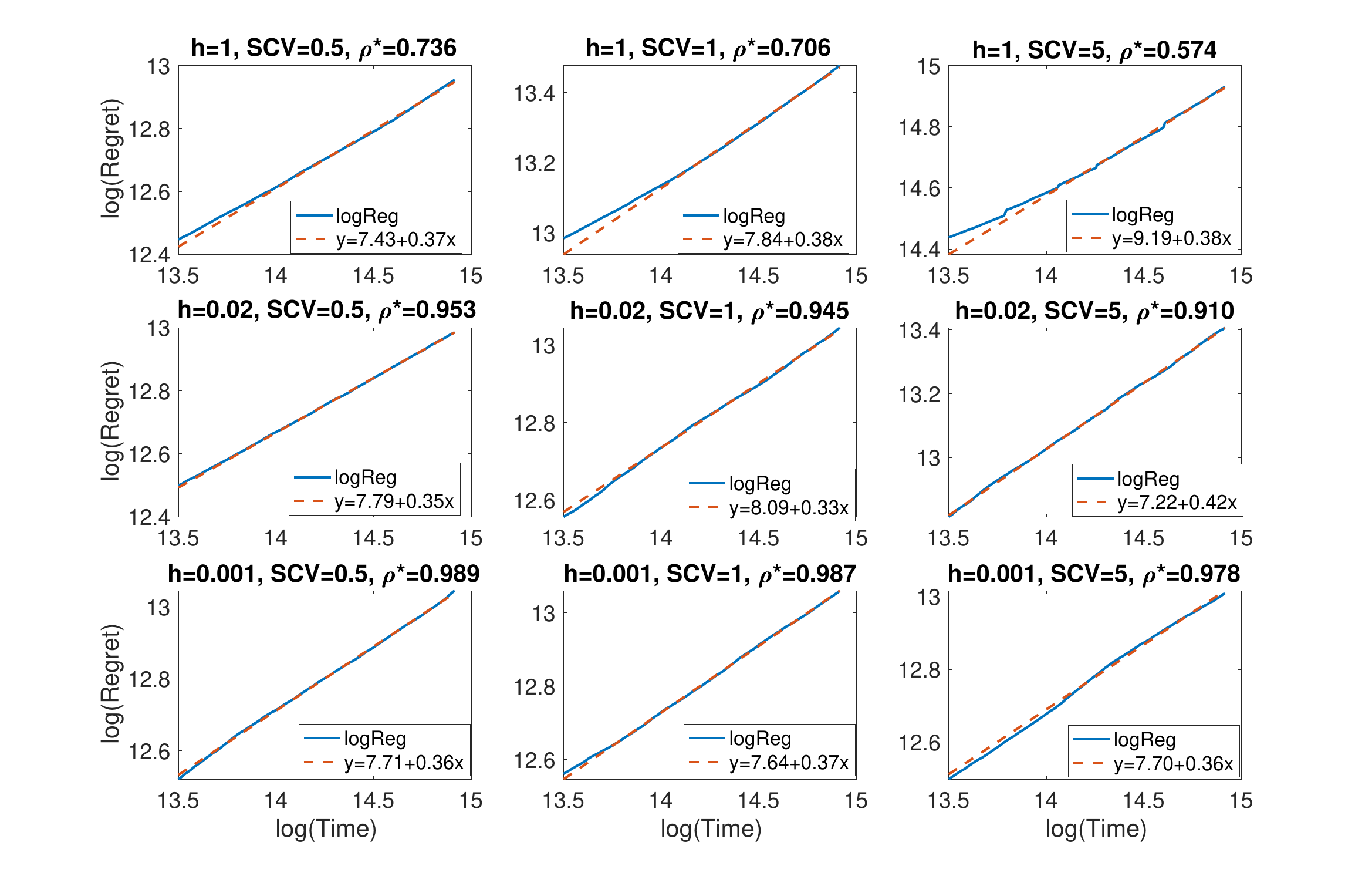}
		\vspace{-0.3in}
		\caption{The regret curve in logarithm scale and a linear fit for the $M/GI/1$ model, under different traffic intensity $\rho^*\in[0.547, 0.989]$ and service-time SCV $c_s^2 =0.5$ ($E_2$ service), $1$ ($M$ service) and $5$ ($H_2$ service). All curves are estimated by averaging 100 independent runs. }
		\label{fig: loglog rubust}
		\vspace{-0.2in}
	\end{figure}
	
	To obtain different values of $\rho^*$, we vary the holding cost $h_0\in\{0.001,0.02,1\}$. For the SCV, we consider $c_s^2=0.5,1,5$  using Erlang-2, exponential and hyperexponential service time distributions. In Figure \ref{fig: loglog rubust} we plot the regret curves in logarithm scale along with their linear fits in all above-mentioned settings. % So there are in total 9 combinations of $(\rho^*,c_S^2)$.
	%In this experiment, we intentionally set the initial decision parameters $(\mu_1,p_1)$ far away from the optimal decision and use the same set of hyperparameters in the  9 different cases: 
	We set $\eta_k=4k^{-1},\delta_k=\min(0.1,0.5k^{-1/3})$, $T_k=200k^{1/3}$ and $\alpha=0.1$. For all 9 cases, we run LiQUAR for $L=1000$ iterations and estimate the regret curve by averaging 100 independent runs. %The algorithm's efficiency is measured by the asymptotic order of regret, which is estimated by a linear fit of log-log curve.
	
	%The log-log regret curves of the 9 cases are reported in Figure \ref{fig: loglog rubust}.
	Note that the optimal traffic intensity $\rho^*$ ranges from $0.547$ to $0.987$. In all the cases, the linear  fitted regret curve has a slope below the theoretic bound $0.5$, ranging in $[0.35,0.42]$. 
	Besides, the intercept (which measures the constant term of the regret) does not increase significantly in $\rho^*$ and ranges in $[7.64, 7.79]$ for $\rho^*>0.95$. The results imply that the performance of LiQUAR is not too sensitive to the traffic intensity $\rho^*$ and service-time SCV. %In addition, as the same set of hyperparameters is used for all cases, the results also indicate that GOLiQ-UDS is easy to tune, which we believe is a valuable feature of GOLiQ-UDS for practical applications. %to some extent, this experiment shows that we do not need to pick the hyperparameters case by case. It is possible to have a hyperparameter with satisfactory performance in different cases. This result may reduce the effort of tuning.
	\subsection{Relaxing the Uniform Stability Condition}\label{sec: LiQUAR-b}
	{In this section, we address the relaxation of the uniform stability assumption by introducing an enhanced version of LiQUAR that adaptively detects and mitigates system instability. Specifically, the refined algorithm incorporates two additional hyperparameters: a workload threshold, $\tau$, and an anchoring price, $p_a$, under which the system is known to remain stable.
		During each cycle, the algorithm continuously monitors the observed workload. If the average workload exceeds the threshold $\tau$, the current cycle is terminated, and a new cycle is initiated. In this new cycle, the service price is updated using a weighted combination of the current price and the anchoring price $p_a$, effectively applying a backtracking mechanism to restore stability. If the average workload remains below $\tau$, the system proceeds identically to the original LiQUAR algorithm.
		For a detailed description of this enhanced method, please refer to Algorithm \ref{alg: backtracking}. We refer to this updated algorithm as LiQUAR with backtracking (LiQUAR-b).}

\begin{figure}
		\centering
		\includegraphics[width=0.8\linewidth]{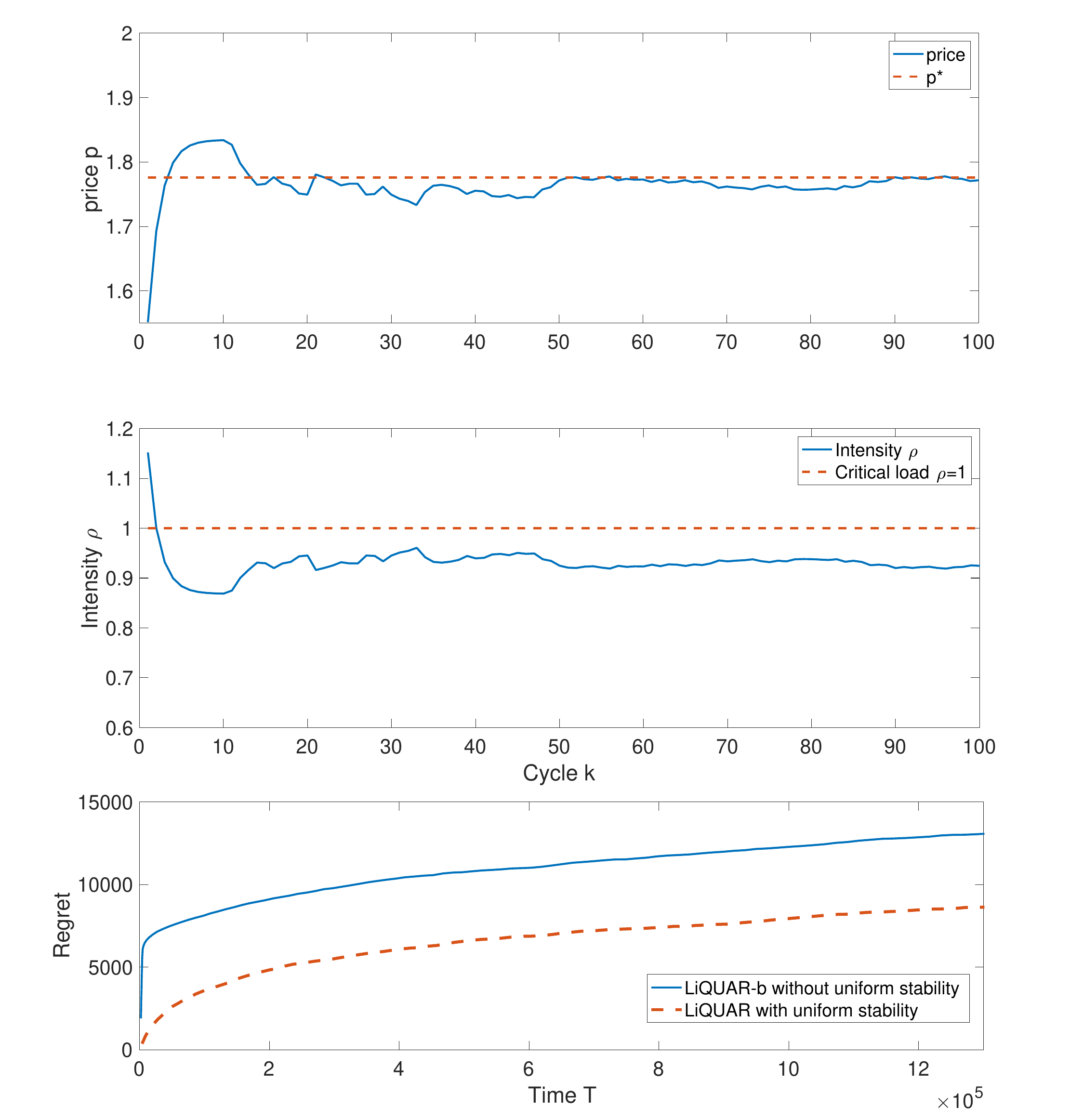}
		\caption{LiQUAR vs. LiGUAR-b:  (i) sample path of price $p_k$ under LiQUAR-b (top panel); (ii) sample path of traffic intensity $\rho_k$ under LiQUAR-b; (iii) regret curves under LiQUAR with $p_k$ subject to uniform stability constraint and  LiQUAR-b with $p_k\in[0,\infty)$. The hyperparameter choices are $T_k=2000k^{1/3},\delta_k=0.07k^{-1/3},\eta_k=0.21k^{-1}$ and $\tau=141$ with the anchoring price $p_a=1.84$ and $h=0.005$. }
		\label{fig: stable}
	\end{figure}

	\begin{algorithm}[h]
		\SetAlgoLined
		\KwIn{number of iterations $L$, workload threshold $\tau$, anchoring price $p_a$\;
			parameters $0<\alpha<1$, and $T_k$, $\eta_k$, $\delta_k$ for $k=1,2,..,L$\;
			initial value $\bar{p}_1$, $W_1(0)=0$\;}
		\For{$k=1,2,...,L$}{
			\vskip 1ex
			\textbf{Set control parameter} $p_{2k-1}= \bar{p}_k-\delta_k /2$ and Stable Sign$=0$\; 
			\While{$\frac{1}{t}\int_{0}^{t}Q_{2k-1}(t)dt<\tau$ for $t<T_k$}{\textbf{Run Cycle $2k-1$:} Run the system under $p_{2k-1}$ \;}
			%\textbf{Run Cycle $2k-1$:} Run the system for $T_k$ units of time under control parameter 
			%Calculate $(\hat{W}_{2k-1}(t))$ according to \red{Procedure ??}\;
			\vskip 1ex
			\If{Cycle $2k-1$ finishes without early-stop}{
				\textbf{Set control parameter} $p_{2k}= \bar{p}_k+\delta_k /2$ \; 
				\While{$\frac{1}{t}\int_{0}^{t}Q_{2k}(t)dt<\tau$ for $t<T_k$}{\textbf{Run Cycle $2k$:} Run the system under $p_{2k}$ \;}
				\If{Cycle $2k$ finishes without early-stop}{Stable Sign $=1$}
			}
			%Calculate $(\hat{W}_{2k}(t))$ according to \red{Procedure ??}\;
			\vskip 1ex
			\If{Stable Sign=1}{
				\textbf{Compute FD gradient estimator:} 
				{\small
					\begin{equation*}
						\begin{aligned}
							H_k=\frac{1}{\delta_k}\Bigg[\frac{h_0}{(1-2\alpha)T_k}&\int_{\alpha T_k}^{(1-\alpha) T_k}\left(Q_{2k}(t)-Q_{2k-1}(t)\right)dt-\frac{p_{2k}N_{2k}-p_{2k-1}N_{2k-1}}{T_k})\Bigg]
						\end{aligned}
					\end{equation*}
				}
				%where the integral term is calculated according to \eqref{eq: integral}\;
				\textbf{Update} $\bar{p}_{k+1}= \Pi_{[0,\infty)} (\bar{p}_k- \eta_k H_k)$.	
				\vskip 1ex
			}
			\Else{\textbf{Backtracking: }$\bar{p}_{k+1}=\frac{\bar{p}_k+p_a}{2}$}
		}
		\caption{LiQUAR with backtracking}
		\label{alg: backtracking}
	\end{algorithm}
	{Next, we test the performance of LiQUAR-b under the heavy-traffic setting in Section \ref{sec: PTO_compare_HT}, this time with the uniformly stable assumption \textbf{relaxed}. Specifically, we consider a pricing problem for $M/M/1$ queue having exponential demand function $\lambda(p)=\exp(a-bp)$ with $a=1+\log(2)$ and $b=1$. In addition, we set coefficient of holding cost $h=0.005$. To see how LiQUAR-b can help maintain a stable system, we relax the feasible domain for $p$ from the uniform stable region to $[0,\infty)$. In addition, we intentionally make the system unstable at $t=0$ by setting the initial price $p_0=1.55$ and thus initial traffic intensity $\rho_0=\lambda(p_0)/\mu=1.15>1$. Following the analysis in Section \ref{sec: heavy traffic}, we set the hyperparameters $T_k=2000k^{1/3},\delta_k=0.07k^{-1/3},\eta_k=0.21k^{-1}$ and $\tau=141$ with the anchoring price $p_a=1.84$. In the top panel of \ref{fig: stable},  we plot the learning curve of price $p_k$ which eventually converges to the optimal $p^*$. The middle panel of \ref{fig: stable} shows that traffic intensity $\rho_k$ is consistently held below 1 after a few iterations.}
	
	To evaluate the impact of this relaxation on regret, we consider two operational scenarios: (i) LiQUAR-b with \( p \in [0, \infty) \), and (ii) LiQUAR with \( p \) constrained to a uniform stable region as used in Section \ref{sec: heavy traffic}. We set the hyperparameters as \( T_k = 2000k^{1/3} \), \( \delta_k = 0.07k^{-1/3} \), and \( \eta_k = 0.21k^{-1} \), and plot the regret curves for LiQUAR-b and LiQUAR in the bottom panel of Figure \ref{fig: stable}. 
	We observe that LiQUAR-b exhibits a steeper initial regret growth, attributable to the algorithm’s need to dedicate initial iterations to steer the policy into the stable domain. This is because \( p_0 \) for LiQUAR-b was deliberately chosen outside the stability domain, while \( p_0 \) for LiQUAR was set within it. 
	Notably, the backtracking mechanism in LiQUAR-b acts as an effective safeguard, ensuring that subsequent \( p_k \) values remain within the stability region. Consequently, after this initial adjustment phase, the regret curve for LiQUAR-b stabilizes and exhibits a growth rate similar to that of LiQUAR.

	{Although this experiment serves only as an initial exploration of relaxing the uniform stability condition in LiQUAR, Figure \ref{fig: stable} demonstrates the promising potential of the LiQUAR-b approach. Several important directions for future research remain:
		\begin{itemize}
			\item
			Theoretical regret analysis: A comprehensive theoretical analysis of the regret for LiQUAR-b is required. This involves developing effective techniques to bound the regret growth during unstable learning cycles and quantifying the impact of accumulated excessive workload on subsequent cycles.
			\item
			Optimizing the detection threshold $\tau$: It is crucial to determine an optimal workload threshold $\tau$ that balances two competing factors: frequent false detections of instability (when $\tau$ is too small) and excessive workload accumulation during unstable cycles (when $\tau$ is too large).
		\end{itemize}
		Further exploration of these aspects will provide deeper insights and strengthen the practical applicability of LiQUAR-b. We leave these extensions for future study.

	}
	
	\section{Details of Algorithm \ref{alg: PG} in Section \ref{sec: RL}}\label{sec: AlgRL}
	In this section, we provide the detailed description for the PG algorithms in Algorithm \ref{alg: PG}, the outline of which is described in Section \ref{sec: RL}. Specifically, in Algorithm \ref{alg: PG}, Policy Gradient algorithm organizes time by cycles, with each cycle containing $L$ episodes. In each episode, the system operates the system following $\pi_\theta$ for episode length $T$ time units. At the end of each episode, a gradient gradient estimator in this episode $\hat{\nabla}_{i,t}$ is calculated using the policy gradient formula  \cite[p.339]{SB18} and the closed form of Gaussian parameterization (line 9 in Algorithm \ref{alg: PG}). Then, at the end of each cycle, an overall policy gradient estimator is obtained by averaging over all the episodic policy gradient estimators in the cycle (line 11 in Algorithm \ref{alg: PG}). The full algorithm is given in Algorithm \ref{alg: PG}.
	
	\begin{algorithm}[h]
		\SetAlgoLined
		\KwIn{normal parameterization $\pi(a|\theta)$, step size $\eta>0$, initial policy parameter $\theta:(\bar{p}_1,\bar{\mu}_1,\sigma_{p,1}^2,\sigma_{\mu,1}^2)$, cycle length $L$ (how many episodes in each episode), episode length $T$ (how many time slots in each episode)\;}
		\For{each cycle}{
			\For{episode $i=1:L$}
			{
				Generate an episode $Q_1,(p_1,\mu_1),R_1,\cdots,Q_{T-1},(p_{T},\mu_{T}),R_T$ following $\pi_\theta$\;
				$\bar{R}=\frac{1}{T}\sum_{t=1}^{T}R_T$\;
				\For{$t=1,\cdots,T$}{
					$G=\sum_{k=t}^{T}R_k-\bar{R}$\;
					$\hat{\nabla}_{i,t}\leftarrow G\cdot  \begin{pmatrix}
						(p_t-\bar{p})/\sigma_p^2\\
						(\mu_t-\bar{\mu})/\sigma_\mu^2\\
						\left[(p_t-\bar{p})^2-\sigma_p^2\right]/\sigma_p^3\\
						\left[(\mu_t-\bar{p})^2-\sigma_\mu^2\right]/\sigma_\mu^3
					\end{pmatrix}$\;
				}
				$\hat{\nabla}_i=\frac{1}{T}\sum_{t=1}^{T}\hat{\nabla}_{i,t}$\;
			}
			$\theta\leftarrow \theta +\eta \cdot \frac{1}{L}\sum_{i=1}^{L}\hat{\nabla}_i$\;
		}
		\caption{Policy Gradient method}
		\label{alg: PG}
	\end{algorithm}

	%%%%%%%%%%%%%%%%%%%%%%%%%%%%%%%%%%%%%%%
	\newpage
	%%%%%%%%%%%%%%%%%%Glossary of notations%%%%%%%%%%%%%%%%%%%%%%%%%%%%				
	\begin{table}
		\centering
		\scriptsize{
			\begin{tabular}{c|ll}
				\hline
				\multicolumn{1}{l}{} & Notation                                     & Description                                                         \\
				\hline
				\multirow{13}{*}{\begin{tabular}[c]{@{}c@{}}Model\\ parameters\\ and\\ functions\end{tabular}} &
				$\mathcal{B}=[\underline{\mu},\bar{\mu}]\times[\underline{p},\bar{p}]$ &
				Feasible region \\
				& $c(\mu)$                                     & Staffing cost                                                \\
				& $c_s^2=Var(S)/\EE[S]^2$                      & Squared coefficient of variation (SCV) of the service times         \\
				& $C=\frac{1+c_s^2}{2}$                        & Variational constant in PK formula                                  \\
				& $f(\mu,p)$                                   & Objective (loss) function                                           \\
				& $h_0$                                        & Holding cost of workload                                            \\
				& $\lambda(p)$                                 & Underlying demand function                                          \\
				& $\mu$                                        & Service rate                                                        \\
				& $p$                                          & Service fee                                                         \\
				&
				$\theta,\gamma_0,\eta$ &
				Parameters of light-tail assumptions (Assumption \ref{assmpt: light tail}) \\
				& $V_n$                                        & Individual workload                                                 \\
				& $W_\infty(\mu,p)$                            & Stationary workload under decision $(\mu,p)$                        \\
				& $\bx^*=(\mu^*,p^*)$                             & Optimal decision service rate and fee                               \\
				\hline
				\multirow{10}{*}{\begin{tabular}[c]{@{}c@{}}Algorithmic\\ parameters\\ and\\ variables\end{tabular}} &
				$\alpha$ &
				Warm-up and overtime rate \\
				& $\delta_k, (\delta_k^h)$                                   & Exploration length in iteration $k$ (of $h^{\text{th}}$ system )                                 \\
				& $\eta_k, (\eta_k^h)$                                     & Step length for gradient update in iteration $k$ (of $h^{\text{th}}$ system )                    \\
				& $\bH_k$                                        & Gradient estimator in iteration $k$                                 \\
				& $\hat{f}^G(\mu_l,p_l)$                       & Estimation of objective function in cycle $l$                       \\
				&$Q_k^h(t)$										& Queue length at time $t$ in cycle $k$ of the $h^{\text{th}}$ system\\
				& $T_k,T_{k(l)}, (T_k^h)$                               & Cycle length of iteration $k$ and cycle $l$  (of $h^{\text{th}}$ system )                         \\
				& $W_l(t) (\hat{W}_l(t))$                      & (Estimated) workload at time $t$ in cycle $l$                       \\
				& $X_l(t)$                                     & Observed busy time at time $t$ in cycle $l$                         \\
				& $\bar{\bx}_k$                                  & Control parameter in iteration $k$                                 \\
				& $\bZ_k$											&Updating direction in iteration $k$									\\
				\hline
				\multirow{17}{*}{\begin{tabular}[c]{@{}c@{}}Constants and\\bounds in\\regret analysis\end{tabular}} &
				$B_k,\mathcal{V}_k$ &
				Bias and Variance upper bound for $H_k$ \\
				& $c$                                          & Constant for noise-free FD error in Lemma \ref{lmm: fd approximation}          \\
				& $c_\eta,c_T,c_\delta$                        & Coefficient of hyperparameters in Theorem \ref{thm: upper bound}                         \\
				& $C$								& Constant in Theorem \ref{thm: heavy traffic} irrelevant to $h$ \\
				& $C_0$                                        & Constant in Lemma \ref{lmm: regret of nonstationarity in warm-up period}                                                 \\
				& $M$                                          & Upper bound for queueing functions in Lemma \ref{lmm: uniform bound}                       \\
				& $\gamma$                                     & Ergodicity rate constant in Lemma \ref{lmm: ergodicity}                                 \\
				& $K_0,K_1$                                    & Convex and smoothness constant of objective function in Lemma \ref{lmm: convexity}     \\
				& $K_2,K_3$                                    & Constants in the proof of Theorem \ref{thm: convergence} in Appendix \ref{appx: proof of convergence theorem}                                    \\
				& $K_4$                                        & Constant in Lemma \ref{lmm: exploration cost}                                                \\
				& $K_5,K_6,K_7$							& Constants in Theorem \ref{thm: lower bound} in Section \ref{sec: lower bound}\\
				& $K_V$                                        & Constant of auto-correlation in Lemma \ref{lmm: autocovariance W}                             \\
				& $K_M$											&Constant of MSE of $\hat{f}^G$ in Proposition \ref{prop: estimation error}												\\
				& $R(L),R_1(L),R_2(L),R_3(L)$	& Total regret, regret of sub-optimality, non-stationarity, finite difference\\
				& $\theta_0$                                   & Constant in Lemma \ref{lmm: uniform bound}                                                 \\
				%& $\theta_1=\frac{\theta_0\underline{\mu}}{2}$ & Constant in Proposition \ref{prop: bound Bk V_k}                                           \\
				& $\theta_1=\min(\gamma,\theta_0\underline{\mu}/2)$             & Constant in Proposition \ref{prop: bound Bk V_k}                                           \\
				& $\bar{W}_l(t)$                               & Stationary workload process coupled from the beginning of cycle $l$ \\
				& $\bar{W}_l^s(t)$                             & Stationary workload process coupled from time $s$ of cycle $l$ (in Appendix)      \\
				& $W^D_l(t),X^D_l(t)$                          & Workload and observed busy time for the dominating queue (in Appendix)          
			\end{tabular}
		}
		\caption{Glossary of key notations}
		\label{tab: notations}
	\end{table}

\end{document}